\documentclass[preprint,10pt]{elsarticle}

\usepackage[english]{babel}
\usepackage{latexsym,amsmath,amssymb,amsbsy,amstext,amscd,amsfonts,amsthm}
\usepackage{graphics,graphicx}
\usepackage{color}
\usepackage{tabularx}
\usepackage{multirow}
\usepackage{booktabs}
\usepackage{float}
\theoremstyle{definition}
\newtheorem{remark}{Remark}
\usepackage{array}

\usepackage{amssymb}

\begin{document}

\begin{frontmatter}



\title{Mathematical modelling of Solid Oxide Fuel Cells revisited - a modified formulation of the problem}


\author{Michal Wrobel}
\ead{wrobel.mc@gmail.com}
\author{ Grzegorz Brus}

\address{Department of Fundamental Research in Energy Engineering,
\\AGH University of Science and Technology
\\ Mickiewicza 30, 30-059 Krakow, Poland}

\begin{abstract}
In the  paper a mathematical model of the PEN  structure (positive electrode - electrolyte -negative electrode) of the Solid Oxide Fuel Cell (SOFC)  is analyzed. It is proved that classical formulation of the problem leads inevitably to locally unphysical effects related to negative values of the activation overpotential. Moreover, the active layers' thicknesses are shown to be components of solution and cannot be predefined in an arbitrary way. A modified mathematical formulation of the problem is proposed which includes this novelty alongside an amended definition of the activation overpotential. A dedicated computational scheme is developed for the cathode sub-problem. The accuracy of computations is investigated by means of a newly introduced analytical benchmark example. The numerical results obtained for LSCF cathode are  used to discuss certain aspects of the modified formulation and the active layer thickness. The new modelling approach is validated by comparison with experimental data.

\end{abstract}


\end{frontmatter}


\section{Introduction}

Solid oxide fuel cells (SOFCs) are devices that convert the chemical energy of fuel directly to electricity. As such, they do not comply with the limitations of the Carnot cycle concerning the maximal efficiency achievable by a classical thermodynamic engine. Recently, SOFCs are becoming increasingly more popular due to their: i) high efficiency, ii) capability to work with a variety of fuels (including hydrogen, carbon monoxide, methane, methanol and ammonia \cite{Zhu_2005,Mozdzierz_2018,Strazza_2010,Tan_2018}), iii) low emission of noxious substances (when operating on pure hydrogen they emit no pollutants at all), iv) multiple geometrical configurations, and v) relatively low cost. On the other hand, high operating temperature results in long start-up time, durability issues related to the thermal stresses and limited selection of constructing materials.

Although a number of experimental techniques can be used to analyze and understand the underlying physics of SOFCs operation,  mathematical modelling still remains a powerful tool for non-invasive investigation of certain aspects of the devices' functioning.
However, numerical (and, generally speaking, mathematical) simulation of SOFC constitutes a formidable task. The main difficulties stem from: i) pronounced multiscale effects, ii) multiphysical character of the problem and complex interactions between component physical fields, iii) highly non-linear nature of the underlying phenomena, iv) uncertainties about description of electrochemical reactions kinetics,  physical and geometrical properties (microstructure) of the materials used in SOFCs, and others.

The mathematical model of SOFC should describe properly the component physical fields and their interactions at different temporal and spatial scales. Respective physical phenomena manifest themselves with different intensities depending on the operating regimes and conditions, and thus each of them can affect the overall cell performance in a variety of ways. For this reason, there exists neither a universal model of SOFC nor a versatile computational method to solve the system of governing equations. A comprehensive review of the modelling techniques and their classification can be found in \cite{Andersson_2010,Kakac_2007}.

When considering the SOFC problem from the scale point of view one can distinguish the following three approaches:
\begin{itemize}
\item{Microscale modelling, which corresponds to the atomic or molecular level. It can be used for example to investigate the microscopic aspects of oxygen diffusion \cite{Frayret_2005} or ionic transport \cite{Cheng_2005}. }
\item{Mesoscale modelling that involves a larger spatial scale than that of the micromodels, nevertheless capable of describing the actual morphology of microstructure \cite{Prokop_2018}. }
\item{Macroscale modelling where the continuum approach is employed. The macroscale models are used to simulate respective electrodes \cite{Miyoshi_2016,Enrico_2014}, the PEN structure (positive electrode - electrolyte- negative electrode) \cite{Brus_et_al_2017}, the whole cell \cite{Ho_2010} or a stack of cells \cite{He_et_al_2017,Al-Masri_2014,Calise_2009} in terms of transport mechanisms and structural mechanics. Depending on the complexity of performed analysis one can encounter here various geometrical models (ranging from 0D to 3D), various descriptions of the  diffusion and fluid flow phenomena, and various mechanisms of interaction between respective physical fields. The macroscale approach has been very useful in understanding the complex phenomena that govern the operation of the cell and in this way improve and optimize the devices' designs. However, accounting for particular microstructure and related properties is not so straightforward here. A popular method to do this is to employ effective coefficients  in the transport equations which convey information from the micro-level \cite{Kishimoto_2011,Brus_et_al_2017,Kishimoto_2013}. This information can be retrieved from the real structural analysis. Recently, the most popular technique for 3D identification and quantification of porous microstructure has been combination of focused ion beam and scanning electron microscope (FIB-SEM) \cite{Joos_2011,Iwai_2010}.
Some macroscale models are calibrated by means of experimental data. Namely, certain values of microstructural and/or reaction kinetic parameters are taken in a way to match the fuel cell polarization curves \cite{Lee_2009,Khaleel_2003}. In this manner the deficiencies inherent to the model can be alleviated. }
\end{itemize}

The above groups of methods can be integrated in the framework of the multiscale approach. The overview of respective strategies can be found in \cite{Andersson_2010}.

The electrochemical reactions take place at the so-called three-phase boundaries (TPBs) - the locations in the electrode where the electron conducting, ion conducting and gaseous phases adjoin. In most studies a model of a single global charge transfer reaction is assumed \cite{Andersson_2010}. It leads to the Butler-Volmer equation which constitutes a relation between the charge transfer current and the local activation overpotential. It is commonly accepted that only a part of each electrode is electrochemically active  with reactions occurring in the immediate vicinity of the electrode-electrolyte interface \cite{Miyawaki_2014,Nam_2017,Kishimoto_2013}. As the thickness of the active (catalyst) layer is of prime importance for the design of electrodes, there have been many attempts to establish it by means of experiment and simulation \cite{Chen_2008,Kong_2007,Primdahl_1997, Miyawaki_2014,Kishimoto_2013,Nam_2017}. According to \cite{Kishimoto_2013}, evolution of the active zone extent is governed mostly by interplay between the ionic conductivity and the electrochemical reaction rate. Magnification of the former mechanism facilitates migration of the ions further away from the electrolyte, which contributes to the extension of the layer. On the other hand, intensification of the reaction rate increases the ions consumption and in this way reduces the catalyst zone thickness. Thus, the active layer dimension is a function of the microstructural properties of the electrode and the operating conditions (cell current, temperature, concentration of species). As such, it changes with stages and regimes of cell operation. Unfortunately, values of the active zone thickness available in the literature are quite different from each other - see Table I in \cite{Kishimoto_2013}.

In the continuous models of SOFC there are various strategies to introduce the electrochemical effects. A simplified approach assumes that electrochemical activity is represented by a pertinent interface transmission condition along the contact surface of the electrode and the electrolyte \cite{Haberman_2004,Lehnert_2000}. However, in such a case an essential feature of the original physical problem is neglected. Another method is to introduce respective source terms in the governing equations. Here, two techniques are used. In the first one, the expression for the charge transfer current is defined over the entire thickness of the electrode \cite{Brus_et_al_2017,Miyoshi_2016}. Note that this is at odds with experimental observations that the reaction zone is  limited to the immediate neighborhood of the electrolyte. In the second technique a predefined extent of the catalyst layer is assumed \cite{He_et_al_2017,Thinh_2009}. Although this approach yields the most precise description of the actual physics of the problem, it is cumbersome in practical implementation. It stems from the fact that one finds it very hard to determine experimentally  the real thickness of the active layer. Moreover, the thickness evolves with differing conditions and regimes of SOFC operation.  It should be accounted for in the numerical computations, which constitutes a real challenge, especially in the transient analysis. Interestingly, even though the aforementioned strategies to introduce the charge transfer current are essentially different from each other and have to affect the mathematical structure of solution in very different ways, no research on the latter can be found in the literature. Furthermore, according to the author's knowledge, no attempt to investigate the existence and uniqueness of solution for the continuous macroscale SOFC models has been made so far (the existence of the weak bounded solution for a mesoscale isothermal formulation has been proved recently in \cite{Al-arydah_2018}).

The basic motivation for this paper was to bridge the gap in the existing literature and provide a detailed analysis of the mathematical structure of solution for the continuous SOFC model. An isothermal 1D model of the PEN structure is considered. In  section \ref{classic_form} we introduce a classical formulation of the problem with the Fickian diffusion  describing the mass transport. All  boundary and transmission conditions are detailed and a desired qualitative behaviour of the dependent variables is presented. In section \ref{mod_form}, when analyzing the standard model, we conclude that it produces locally unphysical effects related to the negative activation overpotential in the immediate vicinity of the external boundaries of the active layers. Remediation in a form of an amended definition of the activation overpotential (or alternatively an amended definition of the charge transfer current) is proposed. Next, we prove analytically that for the predefined conditions of cell operation (described by the cell current, boundary and transmission conditions) there is only one thickness of the catalyst layer for each electrode for which the solution exists. Thus, the active zones thicknesses have to be elements of solution. By accounting for this novelty alongside the new definition of activation overpotential we introduce a modified formulation of the SOFC problem. In section \ref{num_res} we propose a dedicated computational scheme for the cathode sub-problem, based on a system of integral equations (some elements of this approach were originally developed for the heat and mass transfer problems in \cite{Wrobel_2008,Wrobel_2009,Wrobel_2009a}).  The achievable  accuracy of computations is verified against analytical benchmark example given in \ref{benchmark_app}. The numerical solution obtained for LSCF cathode is used to discuss the problem of the active layer thickness. Finally, the validity of the modified SOFC formulation is checked by comparison with experimental data. Discussion and final conclusions are presented in section \ref{disc}.

\section{Mathematical formulation of the problem}
\label{classic_form}
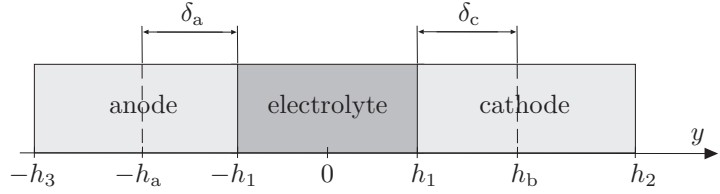
\begin{figure}
\begin{center}
\includegraphics[scale=0.45]{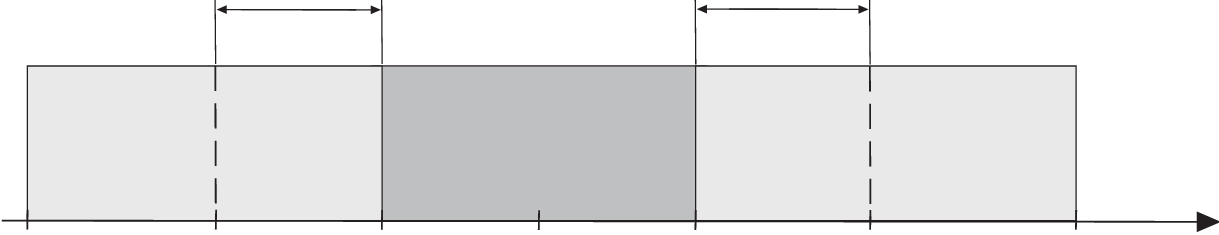}
\put(-10,8){$y$}
\put(-150,-8){$0$}
\put(-116,-8){$h_1$}
\put(-78,-8){$h_\text{b}$}
\put(-33,-8){$h_2$}
\put(-192,-8){$-h_1$}
\put(-228,-8){$-h_\text{a}$}
\put(-268,-8){$-h_3$}
\put(-98,53){$\delta_\text{c}$}
\put(-203,53){$\delta_\text{a}$}
\put(-170,18){$\text{electrolyte}$}
\put(-90,18){$\text{cathode}$}
\put(-230,18){$\text{anode}$}
\caption{The schematic view of the PEN structure.}
\label{PEN_str}
\end{center}
\end{figure}

In our analysis we will use an isothermal 1D model of the PEN structure. Respective subdomains are defined as (see Fig.\ref{PEN_str}): i) anode $-h_3\geq y \geq -h_1$, ii) electrolyte $-h_1\geq y \geq h_1$, iii) cathode $h_1\geq y \geq h_2$. In the anode and the cathode two parts are distinguished: i) the backing layer (anode: $-h_3\geq y \geq -h_\text{a}$, cathode: $h_\text{b}\geq y \geq h_2$), ii) the catalyst (active) layer where the electrochemical reactions take place (anode: $-h_\text{a}\geq y \geq -h_1$, cathode: $h_1\geq y \geq h_\text{b}$). The thicknesses of respective electrodes are: i) anode $h_3-h_1$, ii) electrolyte $2h_1$, iii) cathode $h_2-h_1$. The active layers' spans are denoted as: i) $\delta_\text{a}$ - for the anode, ii) $\delta_\text{c}$ - for the cathode.

We start our analysis by recalling one of the standard systems of governing equations. It describes the charge transfer and the mass transfer (by means of the Fickian diffusion). Respective boundary and transmission conditions are chosen in a way to preserve the basic physical features of the problem. The following mathematical formulation of the problem can be found in many publications (see e.g. \cite{Cannarozzo_2007,He_et_al_2017,Brus_et_al_2017}).

\subsection{The charge transfer equations}

The ions and electrons transport through the cell's structure is governed by the following equations:
\begin{equation}
\label{phi_el_ODE}
\frac{\text{d}}{\text{d}y}\left(\sigma_{\text{el}}\frac{\text{d} \phi_{\text{el}}}{\text{d}y}\right) =
  \begin{cases}
		0       & \quad \text{for } \quad y\in(-h_3,-h_\text{a})\cup(-h_1,h_1)\cup(h_\text{b},h_2),\\
    i^\text{ct}       & \quad \text{for} \quad y\in(-h_\text{a},-h_1),\\
		
    -i^\text{ct}  & \quad \text{for } \quad y\in(h_1,h_\text{b}),
  \end{cases}
\end{equation}
\begin{equation}
\label{phi_ion_ODE}
\frac{\text{d}}{\text{d}y}\left(\sigma_{\text{ion}}\frac{\text{d} \phi_{\text{ion}}}{\text{d}y}\right) =
  \begin{cases}
	0       & \quad \text{for } \quad y\in(-h_3,-h_\text{a})\cup(-h_1,h_1)\cup(h_\text{b},h_2),\\
    -i^\text{ct}       & \quad \text{for} \quad y\in(-h_\text{a},-h_1),\\
    i^\text{ct}  & \quad \text{for } \quad y\in(h_1,h_\text{b}).
  \end{cases}
\end{equation}

In the above formulae $\phi_\text{el}$ [V], $\phi_\text{ion}$ [V] denote the ionic and electric potentials, respectively, while $\sigma_\text{el}$ [S m$^{-1}$], $\sigma_\text{ion}$ [S m$^{-1}$] stand for the effective electron and ion conductivities. $i^\text{ct}$ [A m$^{-3}$] describes the charge transfer current function, to be defined later on. Note that the definition of source terms in equations \eqref{phi_el_ODE}--\eqref{phi_ion_ODE} accounts for the existence of the catalyst (electrochemically active) layers in respective electrodes.

The electron and ion currents' densities (both in [A m$^{-2}$]) are computed as:
\begin{equation}
\label{j_el_def}
j_\text{el}=-\sigma_\text{el} \frac{\text{d}\phi_\text{el}}{\text{d}y}, \quad j_\text{ion}=-\sigma_\text{ion} \frac{\text{d}\phi_\text{ion}}{\text{d}y},
\end{equation}
while the overall cell current density, $j_\text{cell}$ [A m$^{-2}$], yields:
\begin{equation}
\label{j_cell}
j_\text{cell}=j_\text{el}+j_\text{ion}.
\end{equation}

Moreover, from the global balance of electric charge it follows that:
\begin{equation}
\label{j_bal_int}
\int_{-h_\text{a}}^{-h_1} i^\text{ct}dy=\int_{h_1}^{h_\text{b}} i^\text{ct}dy=j_\text{cell}.
\end{equation}

\subsection{The mass transfer equations}

The mass transfer equations describe the diffusion of respective species in the electrodes. In this study we assume that the underlying physical process proceeds according to the Fick's law of diffusion. Being fully aware of the deficiencies of this model \cite{Bertei_2015} we shall not discuss here its applicability. The methodology introduced in the paper can be easily extended to most sophisticated variants of the problem. As such, the simplified diffusion model does not detract from the validity of the general conclusions.

In the anode the diffusion in the fuel mixture is assumed to take place between the hydrogen and the water. Thus, the nitrogen concentration depends on the computed content of the remaining species. Respective equations read:

\begin{equation}
\label{ch2_ODE}
\frac{\text{d}}{\text{d}y}\left(\rho_{\text{f}}D_1\frac{\text{d} C_{\text{H}_2}}{\text{d}y}\right) =
  \begin{cases}
	0       & \quad \text{for } \quad y\in(-h_3,-h_\text{a}),\\
    \frac{M_{\text{H}_2}}{2F}i^\text{ct}  & \quad \text{for } \quad y\in(-h_\text{a},-h_1),
  \end{cases}
\end{equation}

\begin{equation}
\label{ch2o_ODE}
\frac{\text{d}}{\text{d}y}\left(\rho_{\text{f}}D_1\frac{\text{d} C_{\text{H}_2\text{O}}}{\text{d}y}\right) =
  \begin{cases}
	0       & \quad \text{for } \quad y\in(-h_3,-h_\text{a}),\\
    -\frac{M_{\text{H}_2\text{O}}}{2F}i^\text{ct}  & \quad \text{for } \quad y\in(-h_\text{a},-h_1),
  \end{cases}
\end{equation}
where pertinent symbols refer to: $\rho_\text{f}$ [kg m$^{-3}$] - fuel density, $D_1$ [m$^2$ s$^{-1}$] - effective diffusion coefficient of the binary system $\text{O}_2$ - $\text{H}_2\text{O}$, $C_j$ - mass fraction of the $j$-th species, $M_j$ [kg mol$^{-1}$] - molar mass of the $j$-th species, $F$ [C mol$^{-1}$] - the Faraday constant.

From the overall balance of species one has:
\begin{equation}
\label{anode_C_bal}
C_{\text{H}_2}+C_{\text{H}_2\text{O}}+C_{\text{N}_2}=1.
\end{equation}

In the cathode the binary diffusion between the oxygen and the nitrogen is governed by the following ODE:
\begin{equation}
\label{co2_ODE}
\frac{\text{d}}{\text{d}y}\left(\rho_{\text{a}}D_2\frac{\text{d} C_{\text{O}_2}}{\text{d}y}\right) =
  \begin{cases}
	    \frac{M_{\text{O}_2}}{4F}i^\text{ct}  & \quad \text{for } \quad y\in(h_1,h_\text{b}),\\
			0       & \quad \text{for } \quad y\in(h_\text{b},h_2),
  \end{cases}
\end{equation}
where the meaning of respective symbols is analogous to that from \eqref{ch2_ODE}-\eqref{ch2o_ODE} (with $\rho_\text{a}$ [kg m$^{-3}$] being the air density). As:
\begin{equation}
\label{cathode_C_bal}
C_{\text{O}_2}+C_{\text{N}_2}=1,
\end{equation}
we do not write the corresponding ODE for the nitrogen mass fraction. When having the solution of \eqref{co2_ODE} one can easily compute $C_{\text{N}_2}$ from \eqref{cathode_C_bal}.

Conventionally, the flux of $j$-th species is defined as:
\begin{equation}
\label{C_flux_def}
J_j=-\rho D \frac{\text{d}C_j}{\text{d}y},
\end{equation}
where the density and the diffusion coefficient are to be taken accordingly.

\subsection{The boundary and transmission conditions}
\label{BCs_sect}

The credible and efficient simulation of the problem needs careful selection and identification of the boundary and transmission conditions. They should result form the physical specificity of the underlying phenomena and comply with the proper mathematical structure of solution.

Considering the charge transfer, the following boundary conditions hold:
\begin{equation}
\label{j_el_BC}
j_\text{el}(-h_3)=j_\text{cell}, \quad j_\text{el}(h_2)=j_\text{cell},
\end{equation}
\begin{equation}
\label{j_ion_BC}
j_\text{ion}(-h_3)=0, \quad j_\text{ion}(h_2)=0,
\end{equation}
\begin{equation}
\label{fi_el_BC}
\phi_\text{el}(-h_3)=V_3, \quad \phi_\text{el}(h_2)=V_2.
\end{equation}

Note that conditions \eqref{j_el_BC} - \eqref{j_ion_BC} result from the fact  that it is only the electron current that flows in the external circuit of the cell. Moreover, when implementing numerical computations, some of the above conditions constitute components of the solution. Generally, two different strategies are available. When one assumes a predefined value of the voltage drop ($V_3-V_2$) \cite{He_et_al_2017,Thinh_2009}, then the cell current density, $j_\text{cell}$ is to be found computationally. Conversely, when the current density is known \cite{Chan_2001,Matsuzaki_2011}, then the voltage drop is to be determined. Furthermore, as the cell voltage hinges on the difference $V_3-V_2$ rather than on the absolute magnitude of the electron potential,  one of the boundary values of $\phi_\text{el}$ ($V_2$ or $V_3$) can be taken in an arbitrary way.

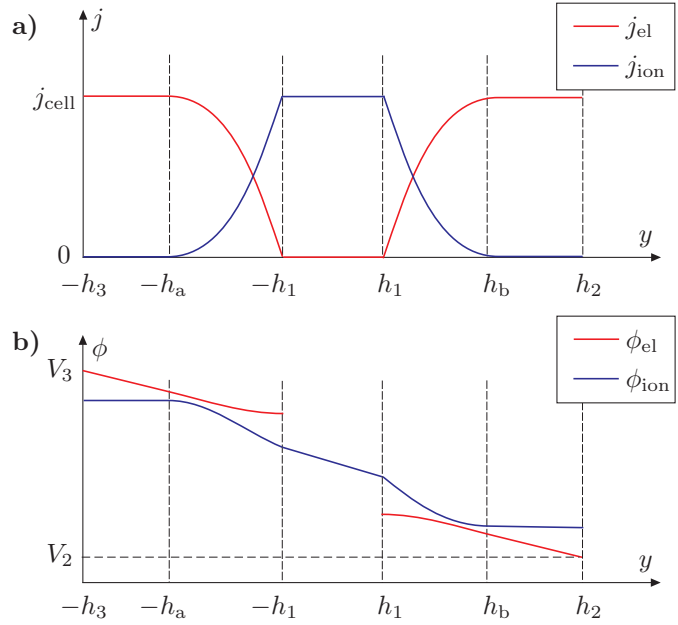
\begin{figure}
\begin{center}
\includegraphics[scale=0.80]{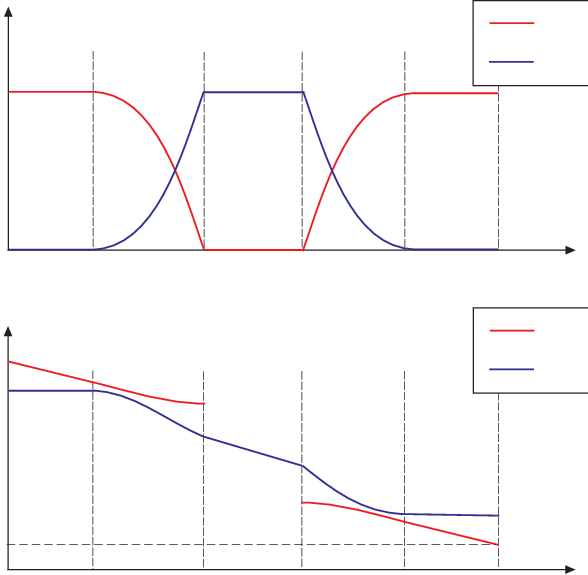}
\put(-250,210){$\textbf{a)}$}
\put(-250,90){$\textbf{b)}$}
\put(-18,210){$j_\text{el}$}
\put(-18,195){$j_\text{ion}$}
\put(-18,91){$\phi_\text{el}$}
\put(-18,76){$\phi_\text{ion}$}
\put(-233,123){$0$}
\put(-232,112){$-h_3$}
\put(-202,112){$-h_\text{a}$}
\put(-160,112){$-h_1$}
\put(-112,112){$h_1$}
\put(-72,112){$h_\text{b}$}
\put(-37,112){$h_2$}
\put(-13,131){$y$}
\put(-242,183){$j_\text{cell}$}
\put(-220,213){$j$}
\put(-220,88){$\phi$}
\put(-13,8){$y$}
\put(-232,-10){$-h_3$}
\put(-202,-10){$-h_\text{a}$}
\put(-160,-10){$-h_1$}
\put(-112,-10){$h_1$}
\put(-72,-10){$h_\text{b}$}
\put(-37,-10){$h_2$}
\put(-237,10){$V_2$}
\put(-237,80){$V_3$}

\caption{The schematic view of the expected qualitative behaviour of: a) electron and ion currents, b) electron and ion potentials. }
\label{j_fi_gen_1}
\end{center}
\end{figure}

Respective boundary conditions and expected distributions of currents and potentials are schematically depicted in Fig.\ref{j_fi_gen_1}.

As for the mass transfer boundary conditions, we specify the mass fractions of respective species at the external interfaces of the electrodes:
\begin{equation}
\label{C_BCs}
C_{\text{H}_2}(-h_3)=C_{\text{H}_2}^\text{bulk}, \quad C_{\text{H}_2\text{O}}(-h_3)=C_{\text{H}_2\text{O}}^\text{bulk}, \quad C_{\text{O}_2}(h_2)=C_{\text{O}_2}^\text{bulk}.
\end{equation}
Obviously, the above conditions are  merely a natural consequence of reducing the cell model to the system of electrodes. When analyzing  a complete problem with the fuel and air channels, conditions \eqref{C_BCs} are to be replaced by respective transmission conditions. Furthermore, as the electrolyte constitutes an impermeable barrier for gas diffusion, the following natural boundary conditions are to be imposed at respective interfaces:
\begin{equation}
\label{J_BCs}
J_{\text{H}_2}(-h_1)=0, \quad J_{\text{H}_2\text{O}}(-h_1)=0, \quad J_{\text{O}_2}(h_1)=0.
\end{equation}

When integrating equations \eqref{ch2_ODE}, \eqref{ch2o_ODE} and \eqref{co2_ODE} over the span of respective electrodes while accounting for \eqref{j_bal_int} and \eqref{J_BCs}, one arrives at the formulae for the mass fluxes at external interfaces of the anode and the cathode:
\begin{equation}
\label{J_ext_BCs}
J_{\text{H}_2}(-h_3)=\frac{M_{\text{H}_2}}{2F}j_\text{cell}, \quad J_{\text{H}_2\text{O}}(-h_3)=-\frac{M_{\text{H}_2\text{O}}}{2F}j_\text{cell}, \quad J_{\text{O}_2}(h_2)=-\frac{M_{\text{O}_2}}{4F}j_\text{cell}.
\end{equation}

Expressions \eqref{J_ext_BCs} can be treated as boundary conditions in that variant of the problem when the cell current density, $j_\text{cell}$, is assumed to be known. However, even when the voltage drop is predefined, \eqref{J_ext_BCs} should be employed in the computational algorithm to provide the consistency of the results. Remark that the fluxes' magnitudes defined by  \eqref{J_ext_BCs} hold throughout the whole backing zones of respective electrodes, as there are no mass sources/sinks in these regions.

\begin{figure}
\begin{center}
\includegraphics[scale=0.80]{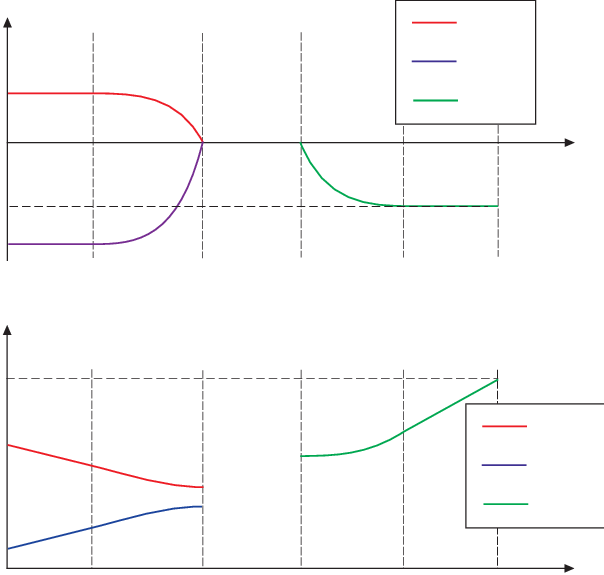}
\put(-250,210){$\textbf{a)}$}
\put(-250,90){$\textbf{b)}$}
\put(-52,210){$J_{\text{H}_2}$}
\put(-52,195){$J_{\text{H}_2\text{O}}$}
\put(-52,180){$J_{\text{O}_2}$}
\put(-25,55){$C_{\text{H}_2}$}
\put(-25,40){$C_{\text{H}_2\text{O}}$}
\put(-25,25){$C_{\text{O}_2}$}
\put(-237,163){$0$}
\put(-240,112){$-h_3$}
\put(-210,112){$-h_\text{a}$}
\put(-168,112){$-h_1$}
\put(-120,112){$h_1$}
\put(-80,112){$h_\text{b}$}
\put(-45,112){$h_2$}
\put(-20,171){$y$}
\put(-264,183){$\frac{M_{\text{H}_2}j_\text{cell}}{2F}$}
\put(-272,140){$-\frac{M_{\text{O}_2}j_\text{cell}}{2F}$}
\put(-278,124){$-\frac{M_{\text{H}_2\text{O}}j_\text{cell}}{2F}$}
\put(-225,208){$J$}
\put(-225,88){$C$}
\put(-20,8){$y$}
\put(-240,-10){$-h_3$}
\put(-212,-10){$-h_\text{a}$}
\put(-168,-10){$-h_1$}
\put(-120,-10){$h_1$}
\put(-80,-10){$h_\text{b}$}
\put(-45,-10){$h_2$}
\put(-257,70){$C_{\text{O}_2}^\text{bulk}$}
\put(-257,50){$C_{\text{H}_2}^\text{bulk}$}
\put(-257,5){$C_{\text{H}_2\text{O}}^\text{bulk}$}

\caption{The schematic view of the expected qualitative behaviour of: a) mass fluxes, b) mass fractions of respective species. }
\label{J_c_gen_1}
\end{center}
\end{figure}
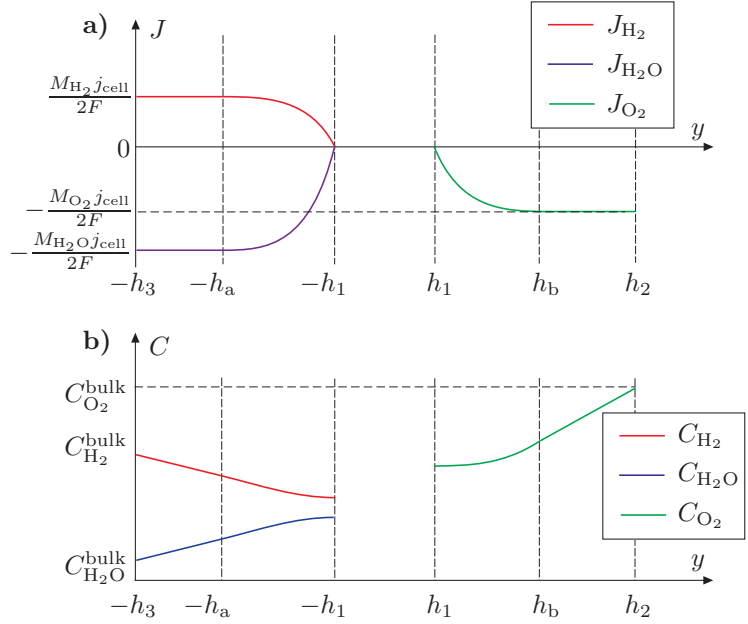

Respective boundary conditions and the schematic distributions of mass fractions and fluxes are depicted in Fig.\ref{J_c_gen_1}.

Provided that the above  boundary conditions hold, proper transmission conditions between respective subdomains are to be strictly satisfied in order to preserve the correct mathematical structure of the solution. Generally, we assume continuity of respective fields and fluxes (currents) at the subdomain's interfaces. Due to the fact that the electrolyte is considered a perfect electron insulator ($\sigma_\text{el}=0$ S$\cdot$m$^{-1}$), equation \eqref{phi_el_ODE} degenerates for $y \in (-h_1,h_1)$, and thus one cannot determine $\phi_\text{el}$ over this section. However, the respective transmission conditions for the electron current (zero current $j_\text{el}$) still hold and are to be employed. The transmission conditions for the electron and ion currents are:
\begin{equation}
\label{j_el_ion_zero}
j_\text{el}(-h_1)=j_\text{el}(h_1)=0, \quad j_\text{ion}(-h_1)=j_\text{ion}(h_1)=j_\text{cell}.
\end{equation}
The jump of electron potential over the electrolyte thickness is to be found as an element of  the solution. On the other hand, the ionic potential, $\phi_\text{ion}$ is a continuous $C^1$- class function in the entire domain. Pertinent continuity conditions can be identified in Figs.\ref{j_fi_gen_1}--\ref{J_c_gen_1}.

\subsection{The electrochemical reactions model}

The standard description of the electrochemical reactions in the electrodes can be found in many publications \cite{Li_2005,Kulikovsky_2010,Brus_et_al_2017,Onaka_2017}. Despite some differences in definition of respective elements, all these models exhibit similar qualitative behavior. Thus, without loss of generality, we will base the following analysis on the very formulation given in \cite{Brus_et_al_2017,Onaka_2017}.

Let the charge transfer current be specified as:
\begin{equation}
\label{i_ct_gen}
i^\text{ct} =
  \begin{cases}
	    i^\text{ct}_\text{ano}  & \quad \text{for } \quad y\in(-h_\text{a},-h_1),\\
			i^\text{ct}_\text{cat}       & \quad \text{for } \quad y\in(h_1,h_\text{b}),
  \end{cases}
\end{equation}
where
\begin{equation}
\label{i_ano_def}
i^\text{ct}_\text{ano}=i^\text{tpb}_\text{0,ano}l^\text{tpb}_\text{ano}\left[\exp\left(\frac{2F\eta^\text{act}_\text{ano}}{RT} \right)-\exp\left(-\frac{F\eta^\text{act}_\text{ano}}{RT} \right) \right],
\end{equation}
\begin{equation}
\label{i_cat_def}
i^\text{ct}_\text{cat}=i^\text{dpb}_\text{0,cat}A^\text{dpb}_\text{cat}\left[\exp\left(\frac{1.2F\eta^\text{act}_\text{cat}}{RT} \right)-\exp\left(-\frac{F\eta^\text{act}_\text{cat}}{RT} \right) \right].
\end{equation}

Pertinent symbols in the above equations denote: $l^\text{tpb}_\text{ano}$ [m$\cdot$m$^{-3}$]- the length of the two-phase boundary, $A^\text{dpb}_\text{cat}$ [m$^2\cdot$m$^{-3}$] - the area of the three-phase boundary, $i^\text{tpb}_\text{0,ano}$  [A$\cdot$m$^{-1}$]- the anode exchange current density, $i^\text{dpb}_\text{0,cat}$  [A$\cdot$m$^{-2}$]- the cathode exchange current density, $R$ [J$\cdot$mol$^{-1}\cdot$ K$^{-1}$]- the universal gas constant, $T$ [K]-the temperature, $\eta^\text{act}_\text{ano}$ [V], $\eta^\text{act}_\text{cat}$ [V]- the activation overpotentials of anode and cathode, respectively. According to \cite{Brus_et_al_2017}, the exchange current densities for the Ni-YSZ anode and LSCF cathode can be computed as:
\begin{equation}
\label{i_0_ano_def}
i^\text{tpb}_\text{0,ano}=32.4\cdot p_{\text{H}_2}^{-0.03}p_{\text{H}_2\text{O}}^{0.4}\exp \left(-\frac{152155}{RT}\right),
\end{equation}
\begin{equation}
\label{i_0_cat_def}
i^\text{dpb}_\text{0,cat}=1.47\cdot 10^6 \cdot p_{\text{O}_2}^{0.2}\exp \left(-\frac{85859}{RT}\right),
\end{equation}
where $p_j$ [Pa] denotes the partial pressure of $j$-th species.

The activation overpotentials are defined by the following formulae:
\begin{equation}
\label{eta_act_ano}
\eta_\text{ano}^\text{act}=\phi_\text{el}-\phi_\text{ion}-\eta_\text{ano}^\text{conc},
\end{equation}
\begin{equation}
\label{eta_act_cat}
\eta_\text{cat}^\text{act}=\phi_\text{ion}-\phi_\text{el}-\eta_\text{cat}^\text{conc},
\end{equation}
with the concentration overpotentials given by:
\begin{equation}
\label{eta_conc_ano}
\eta_\text{ano}^\text{conc}=\frac{RT}{2F}\ln\left(\frac{p_{\text{H}_2}^\text{bulk}}{p_{\text{H}_2}}\frac{p_{\text{H}_2\text{O}}}{p_{\text{H}_2\text{O}}^\text{bulk}}\right),
\end{equation}
\begin{equation}
\label{eta_conc_cat}
\eta_\text{cat}^\text{conc}=\frac{RT}{4F}\ln\left(\frac{p_{\text{O}_2}^\text{bulk}}{p_{\text{O}_2}}\right).
\end{equation}

\begin{remark}
Note that in the above expressions the partial pressures of respective species are used (to have  a direct reference to the cited papers), while the system of governing equations is based on the mass fractions. The partial pressures will be transformed in the computations to the mass fractions by the ideal gas law. Moreover, in relations \eqref{eta_conc_ano}-\eqref{eta_conc_cat} a simple substitution of $p_j$ with $C_j$ is what one needs.
\end{remark}

\section{Modified formulation of the SOFC problem}
\label{mod_form}

The above system of equations constitutes the classical mathematical formulation of the SOFC problem. Although, depending on the approach employed, some elements (e.g. the diffusion equations, formulae for the exchange current densities) can be different from those specified in the previous section, the general characteristics of the problem remain the same. In the following we will prove that this standard formulation leads to locally unphysical effects. Moreover, we will show that the definition of the active layer thickness plays a very important role and cannot be assumed in an arbitrary way. In order to tackle these issues, pertinent modifications to the classical system will be proposed.

\subsection{Comments on application of the Butler-Volmer equation}
\label{BV_comments}

As results from the system of governing equations \eqref{phi_el_ODE}, \eqref{phi_ion_ODE}, \eqref{ch2_ODE}, \eqref{ch2o_ODE}, \eqref{co2_ODE}, respective electrochemical reactions are assumed to take place only in the active layers of the electrodes (i.e. $y \in(-h_\text{a},-h_1)$ for the anode and $y \in(h_1,h_\text{b})$ for the cathode).

In some papers one can find an approach where the active layer is a priori stretched over the whole electrode thickness. However, such a method does not agree with observations (it is commonly accepted that the electrochemical reactions occur only in the immediate vicinity of the electrode/electrolyte interface \cite{Nam_2017,Miyawaki_2014,Kishimoto_2013}), producing unphysical effects such as conduction of ionic current up to the external boundary of the electrode. Moreover, the expression \eqref{i_ct_gen} constitutes an injective function of the activation overpotential (hence, for monotonic $\eta^\text{act}$ an injective function of $y$ as well), and thus one cannot expect it to suppress the charge transfer current over the backing (inactive) part of the electrode.

For this reason, it is important to make a clear distinction between the backing (electrochemically inactive) and catalytic (electrochemically active) layers of respective electrodes, with the source terms in the governing ODEs specified accordingly. At the interfaces of these layers the following conditions have to be satisfied:
\begin{equation}
\label{eta_act_an_con}
\eta_\text{ano}^\text{act}(-h_\text{a})=0,
\end{equation}
\begin{equation}
\label{eta_act_cat_con}
\eta_\text{cat}^\text{act}(h_\text{b})=0.
\end{equation}

Note that \eqref{eta_act_an_con}--\eqref{eta_act_cat_con} are the necessary conditions for the source terms in the governing ODEs (RHSs of \eqref{phi_el_ODE}, \eqref{phi_ion_ODE}, \eqref{ch2_ODE}, \eqref{ch2o_ODE}, \eqref{co2_ODE}) to change continuously with $y$. Moreover, when the above separation of the electrodes' layers is introduced alongside the boundary conditions \eqref{j_ion_BC}, the ionic current turns identically to zero in the backing zones. This effectively means that:
\begin{equation}
\label{fi_ion_an_back}
\frac{\text{d}\phi_\text{ion}}{\text{d}y}=0, \quad \phi_\text{ion}=\text{const.}=\phi_\text{a}, \quad \text{for} \quad y\leq -h_\text{a},
\end{equation}
\begin{equation}
\label{fi_ion_cat_back}
\frac{\text{d}\phi_\text{ion}}{\text{d}y}=0, \quad \phi_\text{ion}=\text{const.}=\phi_\text{b}, \quad \text{for} \quad y\geq h_\text{b}.
\end{equation}

\begin{remark}

If one assumes that whole electrodes are occupied by the active zones ($h_\text{a}=h_3$, $h_\text{b}=h_2$), then conditions \eqref{eta_act_an_con}--\eqref{eta_act_cat_con} combined with the boundary conditions \eqref{C_BCs} yield the following identities:
\begin{equation}
\label{fi_ion_el_BC}
\phi_\text{a}=V_3, \quad \phi_\text{b}=V_2.
\end{equation}
However, as mentioned previously, in such a case some unphysical effects appear.
\end{remark}

Now, let us focus on the behaviour of the charge-transfer current function, $i^\text{ct}$, in the neighbourhood of the interfaces between the catalyst and backing layers. When analyzing respective problem for the cathode, it shows that as moving away from $y=h_\text{b}$ towards the electrolyte ($y<h_\text{b}$), the concentration overpotential grows due to consumption of oxygen (compare \eqref{eta_conc_cat} - for monotonic change of the oxygen mass fraction - see also Fig.\ref{J_c_gen_1}). Simultaneously, the electron potential increases (as for any positive value of the cell current $\phi_\text{el}'<0$). By recalling the condition \eqref{fi_ion_cat_back}$_1$ we come to a conclusion that over at least some small distance in the proximity of the interface the activation overpotential is negative because:
\begin{equation}
\label{eta_neg_1}
\phi_\text{ion}<\phi_\text{el}+\eta_\text{cat}^\text{conc}, \quad \text{for} \quad h_\text{b}-\varepsilon_\text{b} <y < h_\text{b}, \quad \varepsilon_\text{b} \ll h_\text{b}-h_1.
\end{equation}
Thus, also the charge-transfer current assumes in this range negative values:
\begin{equation}
\label{ic_cat_neg}
i^\text{ct}<0, \quad \text{for} \quad h_\text{b}-\varepsilon_\text{b} <y < h_\text{b}, \quad \varepsilon_\text{b} \ll h_\text{b}-h_1.
\end{equation}
Note that this effect is inevitable when imposing respective transmission conditions at $h_\text{b}$. It can be avoided only when $i^\text{ct}$ turns to zero over the whole domain (no current flow).

When consequently employing the same manner of reasoning for the anode (the interface defined by $y=-h_\text{a}$), we see again that  the concentration overpotential grows for growing $y$ (i.e. when moving towards the electrolyte layer). Simultaneously, the electron potential declines (as for any positive cell current $\phi_\text{el}'<0$). Finally, considering the condition \eqref{fi_ion_an_back}$_1$ we come to the following result:
\begin{equation}
\label{eta_neg_11}
\phi_\text{el}<\phi_\text{ion}+\eta_\text{ano}^\text{conc}, \quad \text{for} \quad -h_\text{a} <y < -h_\text{a}+\varepsilon_\text{a}, \quad \varepsilon_\text{a} \ll h_\text{a}-h_1,
\end{equation}
which implies that:
\begin{equation}
\label{ic_n_neg}
i^\text{ct}<0, \quad \text{for} \quad -h_\text{a} <y < -h_\text{a}+\varepsilon_\text{a}, \quad \varepsilon_\text{a} \ll h_\text{a}-h_1.
\end{equation}
As mentioned previously, this peculiarity is a direct consequence of the accepted mathematical model. It  does not appear only when the cell current is zero.

\begin{figure}
\begin{center}
\includegraphics[scale=0.80]{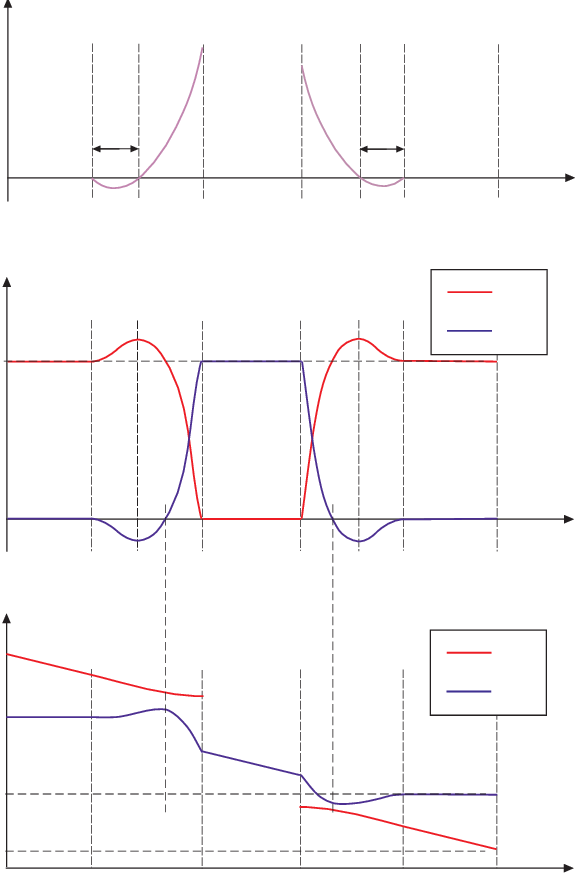}
\put(-250,330){$\textbf{a)}$}
\put(-250,210){$\textbf{b)}$}
\put(-250,90){$\textbf{c)}$}
\put(-215,330){$\eta^\text{act}$}
\put(-10,274){$y$}
\put(-228,249){$-h_3$}
\put(-197,249){$-h_\text{a}$}
\put(-155,249){$-h_1$}
\put(-108,249){$h_1$}
\put(-67,249){$h_\text{b}$}
\put(-32,249){$h_2$}
\put(-28,222){$j_\text{el}$}
\put(-28,207){$j_\text{ion}$}
\put(-28,83){$\phi_\text{el}$}
\put(-28,68){$\phi_\text{ion}$}
\put(-230,132){$0$}
\put(-228,116){$-h_3$}
\put(-197,116){$-h_\text{a}$}
\put(-155,116){$-h_1$}
\put(-108,116){$h_1$}
\put(-67,116){$h_\text{b}$}
\put(-32,116){$h_2$}
\put(-13,141){$y$}
\put(-242,197){$j_\text{cell}$}
\put(-212,218){$j$}
\put(-180,284){$\varepsilon_\text{a}$}
\put(-80,284){$\varepsilon_\text{b}$}
\put(-213,92){$\phi$}
\put(-13,8){$y$}
\put(-228,-8){$-h_3$}
\put(-197,-8){$-h_\text{a}$}
\put(-155,-8){$-h_1$}
\put(-108,-8){$h_1$}
\put(-67,-8){$h_\text{b}$}
\put(-32,-8){$h_2$}
\put(-237,6){$V_2$}
\put(-237,83){$V_3$}
\put(-237,30){$\phi_\text{b}$}
\put(-237,60){$\phi_\text{a}$}

\caption{The schematic view of the qualitative behaviour of: a) activation overpotential, b) electron and ion currents, c) electron and ion potentials. The zones of negative activation overpotentials are defined by $\varepsilon_\text{a}$ and $\varepsilon_\text{b}$, respectively.}
\label{j_fi_2}
\end{center}
\end{figure}

\begin{figure}
\begin{center}
\includegraphics[scale=0.80]{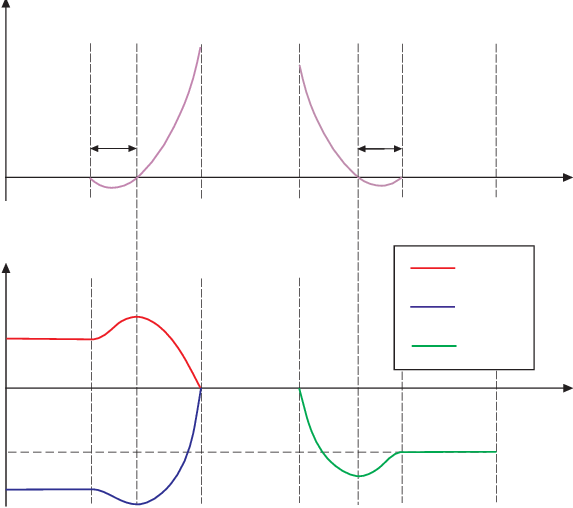}
\put(-250,190){$\textbf{a)}$}
\put(-250,90){$\textbf{b)}$}
\put(-180,144){$\varepsilon_\text{a}$}
\put(-80,144){$\varepsilon_\text{b}$}
\put(-213,185){$\eta^\text{act}$}
\put(-12,135){$y$}
\put(-230,106){$-h_3$}
\put(-198,106){$-h_\text{a}$}
\put(-156,106){$-h_1$}
\put(-110,106){$h_1$}
\put(-67,106){$h_\text{b}$}
\put(-32,106){$h_2$}
\put(-40,89){$J_{\text{H}_2}$}
\put(-40,74){$J_{\text{H}_2\text{O}}$}
\put(-40,59){$J_{\text{O}_2}$}
\put(-230,125){$0$}
\put(-260,63){$\frac{M_{\text{H}_2}j_\text{cell}}{2F}$}
\put(-268,20){$-\frac{M_{\text{O}_2}j_\text{cell}}{4F}$}
\put(-274,4){$-\frac{M_{\text{H}_2\text{O}}j_\text{cell}}{2F}$}
\put(-213,85){$J$}
\put(-12,53){$y$}
\put(-230,-10){$-h_3$}
\put(-198,-10){$-h_\text{a}$}
\put(-176,-10){$-h_1$}
\put(-110,-10){$h_1$}
\put(-85,-10){$h_\text{b}$}
\put(-32,-10){$h_2$}

\caption{The schematic view of the qualitative behaviour of: a) activation overpotential, b) mass fluxes. The zones of negative activation overpotentials are defined by $\varepsilon_\text{a}$ and $\varepsilon_\text{b}$, respectively.}
\label{J_c_2}
\end{center}
\end{figure}
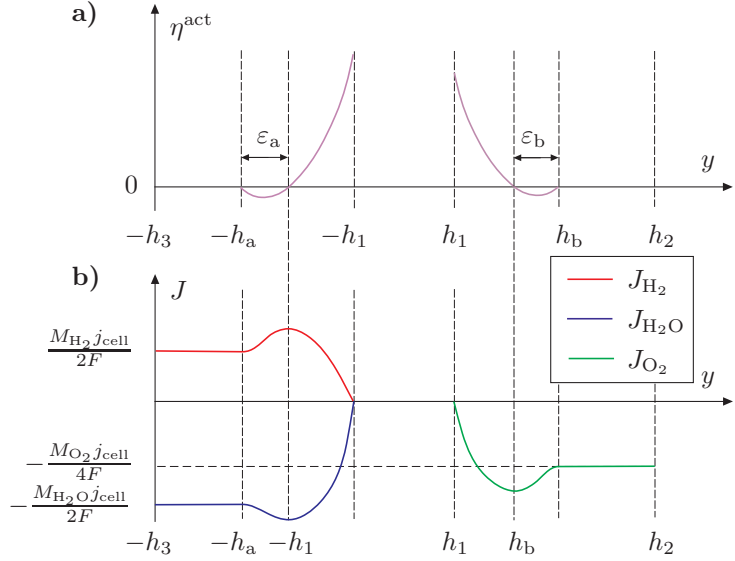

The above observations about the charge-transfer current have very serious consequences for the mathematical structure of solution and its physical interpretation. As for the latter aspect, it means that over some  section of each electrode the backward mechanism of respective electrochemical reaction prevails over its forward component (production of hydrogen and oxygen instead of consumption, and absorption of water in lieu of its generation). It is even odder when one realizes that in the immediate vicinity of the interfaces between the backing and catalyst layers, the electron current is supposed to exceed the magnitude of overall cell current ($j_\text{cell}$), while the ionic current has to assume negative values (reversion of ions' flow !). The consequences of the negative activation overpotential (and thus the negative charge transfer current) are schematically depicted in Figs. \ref{j_fi_2}-\ref{J_c_2}. Respective zones of $i^\text{ct}$ have been intentionally rescaled to better present the  related effects. Note that the ion and electron current functions are not monotonic. Moreover, due to the negative values of the ionic current in the proximity of $h=-h_\text{a}$ and $y=h_\text{b}$, the ionic potential is not a monotonic function as well. As shown in Fig.\ref{J_c_2}, the mass fluxes of corresponding species exceed their boundary values \eqref{J_ext_BCs}, producing counterintuitive and unphysical results. Admittedly, the mass fraction functions remain monotonic, for the pertinent fluxes being of the same sign (unidirectional mass flow) over corresponding electrodes  (for this reason we do not depict the mass fractions of species in the graph).

Let us emphasize one more time that these locally unphysical effects are inherent to the classical model of the cell and cannot be avoided if one wants to satisfy the pertinent boundary and transmission conditions. Moreover, the negative values of $i^\text{ct}$ will be obtained even if more advanced diffusion models are employed (e.g. the Dusty Gas Model) and thus, the above general considerations still hold in such cases.

Note that the negative charge transfer current can be eliminated at a cost of relaxing some boundary conditions, which however leads to even less physically reliable results (e.g. non-zero ionic current at the interconnects).
If one assumes that the classical model, although locally unphysical, is anyway credible in a global sense, there still exists a need to implement the aforementioned peculiarities in the computational schemes as any accurate and efficient numerical simulation of the problem needs to account for the proper solution structure.

Yet in our analysis we shall introduce a modification of the classical model. Namely, we assume that the activation overpotential has always a non-negative value:
\begin{equation}
\label{eta_act_ano_abs}
\eta_\text{ano}^\text{act}=|\phi_\text{el}-\phi_\text{ion}-\eta_\text{ano}^\text{conc}|,
\end{equation}
\begin{equation}
\label{eta_act_cat_abs}
\eta_\text{cat}^\text{act}=|\phi_\text{ion}-\phi_\text{el}-\eta_\text{cat}^\text{conc}|.
\end{equation}
Although this amendment is only an ad hoc solution, it completely eliminates the unphysical effects reported above. Similar result could be achieved by taking the absolute value of the charge-transfer current instead of \eqref{i_ct_gen}:
\begin{equation}
\label{i_ct_gen_abs}
i^\text{ct} =
  \begin{cases}
	    |i^\text{ct}_\text{ano}|  & \quad \text{for } \quad y\in(-h_\text{a},-h_1),\\
			|i^\text{ct}_\text{cat}|       & \quad \text{for } \quad y\in(h_1,h_\text{b}).
  \end{cases}
\end{equation}
Both approaches are equivalent for small activation overpotentials, i.e. when $\varepsilon_\text{a}$ and $\varepsilon_\text{b}$ tend to zero. When employing any of the above modifications, the activation overpotentials and the charge-transfer current assume zero values at the external boundaries of the catalyst layer ($y=-h_\text{a}$ and $y=h_\text{b}$) and at some distance inside these layers ($y=-h_\text{a}+\varepsilon_\text{a}$ and $y=h_\text{b}-\varepsilon_\text{b}$). However, all the dependent variables are now monotonic functions of $y$, while the mass and charge transfer becomes solely unidirectional.

\begin{remark}
It may look tempting to use another modification instead of \eqref{eta_act_ano_abs}-\eqref{eta_act_cat_abs} or \eqref{i_ct_gen_abs}, namely to set $i^\text{ct}=0$ over this part of the electrode where the negative charge transfer current is obtained. However, if such an amendment is employed, it introduces  essential violation of the system of governing equations, as neither of them can be satisfied in the proximity of the external interfaces of the electrodes. Moreover, the global balances of mass and charge will not be fulfilled either. Thus, it is more reasonable to use the standard formulation with locally unphysical effects than to utilize the aforementioned modification.
\end{remark}

\subsection{The catalyst layers thickness - the problem of uniqueness of solution}
\label{cat_sub}

The catalyst (active) layer is a part of respective electrode in which the electrochemical reactions take place. Its thickness is not constant and depends on many parameters including material and microstructural properties and operating conditions. As stated previously, in numerical computations there are usually two strategies adopted to introduce the layer (we do not discuss here the simplified model where the electrochemical reactions are assumed to take place only at the interfaces of electrolyte and respective electrodes). It is either stretched across the  whole electrode thickness  or a predefined size of the active zone (established experimentally) is imposed. In the first case it is accepted that the 'real' layer's thickness pertains to the part of the electrode over which a certain amount of the total cell current (one can find the following figures in literature: $90\%$, $95\%$, $99\%$ - see e.g. \cite{Nam_2017}) is produced.
The second approach seems much more reasonable as it accounts for a real physical feature of the problem, even if experimental determination of the active layer thickness and its changes constitutes a formidable task. However, in the following we shall prove that when accepting the most straightforward definition of the active layer (the exclusive zone in which the electrochemical reactions occur) alongside respective boundary and transmission conditions (see Section \ref{BCs_sect}), the thickness of the catalyst layer cannot be arbitrary and is a component of the solution.

\begin{figure}
\begin{center}
\includegraphics[scale=0.80]{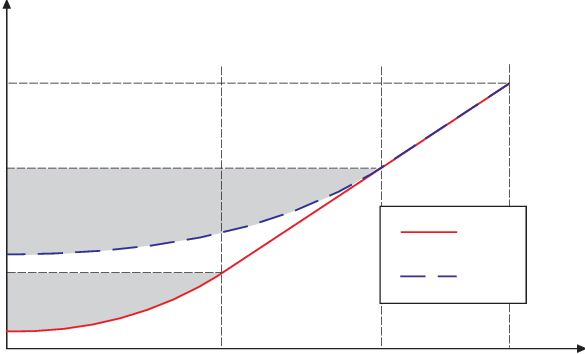}
\put(-44,44){$C_{\text{O}_2}^{\text{(I)}}$}
\put(-44,26){$C_{\text{O}_2}^{\text{(II)}}$}
\put(-220,125){$C_{\text{O}_2}$}
\put(-252,102){$C_{\text{O}_2}^\text{bulk}$}
\put(-252,72){$C_{\text{O}_2}^{\text{b(II)}}$}
\put(-252,30){$C_{\text{O}_2}^{\text{b(I)}}$}
\put(-12,8){$y$}
\put(-225,-10){$h_1$}
\put(-82,-10){$h_\text{b}'$}
\put(-145,-10){$h_\text{b}$}
\put(-32,-10){$h_2$}

\caption{The schematic view of the mass fraction of oxygen in the cathode.  Solid line refers to the variant when the catalyst layer is defined by $h_\text{b}$, while the dashed one pertains to the active zone being delimited by  $h_\text{b}'$. The same magnitude of $C_{\text{O}_2}^\text{bulk}$ is assumed. The shaded areas refer to the amount of oxygen consumed in respective cases.}
\label{co2_int}
\end{center}
\end{figure}
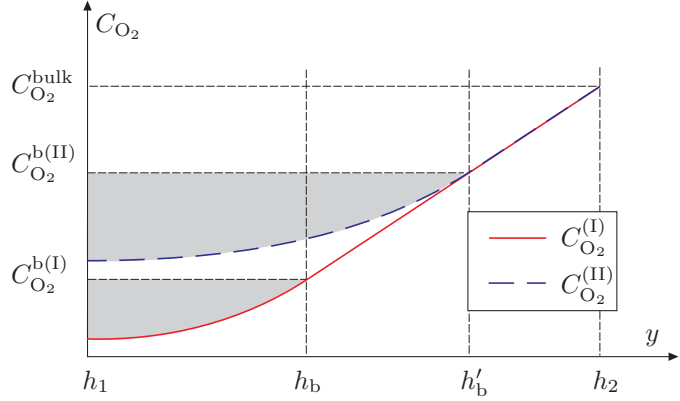

To this end let us consider first the respective sub-problem for the cathode. Assume that for some thickness of the active layer, defined by the coordinate $h_\text{b}$, there exists a unique solution to the problem. This solution will henceforth be marked  by the superscript $(\text{I})$. Now, imagine that the active layer thickness has increased to the dimension bounded by the coordinate $h_\text{b}'>h_\text{b}$ (see Fig.\ref{co2_int}), with all the remaining parameters (i.e. the cell voltage and current, boundary conditions) the same as previously. The new solution will be singled out by the superscript $(\text{II})$. The necessary condition for the new solution to produce the same cell current, $j_\text{cell}$, is:
\begin{equation}
\label{i_ct_cond_proof}
\int_{h_1}^{h_\text{b}}i^{\text{ct(I)}}dy=\int_{h_1}^{h'_\text{b}}i^{\text{ct(II)}}dy,
\end{equation}
which implies that:
\begin{equation}
\label{i_ct_inq}
i^\text{ct(I)}>i^{\text{ct(II)}},
\end{equation}
over most of the interval $h_1<y<h_\text{b}$.
In the following analysis we will verify to what degree the condition \eqref{i_ct_inq} can be satisfied, if the catalyst layer thickness is changed, without any interference in the remaining problem parameters.

In Fig. \ref{co2_int} we have depicted symbolically the respective solutions for the mass fraction of oxygen ($C_{\text{O}_2}^{\text{(I)}}$ and $C_{\text{O}_2}^{\text{(II)}}$). Note that:
\begin{equation}
\label{Co_I_II}
C_{\text{O}_2}^{\text{(I)}}=C_{\text{O}_2}^{\text{(II)}},\quad \text{for} \quad y\geq h_\text{b}',
\end{equation}
since in the backing zone of the cathode both functions have to satisfy the same ODE (diffusion without mass sources) and boundary conditions, regardless of the diffusion model employed.  Moreover, bearing in mind that the fluxes of respective species outside the catalyst layers are constant, one can infer from \eqref{J_ext_BCs} that:
\begin{equation}
\label{DC_eq}
\frac{\text{d}C_{\text{O}_2}^{\text{(I)}}}{\text{d}y}\big|_{y=h_\text{b}}=\frac{\text{d}C_{\text{O}_2}^{\text{(II)}}}{\text{d}y}\big|_{y=h'_\text{b}}.
\end{equation}

Remark that for the same value of the cell total current, $j_\text{cell}$, the amount of oxygen consumed should be identical in both cases, and thus respective shaded areas in Fig.\ref{co2_int} need to be equal.

In the following analysis we assume that mass fraction functions are monotonic and:
\begin{equation}
\label{C_I_II_inq}
C_{\text{O}_2}^{\text{(I)}}<C_{\text{O}_2}^{\text{(II)}}, \quad \text{for} \quad y \in(h_1,h_\text{b}).
\end{equation}
This assumption results from the fact that (compare \eqref{C^b_def}):
\begin{equation}
\label{C_I_II_b_dif}
C_{\text{O}_2}^{\text{b(II)}}-C_{\text{O}_2}^{\text{b(I)}}=\frac{M_{\text{O}_2}j_\text{cell}}{4F\rho_\text{a}D_2}(h'_\text{b}-h_\text{b})>0,
\end{equation}
while
\begin{equation}
\label{dC_inq}
\frac{\text{d}C_{\text{O}_2}^{\text{(I)}}}{\text{d}y} \geq \frac{\text{d}C_{\text{O}_2}^{\text{(II)}}}{\text{d}y}, \quad \text{for} \quad y\in(h_1,h_b).
\end{equation}
Equation \eqref{dC_inq} can be deduced from \eqref{cat_ax_2} when using the electric current function in the left hand side instead of the potential derivative.

The charge transfer current at the cathode is proportional to the following product:
\begin{equation}
\label{i_cat_pr}
i^\text{ct}_\text{cat} \sim C_{\text{O}_2}^{0.2}\left[\exp\left(\frac{1.2F\eta^\text{act}_\text{cat}}{RT} \right)-\exp\left(-\frac{F\eta^\text{act}_\text{cat}}{RT} \right) \right].
\end{equation}
As results from \eqref{C_I_II_inq}, the first term of the above product is smaller for the solution $\text{(I)}$. Thus, if \eqref{i_ct_inq} is to be satisfied at some distance away from $h_\text{b}$ (it cannot be fulfilled in the immediate vicinity of the interface, as $i^\text{ct(I)}(h_\text{b})=0$), then the second term should overweight the decline of the first one (note that, due to \eqref{eta_conc_cat} and \eqref{C_I_II_inq}: $\eta_\text{cat}^\text{conc(I)}>\eta_\text{cat}^\text{conc(II)}$ for $y\in(h_1,h_\text{b})$) . This can be achieved only if:
\begin{equation}
\label{eta_cat_inq}
\eta^\text{act(I)}_\text{cat}>\eta^\text{act(II)}_\text{cat},
\end{equation}
over most of the span of the catalyst layer. In order to verify the above condition let us analyze the expression:
\begin{equation}
\label{delta_eta_cat}
\Delta \eta^\text{act}_\text{cat}=\eta^\text{act(II)}_\text{cat}-\eta^\text{act(I)}_\text{cat}=\phi_\text{ion}^\text{(II)}-\phi_\text{ion}^\text{(I)}+\phi_\text{el}^\text{(I)}-\phi_\text{el}^\text{(II)}+\eta^\text{conc(I)}_\text{cat}-\eta^\text{conc(II)}_\text{cat}.
\end{equation}
For the time being we neglect the absolute value in the activation overpotential introduced in \eqref{eta_act_cat_abs} as we are primarily interested in the bahviour of $\Delta \eta^\text{act}_\text{cat}$ away from the interface between the catalyst and the backing layer. However, this issue will be commented later on. By combining equations \eqref{eta_conc_cat}, \eqref{fi_el_cat_fin} and \eqref{fi_ion_cat} according to \eqref{delta_eta_cat} one arrives at the following formula for $\Delta \eta^\text{act}_\text{cat}$:
\begin{equation}
\label{delta_eta_cat_fin}
\begin{split}
\Delta \eta^\text{act}_\text{cat}=\frac{4F\rho_\text{a}D_\text{2}}{M_{\text{O}_2}}\left(\frac{1}{\sigma_\text{el}}+\frac{1}{\sigma_\text{ion}} \right)\left( C_{\text{O}_2}^\text{(II)}(y)-C_{\text{O}_2}^\text{(I)}(y)\right)+\frac{RT}{4F}\ln\left( \frac{C_{\text{O}_2}^\text{(II)}(y)}{C_{\text{O}_2}^\text{(I)}(y)}\right)+\\
\frac{RT}{4F}\ln\left( \frac{C_{\text{O}_2}^\text{b(I)}(y)}{C_{\text{O}_2}^\text{b(II)}(y)}\right)+\frac{j_\text{cell}}{\sigma_\text{el}}(h_\text{b}-h'_\text{b}).
\end{split}
\end{equation}
The first two terms in the above equation are positive and monotonically decreasing functions of $y$. The second two terms form a negative constant. Note that if we assume that activation overpotential is positive throughout the whole catalyst layer (as in the modification \eqref{eta_act_cat_abs}), then:
\begin{equation}
\label{eta_act_hb}
\Delta \eta^\text{act}_\text{cat}(h_\text{b})=\eta^\text{act(II)}_\text{cat}(h_\text{b})>0,
\end{equation}
which implies that:
\begin{equation}
\label{eta_act_pos}
\Delta \eta^\text{act}_\text{cat}>0, \quad \text{for} \quad y\in(h_1,h_\text{b}).
\end{equation}
Thus, the condition \eqref{eta_cat_inq} is not satisfied at any point of the active layer and:
\begin{equation}
\label{i_ct_cond_int}
\int_{h_1}^{h_\text{b}}i^{\text{ct(I)}}dy<\int_{h_1}^{h'_\text{b}}i^{\text{ct(II)}}dy.
\end{equation}
This leads us to a non-trivial observation. Namely, the estimation \eqref{i_ct_cond_int} does not mean than the solution $\text{(II)}$ produces higher cell current that solution $\text{(I)}$. It basically means that \textbf{solution $\text{(II)}$ does not exist}, because respective equations in the backing zone are not satisfied for the considered case if \eqref{i_ct_cond_int} holds true. As a result, for a predefined cell current and respective boundary conditions, the thickness of the active layer cannot be arbitrary and should be found as a component of solution.

\begin{remark}
\label{rem_neg_eta}
When the modification \eqref{eta_act_cat_abs} is neglected, the activation overpotential assumes negative values in the immediate vicinity of the interface between the backing and catalyst layers (compare subsection \ref{BV_comments}). Consequently, condition \eqref{eta_act_hb} may not be satisfied for some values of $h_\text{b}$. However, if the location of $h_\text{b}$ is taken out of the negative zone of $\eta_\text{cat}^\text{act(II)}$,  all of the above conclusions are true. Thus, even in the classical case, the catalyst layer thickness should be treated as an element of solution.
\end{remark}

\begin{remark}
\label{rem_oth_mod}
The above formal proof of uniqueness of solution with respect to the active layer thickness was derived for the Fickian diffusion model. However, respective balance equations in the active/backing layers are to be satisfied regardless of the diffusion model employed. As mentioned previously, the shaded areas in  Fig.\ref{co2_int} reflect the amount of consumed oxygen. From our analysis it follows that for the conditions given above they cannot be equal. Thus, for any other diffusion model, even if the shapes of respective curves are different from those obtained for the binary diffusion, the areas describing the consumed oxygen still need to be different from each other (as the overall mass balance has to be satisfied). Thus, we believe that for other diffusion models the active layer thickness should be a component of solution as well.
\end{remark}

\begin{figure}
\begin{center}
\includegraphics[scale=0.80]{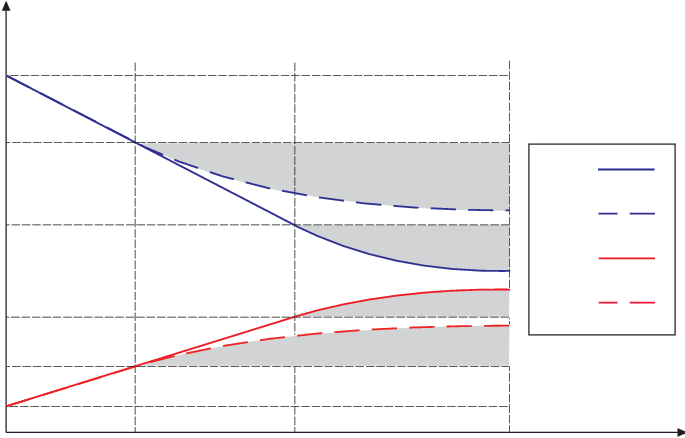}
\put(-56,100){$C_{\text{H}_2}^{\text{(I)}}$}
\put(-56,83){$C_{\text{H}_2}^{\text{(II)}}$}
\put(-56,63){$C_{\text{H}_2\text{O}}^{\text{(I)}}$}
\put(-56,46){$C_{\text{H}_2\text{O}}^{\text{(II)}}$}
\put(-258,155){$C_{\text{H}_2}, C_{\text{H}_2\text{O}}$}
\put(-290,112){$C_{\text{H}_2}^{\text{a(II)}}$}
\put(-290,80){$C_{\text{H}_2}^{\text{a(I)}}$}
\put(-290,138){$C_{\text{H}_2}^\text{bulk}$}
\put(-290,10){$C_{\text{H}_2\text{O}}^\text{bulk}$}
\put(-290,45){$C_{\text{H}_2\text{O}}^{\text{a(I)}}$}
\put(-290,25){$C_{\text{H}_2\text{O}}^{\text{a(II)}}$}
\put(-12,8){$y$}
\put(-270,-10){$-h_3$}
\put(-222,-10){$-h_\text{a}'$}
\put(-160,-10){$-h_\text{a}$}
\put(-78,-10){$-h_1$}

\caption{The schematic view of the mass fractions of hydrogen and water in the anode.  Solid lines refer to the variant when the catalyst layer is defined by $h_\text{a}$, while the dashed ones pertain to the active zone being delimited by  $h_\text{a}'$. The same magnitudes of $C_{\text{H}_2}^\text{bulk}$ and $C_{\text{H}_2\text{O}}^\text{bulk}$ are assumed. The shaded areas refer to the amount of hydrogen consumed and water produced in respective cases.}
\label{anoda_dowod}
\end{center}
\end{figure}
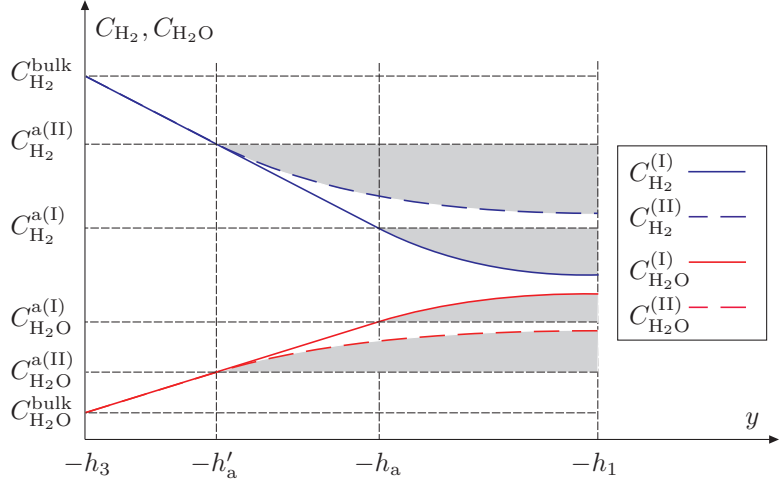

Now, let us conduct similar analysis for the anode. The schematic view of the compared solutions for the hydrogen and water mass fractions are shown in Fig.\ref{anoda_dowod}. As previously, we assume that solution $\text{(I)}$, obtained for the active layer bounded by coordinate $h_\text{a}$,  satisfies all the governing equations, boundary and transmission conditions. In the following we investigate whether the solution $\text{(II)}$, computed for the enlarged catalyst zone (defined by $h'_\text{a}>h_\text{a}$) and all the remaining parameters unchanged, can exist.

Clearly, the necessary condition to have the same cell current produced by both solutions is:
\begin{equation}
\label{i_cell_an}
\int_{-h_\text{a}}^{h_1}i^{\text{ct(I)}}_\text{ano}=\int_{-h'_\text{a}}^{h_1}i^{\text{ct(II)}}_\text{ano},
\end{equation}
which implies that:
\begin{equation}
\label{i_ct_ano_inq}
i^{\text{ct(I)}}_\text{ano}>i^{\text{ct(II)}}_\text{ano}
\end{equation}
over most of the interval $y\in(-h_\text{a},-h_1)$. When employing the same reasoning as that for the cathode, we have:
\begin{equation}
\label{C_I_II_an_eq}
C_{\text{H}_2}^{\text{(I)}}=C_{\text{H}_2}^{\text{(II)}},\quad C_{\text{H}_2\text{O}}^{\text{(I)}}=C_{\text{H}_2\text{O}}^{\text{(II)}}, \quad \text{for} \quad y \in(-h_3,-h'_\text{a}),
\end{equation}
\begin{equation}
\label{DC_eq_ano}
\frac{\text{d}C_{\text{H}_2}^{\text{(I)}}}{\text{d}y}\big|_{y=-h_\text{a}}=\frac{\text{d}C_{\text{H}_2}^{\text{(II)}}}{\text{d}y}\big|_{y=-h'_\text{a}}, \quad \frac{\text{d}C_{\text{H}_2\text{O}}^{\text{(I)}}}{\text{d}y}\big|_{y=-h_\text{a}}=\frac{\text{d}C_{\text{H}_2\text{O}}^{\text{(II)}}}{\text{d}y}\big|_{y=-h'_\text{a}},
\end{equation}
\begin{equation}
\label{C_I_II_inq_ano}
C_{\text{H}_2}^{\text{(I)}}<C_{\text{H}_2}^{\text{(II)}},\quad C_{\text{H}_2\text{O}}^{\text{(I)}}>C_{\text{H}_2\text{O}}^{\text{(II)}}, \quad \text{for} \quad y \in(-h_\text{a},-h_1).
\end{equation}

The charge transfer current is proportional to the following product (compare \eqref{i_ano_def} and \eqref{i_0_ano_def}):
\begin{equation}
\label{i_ano_pr}
i^\text{ct}_\text{ano} \sim \frac{C_{\text{H}_2\text{O}}^{0.4}}{C_{\text{H}_2}^{0.03}}\left[\exp\left(\frac{2F\eta^\text{act}_\text{ano}}{RT} \right)-\exp\left(-\frac{F\eta^\text{act}_\text{ano}}{RT} \right) \right].
\end{equation}
Note that this time the first component of the product has greater value for the variant $\text{(I)}$ than $\text{(II)}$:
\begin{equation}
\label{c_h2_h2o_inq}
 \frac{{C_{\text{H}_2\text{O}}^{\text(I)}}^{0.4}}{{C_{\text{H}_2}^{\text(I)}}^{0.03}}> \frac{{C_{\text{H}_2\text{O}}^{\text{(II)}}}^{0.4}}{{C_{\text{H}_2}^\text{(II)}}^{0.03}}, \quad \text{for} \quad y \in(-h_\text{a},-h_1).
\end{equation}
However, we start the analysis by investigating the change of activation overpotential defined as:
\begin{equation}
\label{delta_eta_ano}
\Delta \eta^\text{act}_\text{ano}=\eta^\text{act(II)}_\text{ano}-\eta^\text{act(I)}_\text{ano}=\phi_\text{el}^\text{(II)}-\phi_\text{el}^\text{(I)}+\phi_\text{ion}^\text{(I)}
-\phi_\text{ion}^\text{(II)}+\eta^\text{conc(I)}_\text{ano}-\eta^\text{conc(II)}_\text{ano}.
\end{equation}
From \eqref{eta_conc_ano} and \eqref{c_h2_h2o_inq} it follows that: $\eta^\text{conc(I)}_\text{ano}>\eta^\text{conc(II)}_\text{ano}$. When substituting \eqref{eta_conc_ano}, \eqref{fi_el_an_fin} and \eqref{fi_ion_an_fin} into \eqref{delta_eta_ano} we arrive at the following expression for $\Delta \eta^\text{act}_\text{ano}$:
\begin{equation}
\label{delta_eta_ano_fin}
\begin{split}
\Delta \eta^\text{act}_\text{ano}=\frac{2F\rho_\text{f}D_1}{M_{\text{H}_2}}\left(C_{\text{H}_2}^\text{(II)}(y)-C_{\text{H}_2}^\text{(I)}(y)\right)
\left(\frac{1}{\sigma_\text{ion}}+\frac{1}{\sigma_\text{el}}\right)+\frac{RT}{2F}\ln\left(\frac{C_{\text{H}_2}^\text{(II)}(y)}{C_{\text{H}_2}^\text{(I)}(y)}\frac{C_{\text{H}_2\text{O}}^\text{(I)}(y)}{C_{\text{H}_2\text{O}}^\text{(II)}(y)}\right)
\\
-\frac{j_\text{cell}}{\sigma_\text{el}}(h'_\text{a}-h_\text{a})+
\frac{RT}{2F}\ln\left(\frac{C_{\text{H}_2\text{O}}^\text{a(II)}}{C_{\text{H}_2\text{O}}^\text{a(I)}}\frac{C_{\text{H}_2}^\text{a(I)}}{C_{\text{H}_2}^\text{a(II)}}\right).
\end{split}
\end{equation}
In the above expression the first two terms are positive  monotonically increasing functions of $y$, while the remaining two yield a negative constant. Bearing in mind that for a positive activation overpotential:
\begin{equation}
\label{eta_act_ha}
\Delta \eta^\text{act}_\text{ano}(-h_\text{a})=\eta^\text{act(II)}_\text{ano}(-h_\text{a})>0,
\end{equation}
we come to the following conclusion:
\begin{equation}
\label{eta_act_ano_pos}
\Delta \eta^\text{act}_\text{ano}>0, \quad \text{for} \quad y\in(-h_\text{a},-h_1).
\end{equation}
As proved above, $\eta^\text{act(II)}_\text{ano}>\eta^\text{act(I)}_\text{ano}$ and thus, the second term of the product in the RHS of \eqref{i_ano_pr} is greater for solution $\text{(II)}$ (Remark \ref{rem_neg_eta} holds also in the case of anode). Now let us investigate if the increase in the value of this term can be counteracted by such a decline of the ratio ${C_{\text{H}_2\text{O}}}^{0.4}/{C_{\text{H}_2}}^{0.03}$ that keeps the cell current, $j_\text{cell}$, at an unchanged level (and in this way makes solution $\text{(II)}$ a permissible one). To this end, let us assume that the change of activation overpotential that results from the variation of the active layer thickness is small. In the following we accept the notation:
\begin{equation}
\label{S_def}
S(y)=\frac{i^\text{ct(II)}_\text{ano}}{i^\text{ct(I)}_\text{ano}}=\left(\frac{C^\text{(II)}_{\text{H}_2\text{O}}}{C^\text{(I)}_{\text{H}_2\text{O}}}\right)^{0.4}
\left(\frac{C^\text{(I)}_{\text{H}_2}}{C^\text{(II)}_{\text{H}_2}}\right)^{0.03}\frac{B^\text{(II)}}{B^\text{(I)}}=\left(\frac{C^\text{(II)}_{\text{H}_2\text{O}}}{C^\text{(I)}_{\text{H}_2\text{O}}}\frac{C^\text{(I)}_{\text{H}_2}}{C^\text{(II)}_{\text{H}_2}}\right)^{0.4} \left(\frac{C^\text{(II)}_{\text{H}_2}}{C^\text{(I)}_{\text{H}_2}}\right)^{0.37}\frac{B^\text{(II)}}{B^\text{(I)}},
\end{equation}
where:
\begin{equation}
\label{R_def}
B=\exp\left(\frac{2F\eta^\text{act}_\text{ano}}{RT} \right)-\exp\left(-\frac{F\eta^\text{act}_\text{ano}}{RT} \right).
\end{equation}
Obviously, the necessary condition to satisfy \eqref{i_ct_ano_inq} is that $S<1$ over the prevailing part of the active layer.

Remark  that in \eqref{S_def} only the first term of the product in the RHS  is less than 1, while the remaining part of the expression yields a value greater than 1.
For small magnitudes of $\Delta\eta^\text{act}_\text{ano}$ one has:
\begin{equation}
\begin{split}
\label{R_exp}
B^\text{(II)}=\exp\left(\frac{2F(\eta^\text{act(I)}_\text{ano}+\Delta\eta^\text{act}_\text{ano})}{RT} \right)-\exp\left(-\frac{F(\eta^\text{act(I)}_\text{ano}+\Delta\eta^\text{act}_\text{ano})}{RT} \right)\simeq\\
 \simeq B^\text{(I)}+
\frac{F}{RT}\left[2\exp\left(\frac{2F\eta^\text{act(I)}_\text{ano}}{RT} \right)+\exp\left(-\frac{F\eta^\text{act}_\text{ano}}{RT} \right)\right]\Delta\eta^\text{act}_\text{ano}.
\end{split}
\end{equation}
The following estimation holds for \eqref{R_exp}:
\begin{equation}
\label{R_est}
B^\text{(II)}\geq B^\text{(I)}+\frac{2F}{RT}B^\text{(I)}\Delta\eta^\text{act}_\text{ano}.
\end{equation}
In order to stay 'on the safe side' we take the lower bound of \eqref{R_est}, which yields:
\begin{equation}
\label{R_rat}
\frac{B^\text{(II)}}{B^\text{(I)}}=1+\frac{2F}{RT}\Delta\eta^\text{act}_\text{ano}.
\end{equation}

Next, from the formula for the activation overpotential \eqref{eta_act_ano} one has:
\begin{equation}
\label{c_rat_ax}
\begin{split}
\frac{C^\text{(II)}_{\text{H}_2\text{O}}}{C^\text{(I)}_{\text{H}_2\text{O}}}\frac{C^\text{(I)}_{\text{H}_2}}{C^\text{(II)}_{\text{H}_2}}=\exp\left(\frac{2F}{RT}\left(\phi_\text{el}^\text{(II)}-\phi_\text{el}^\text{(I)}+\phi_\text{ion}^\text{(I)}
-\phi_\text{ion}^\text{(II)}+\eta^\text{act(I)}_\text{ano}-\eta^\text{act(II)}_\text{ano}\right)\right)=\\
=\exp\left(-\frac{2F}{RT}\Delta\eta_{\text{ano}}^\text{act}\right)\exp\left(\frac{2F}{RT}L(y)\right),
\end{split}
\end{equation}
where
\begin{equation}
\label{L_def}
L(y)=\phi_\text{el}^\text{(II)}-\phi_\text{el}^\text{(I)}+\phi_\text{ion}^\text{(I)}
-\phi_\text{ion}^\text{(II)}.
\end{equation}
For small values of $\Delta\eta^\text{act}_\text{ano}$ one can write:
\begin{equation}
\label{rat_ax}
\left(\frac{C^\text{(II)}_{\text{H}_2\text{O}}}{C^\text{(I)}_{\text{H}_2\text{O}}}\frac{C^\text{(I)}_{\text{H}_2}}{C^\text{(II)}_{\text{H}_2}}\right)^{0.4}\simeq \left(1-\frac{0.8F}{RT}\Delta\eta^\text{act}_\text{ano} \right)\exp\left(\frac{0.8F}{RT}L(y)\right).
\end{equation}
After substitution of \eqref{R_rat} and \eqref{rat_ax} into \eqref{S_def} we obtain:
\begin{equation}
\label{S_fin_1}
S=\left(1+\frac{2F}{RT}\Delta\eta^\text{act}_\text{ano}\right)\left(1-\frac{0.8F}{RT}\Delta\eta^\text{act}_\text{ano} \right)\left(\frac{C^\text{(II)}_{\text{H}_2}}{C^\text{(I)}_{\text{H}_2}}\right)^{0.37}\exp\left(\frac{0.8F}{RT}L(y)\right),
\end{equation}
which can be rewritten with the accuracy order $O\left((\Delta\eta^\text{act}_\text{ano})^2\right)$ as:
\begin{equation}
\label{S_fin_2}
S=\left(1+\frac{1.2F}{RT}\Delta\eta^\text{act}_\text{ano}\right)\left(\frac{C^\text{(II)}_{\text{H}_2}}{C^\text{(I)}_{\text{H}_2}}\right)^{0.37}\exp\left(\frac{0.8F}{RT}L(y)\right).
\end{equation}
The first two terms of the above product are greater than 1, while the value of the third one depends on the sign of $L(y)$. By using formulae \eqref{fi_el_an_fin} and \eqref{fi_ion_an_fin} in \eqref{L_def} we arrive at an alternative definition of $L(t)$:
\begin{equation}
\label{L_alt}
L(t)=\frac{2F\rho_\text{f}D_1}{M_{\text{H}_2}}\left(C_{\text{H}_2}^\text{(II)}(y)-C_{\text{H}_2}^\text{(I)}(y)\right)\left(\frac{1}{\sigma_\text{ion}}+\frac{1}{\sigma_\text{el}}\right)-\frac{j_\text{cell}}{\sigma_\text{el}}(h'_\text{a}-h_\text{a})+\frac{RT}{2F}\ln\left(\frac{C_{\text{H}_2\text{O}}^\text{a(II)}}{C_{\text{H}_2\text{O}}^\text{a(I)}}\frac{C_{\text{H}_2}^\text{a(I)}}{C_{\text{H}_2}^\text{a(II)}}\right).
\end{equation}
The first term in the above expression is a positive increasing function of $y$, while the remaining two form a negative constant. Now, by substituting \eqref{L_alt} into \eqref{S_fin_2} and recalling that for positive activation overpotentials:
\begin{equation}
\label{S_ha}
S(h_\text{a})>1,
\end{equation}
we come to a conclusion that:
\begin{equation}
\label{S_fin_est}
S(y)>1, \quad \text{for} \quad y \in(-h_\text{a}-h_1).
\end{equation}
This effectively means that condition \eqref{i_ct_ano_inq} is not satisfied and solution $\text{(II)}$ does not exist as:
\begin{equation}
\label{i_cell_an_inq}
\int_{-h_\text{a}}^{h_1}i^{\text{ct(I)}}_\text{ano}<\int_{-h'_\text{a}}^{h_1}i^{\text{ct(II)}}_\text{ano}.
\end{equation}

Similarly as was done for the cathode, we have proved that for predefined conditions of the electrode operation (i.e. the boundary conditions and the cell current) there is only one active layer thickness for which the governing equations can be satisfied. Thus, the catalyst layer thickness should be a component of solution, unless some other parameters are relaxed. Also here Remarks \ref{rem_neg_eta}-\ref{rem_oth_mod} hold true.

The above observations have very serious ramifications for mathematical modelling of SOFCs. Note that the active layer thickness as described above is a parameter necessary to secure existence, uniqueness and proper mathematical structure of solution. As such, it does not have to coincide with the real physical catalyst layer. In fact, the more credible mathematical model is employed (e.g. the more advanced the diffusion model), the closer to each other the respective thicknesses of active layers are expected to be.

All the theoretical considerations carried out above led us to an essential modification of the original mathematical description of the SOFC problem. It includes: i) amendment of the electrochemical reactions models by means of formulae \eqref{eta_act_ano_abs}-\eqref{eta_act_cat_abs} (or alternatively \eqref{i_ct_gen_abs}) and ii) searching for the active layers thicknesses as elements of the solution (respective equation to find the active thickness for the cathode sub-problem will be specified later on). For this reason, we will henceforth call our formulation of the problem the 'modified formulation' (as opposed to the classical one).

\section{Numerical results}
\label{num_res}

In this section, the theoretical results derived above will be implemented in a numerical scheme to simulate the component problem for the cathode. A dedicated integral solver will be developed. The accuracy of computations will be verified against a newly introduced analytical benchmark solution for the combined physical fields. The problem of the active layer thickness will be analyzed further numerically. Finally, a comparison between the numerical simulations performed for the modified formulation and the experimental data will be given.

\subsection{The computational relations}
\label{comp_rels}

Below we derive a system of computational relations that will be used later on to construct the final numerical scheme. Note that, due to the considered model being isothermal, respective ODEs in the backing zone of the electrode can be solved analytically. However, if one wants to analyze the temperature-dependent version of the problem, respective integral relations can be easily derived and employed in the framework of the proposed algorithm.

When integrating respective ODEs in the backing zone (equations \eqref{phi_el_ODE}, \eqref{phi_ion_ODE}, \eqref{co2_ODE}) under the pertinent boundary conditions (\eqref{j_el_BC}$_2$, \eqref{j_ion_BC}$_2$, \eqref{fi_el_BC}$_2$, \eqref{C_BCs}$_3$, \eqref{J_ext_BCs}$_3$) one obtains the following solution:
\begin{equation}
\label{j_el_sol_back}
j_\text{el}(y)=j_\text{cell}, \quad y\in(h_\text{b},h_2),
\end{equation}
\begin{equation}
\label{j_ion_sol_back}
j_\text{ion}(y)=0, \quad y\in(h_\text{b},h_2),
\end{equation}
\begin{equation}
\label{J_o2_sol_back}
J_{\text{O}_2}(y)=-\frac{M_{\text{O}_2}j_\text{cell}}{4F}, \quad y\in(h_\text{b},h_2),
\end{equation}
\begin{equation}
\label{fi_el_sol_back}
\phi_\text{el}(y)=V_2+\frac{j_\text{cell}}{\sigma_\text{el}}(h_2-y),\quad y\in(h_\text{b},h_2),
\end{equation}
\begin{equation}
\label{fi_ion_sol_back}
\phi_\text{ion}(y)=\text{const}=\phi_\text{b},\quad y\in(h_\text{b},h_2),
\end{equation}
\begin{equation}
\label{C_o2_sol_back}
C_{\text{O}_2}(y)=C_{\text{O}_2}^\text{bulk}-\frac{M_{\text{O}_2}j_\text{cell}}{4F\rho_\text{a}D_2}(h_2-y), \quad y\in(h_\text{b},h_2),
\end{equation}
where constant $\phi_\text{b}$ is computed according to \eqref{fi_ion_b}.

Let us introduce the following function:
\begin{equation}
\label{lambda_def}
\Lambda(y)=\int_{h_1}^y i^\text{ct}_\text{cat}(\xi)d\xi,
\end{equation}
whose values and derivatives at the boundaries of the active layer are:
\begin{equation}
\label{lambda_prop}
\Lambda(h_1)=0,\quad \Lambda(h_2)=j_\text{cell}, \quad \frac{\text{d}\Lambda}{\text{d}y}\big|_{h_1}=i^\text{ct}_\text{cat}(h_1), \quad \frac{\text{d}\Lambda}{\text{d}y}\big|_{h_\text{b}}=0.
\end{equation}
From \eqref{phi_el_ODE}, \eqref{phi_ion_ODE}, \eqref{co2_ODE} and the conditions \eqref{J_BCs}$_3$, \eqref{j_el_ion_zero} one can prove that the solution for respective currents and oxygen flux in the catalyst zone has the form:
\begin{equation}
\label{j_el_sol_cat}
j_\text{el}(y)=\Lambda(y), \quad y\in(h_1,h_\text{b}),
\end{equation}
\begin{equation}
\label{j_ion_sol_cat}
j_\text{ion}(y)=j_\text{cell}-\Lambda(y), \quad y\in(h_1,h_\text{b}),
\end{equation}
\begin{equation}
\label{J_o2_sol_cat}
J_{\text{O}_2}(y)=-\frac{M_{\text{O}_2}}{4F}\Lambda(y), \quad y\in(h_1,h_\text{b}).
\end{equation}
Note that the above representation satisfies identically respective transmission conditions, as well as charge and mass balances.
By integration of \eqref{j_el_sol_cat}-\eqref{J_o2_sol_cat} from $y$ to $h_\text{b}$ we arrive at the following formulae for the potentials and mass fraction of oxygen:
\begin{equation}
\label{fi_el_sol_cat}
\phi_\text{el}(y)=V_\text{b}+\int_y^{h_\text{b}}\frac{\Lambda(\xi)}{\sigma_\text{el}}d\xi, \quad y\in(h_1,h_\text{b}),
\end{equation}
\begin{equation}
\label{fi_ion_sol_cat}
\phi_\text{ion}(y)=\phi_\text{b}+\int_y^{h_\text{b}}\frac{j_\text{cell}-\Lambda(\xi)}{\sigma_\text{ion}}d\xi, \quad y\in(h_1,h_\text{b}),
\end{equation}
\begin{equation}
\label{C_o2_sol_cat}
C_{\text{O}_2}(y)=C^\text{b}_{\text{O}_2}-\frac{M_{\text{O}_2}}{4F}\int_y^{h_\text{b}}\frac{\Lambda(\xi)}{\rho_\text{a}D_2}d\xi, \quad y\in(h_1,h_\text{b}),
\end{equation}
where $V_\text{b}$ and $C^\text{b}_{\text{O}_2}$ are given by \eqref{cat_ax_3} and \eqref{C^b_def}, correspondingly. Similarly as before with \eqref{j_el_sol_cat}-\eqref{J_o2_sol_cat}, respective transmission and balance conditions are satisfied automatically.

\subsection{The computational algorithm}
\label{comp_alg}

The proposed computational algorithm utilizes the above relations. Respective integrals in \eqref{lambda_def} and \eqref{fi_el_sol_cat}-\eqref{C_o2_sol_cat} are computed by cubic spline interpolation of the integrands and subsequent analytical integration. The following rescaling of the integration interval is introduced:
\begin{equation}
\label{z_def}
z=\frac{y-h_1}{\delta_\text{c}},
\end{equation}
due to which \eqref{lambda_def} can be converted to:
\begin{equation}
\label{lambda_def_z}
\Lambda(y)=\Lambda(\delta_\text{c}z+h_1)=\delta_\text{c}\int_{0}^z \tilde i^\text{ct}_\text{cat}(\xi)d\xi,
\end{equation}
where:
\begin{equation}
\label{i_ct_z}
\tilde  i^\text{ct}_\text{cat}(z)=i^\text{ct}_\text{cat}(\delta_\text{c}z+h_1).
\end{equation}
According to \eqref{lambda_prop}$_2$, equation \eqref{lambda_def_z} can be used to find the thickness of the catalyst layer:
\begin{equation}
\label{delta_c_obl}
\delta_\text{c}=\frac{j_\text{cell}}{\int_{0}^1 \tilde i^\text{ct}_\text{cat}(\xi)d\xi}.
\end{equation}
In our computations we assume the total cell current, $j_\text{cell}$, to be known and consequently the voltage drop is to be found as an element of solution. However, even in the reverse variant (the voltage drop imposed as a boundary condition) \eqref{delta_c_obl} should be employed.

The computational algorithm is composed of the following stages:
\begin{enumerate}
\item{At the preconditioning stage, the first approximations of respective dependent variables are defined in a way so as to preserve respective boundary and transmission conditions.}
\item{A set of dependent variables for $i$-th iteration is prepared: $\phi_\text{el}^{(i)}$, $\phi_\text{ion}^{(i)}$, $C_{\text{O}_2}^{(i)}$.}
\item{The activation overpotential, $\eta^\text{act}_\text{cat}$, and the charge transfer current, $\tilde i^\text{ct}_\text{cat}$, for the next iteration are computed from \eqref{eta_act_cat_abs} and \eqref{i_cat_def}, respectively.}
\item{The catalyst layer thickness, $\delta_\text{c}$, is calculated according to \eqref{delta_c_obl}. At this stage the global charge balance is satisfied for the previously accepted set of dependent variables. }
\item{The next approximation of the oxygen mass fraction, $C_{\text{O}_2}^{(i+1)}$, is obtained from \eqref{C_o2_sol_back} and \eqref{C_o2_sol_cat}.}
\item{The concentration overpotential, $\eta_\text{cat}^\text{conc}$, is updated with equation \eqref{eta_conc_cat}.}
\item{The electron and ion potentials, $\phi_\text{el}^{(i+1)}$ and $\phi_\text{ion}^{(i+1)}$, are computed from \eqref{fi_el_sol_back}, \eqref{fi_el_sol_cat} and \eqref{fi_ion_sol_back}, \eqref{fi_ion_sol_cat}, respectively.}
\item{The solution obtained in the iteration $(i+1)$ (points 3 -7) is compared with that from the previous step. If the difference does not exceed a prescribed value, then the computational process is terminated. Otherwise, we come back to the point 2, where  respective dependent variables are taken from points 5, 7 and the whole routine is repeated.  }
\end{enumerate}

Note that, when converged, the final solution satisfies the global and local charge and mass balances, as well as respective transmission conditions,  up to the level of accuracy of numerical integration. Pertinent boundary conditions are satisfied identically.

\subsection{Accuracy analysis}

The accuracy of computations provided by the above algorithm will be verified against an analytical benchmark example introduced in Appendix \ref{benchmark_app}. We set the following parameters of the benchmark solution: $h_\text{cat}=5\cdot 10^{-5}$ m, $j_\text{cell}= 10^5 \text{A}\cdot \text{m}^{-2}$, $V_2=0$ V, $T=700$$^\circ$C, $\alpha=8$. The material and electrochemical properties of the cathode are taken from \cite{Brus_et_al_2017}. The active layer thickness obtained under these conditions is $\delta_\text{c}=1.294\cdot 10^{-5} \text{m}$.

In the following, the solution accuracy will be described by the relative errors of: the potentials  of electron ($\Delta \phi_\text{el}$) and ion ($\Delta \phi_\text{ion}$) phases, the electron current ($\Delta j_\text{el}$ - note that the error of the ion potential can be recreated here from \eqref{j_cell}), the mass fraction of oxygen ($\Delta C_{\text{O}_2}$) and thickness of the active layer ($\Delta \delta_\text{c}$).

\begin{figure}[h!]

    \includegraphics [scale=0.40]{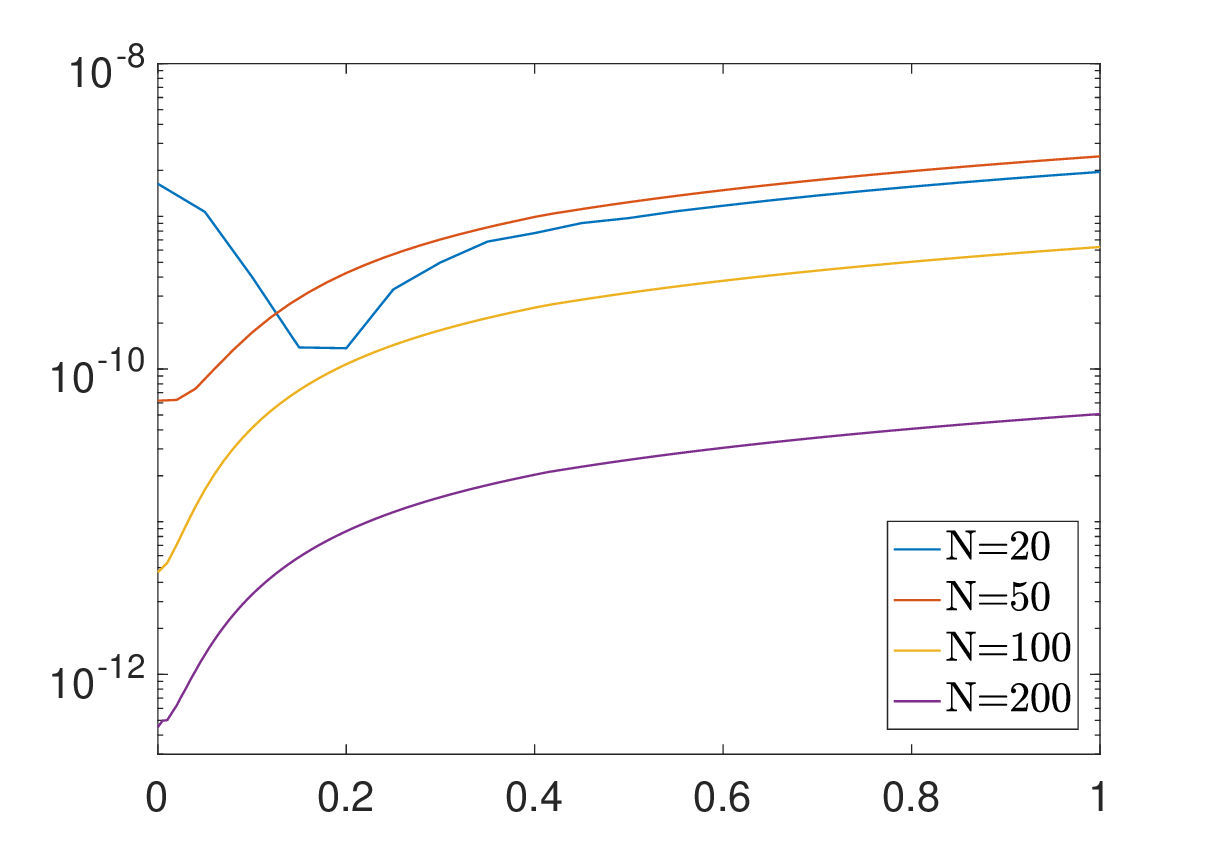}
    \put(-115,-5){$z$}
    \put(-240,75){$\Delta \phi_\text{el}$}
    \put(-230,160){$\textbf{a)}$}
    \hspace{2mm}
    \includegraphics [scale=0.40]{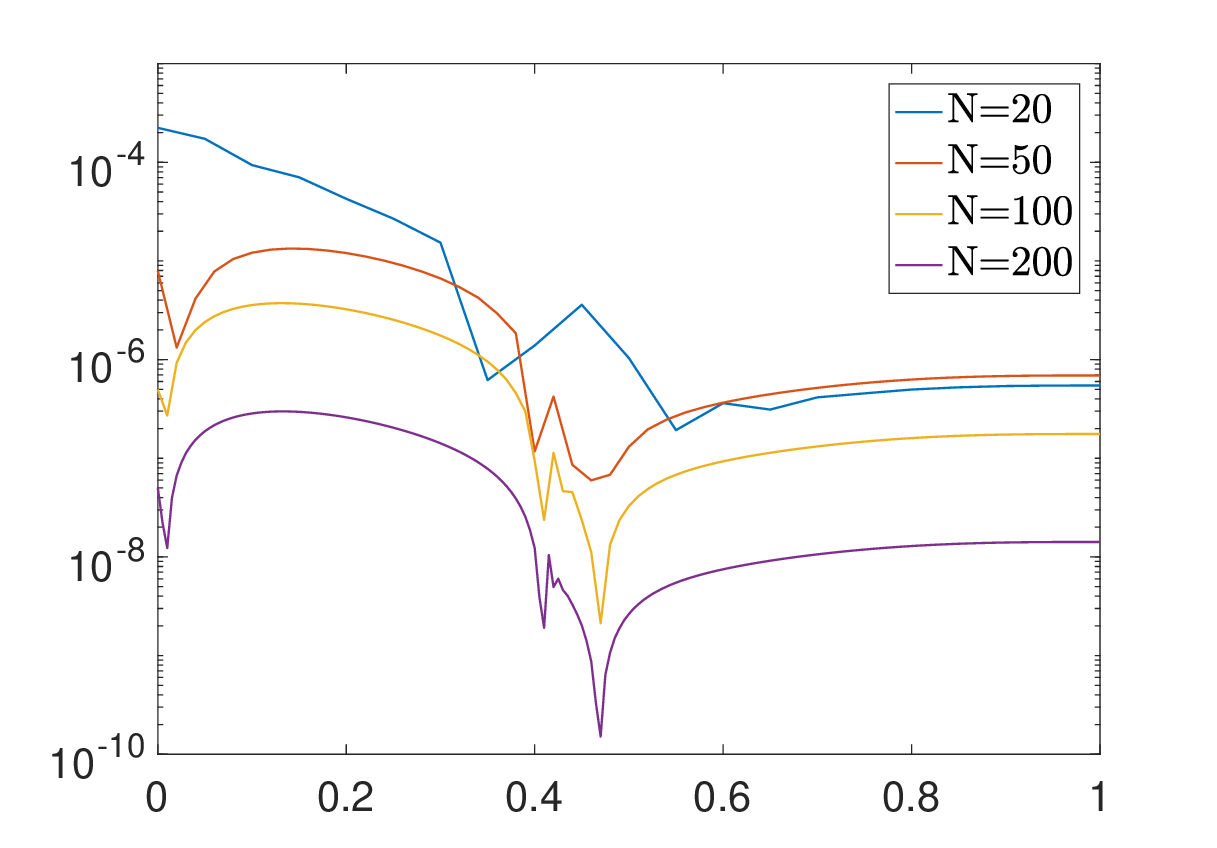}
    \put(-115,-5){$z$}
    \put(-240,75){$\Delta \phi_\text{ion}$}
    \put(-230,160){$\textbf{b)}$}

    \caption{Spatial distributions of relative errors for: a) the potential of electron phase, $\phi_\text{el} $, b) the potential of ion phase, $\phi_\text{ion} $, for different numbers of nodal points, $N$. }

\label{bledy_1}
\end{figure}

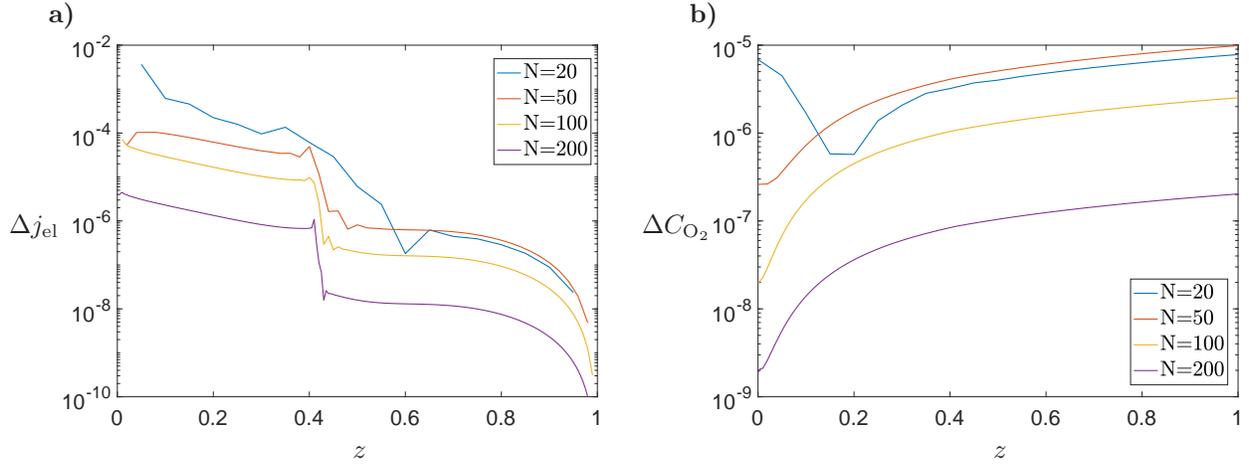
\begin{figure}[h!]

    \includegraphics [scale=0.40]{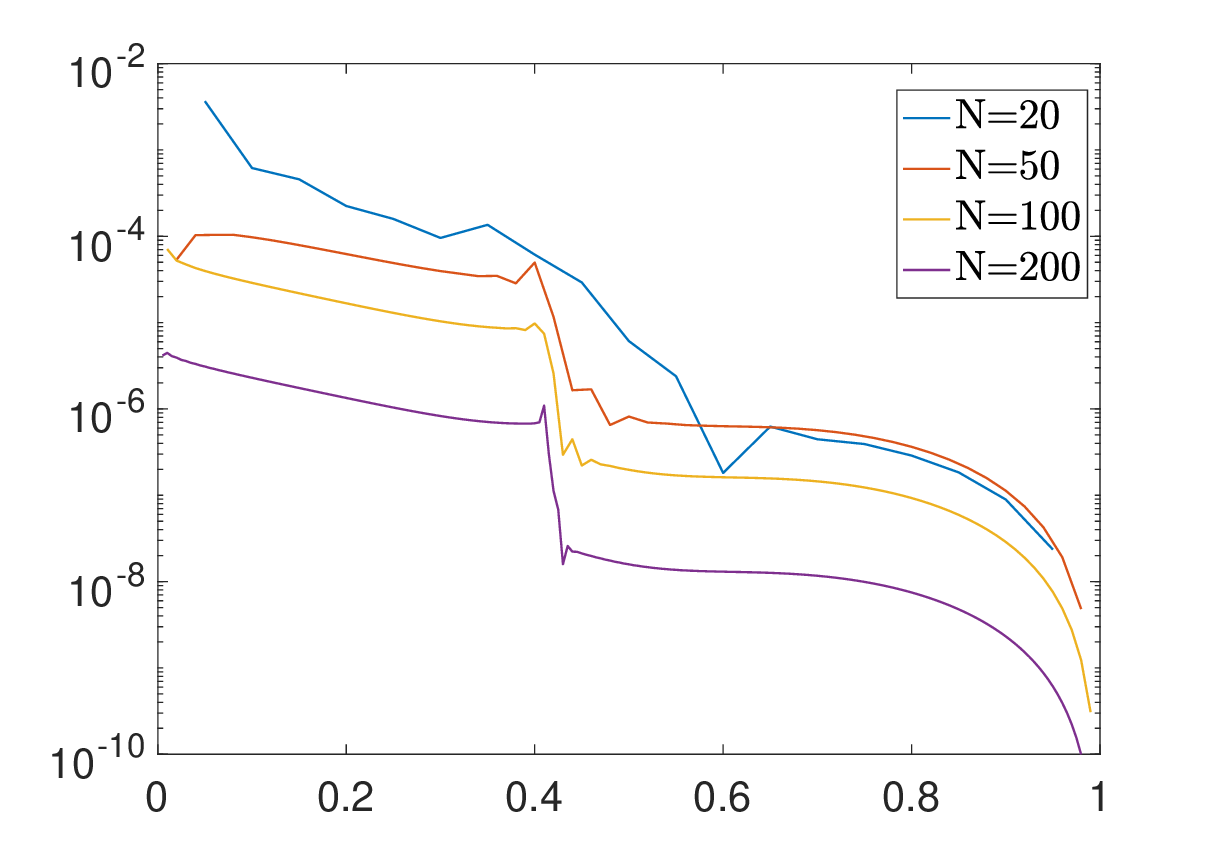}
    \put(-115,-5){$z$}
    \put(-245,80){$\Delta j_\text{el}$}
    \put(-230,160){$\textbf{a)}$}
    \hspace{2mm}
    \includegraphics [scale=0.40]{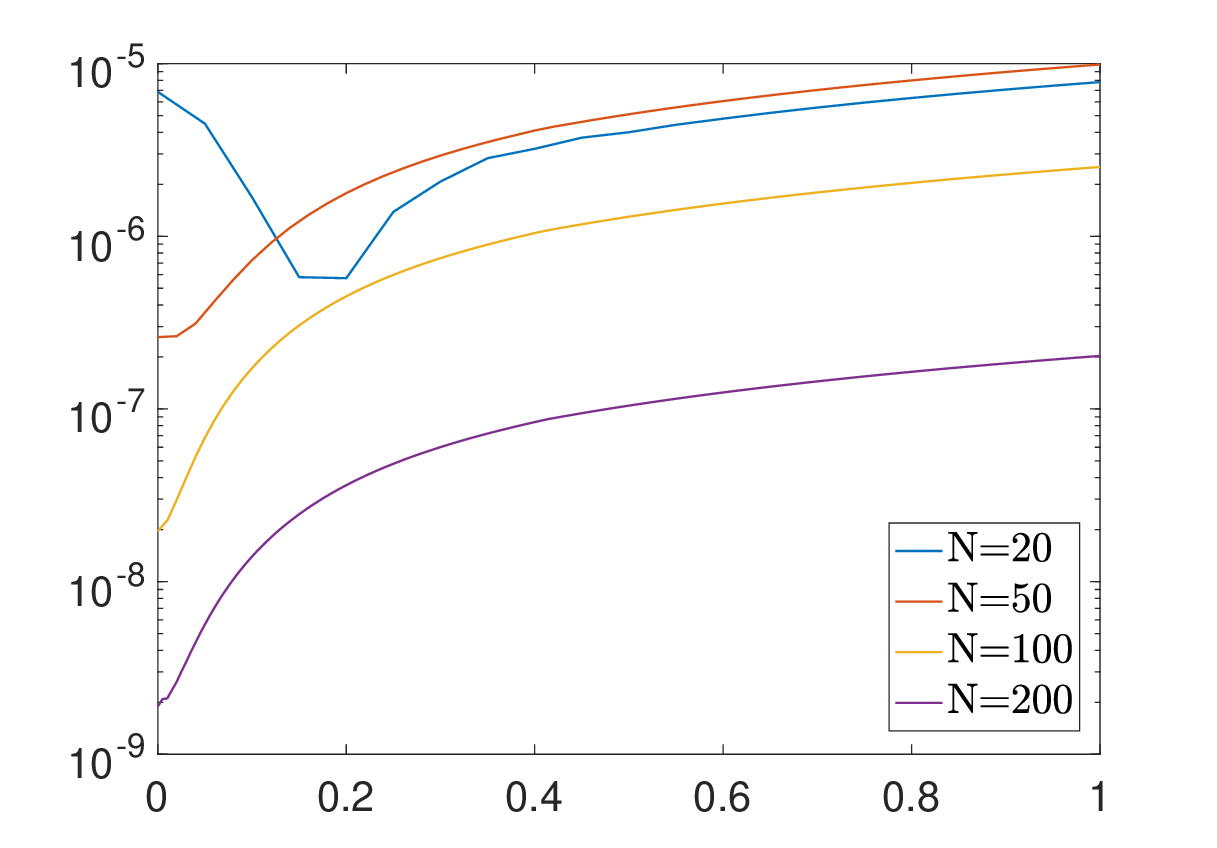}
     \put(-115,-5){$z$}
    \put(-248,80){$\Delta C_{\text{O}_2}$}
    \put(-230,160){$\textbf{b)}$}

    \caption{Spatial distributions of relative errors for: a) the electron current, $j_\text{el} $, b) the mass fraction of oxygen, $C_{\text{O}_2} $, for different numbers of nodal points, $N$.  }

\label{bledy_2}
\end{figure}

In Figs.\ref{bledy_1}-\ref{bledy_2} we show spatial distributions of the above errors for different numbers of nodal points $N$ ($N=\{20,50,100,200\}$) of a regular mesh. As can be seen, very good accuracy of the solution is obtained even with a mere twenty nodal points.
The accuracy of $\phi_\text{el}$ is much better than that of $\phi_\text{ion}$ despite using the same type of integral operator in computations (see \eqref{fi_el_sol_cat} and \eqref{fi_ion_sol_cat}). This is due to the fact that $\phi_\text{b}$ contains an additional error transferred from the component solution for the oxygen mass fraction, while for the isothermal model $V_\text{b}$ has a purely analytical value. Very low values of $\Delta j_\text{el}$ are obtained even in the immediate vicinity of the electrolyte interface, where the electron current turns to zero. A general trend of error reduction with growing number of nodal points is observed.

\begin{figure}[h!]

    \includegraphics [scale=0.40]{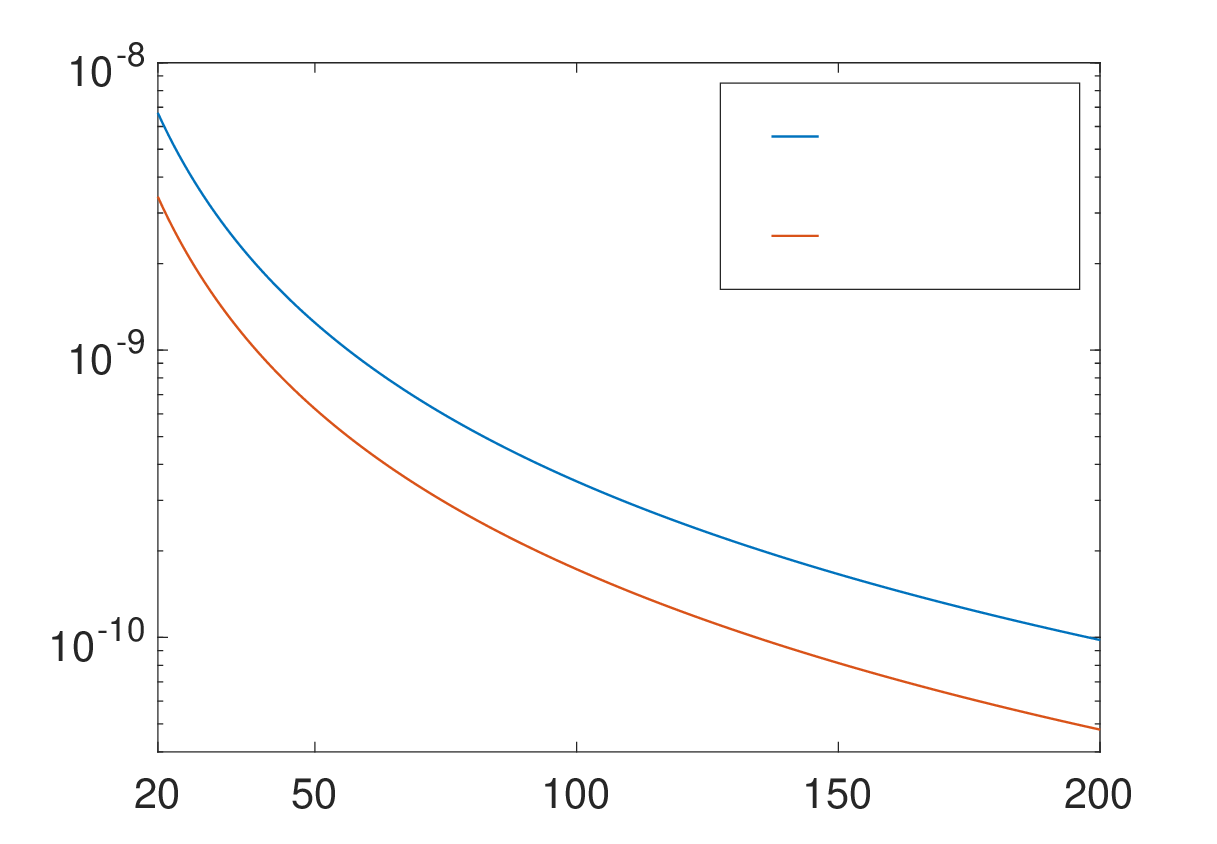}
    \put(-115,-5){$N$}
    \put(-245,80){$\Delta \phi_\text{el}$}
    \put(-72,134){$\Delta \phi_\text{el(max)}$}
    \put(-72,114){$\Delta \phi_\text{el(mean)}$}
    \put(-230,160){$\textbf{a)}$}
    \hspace{2mm}
    \includegraphics [scale=0.40]{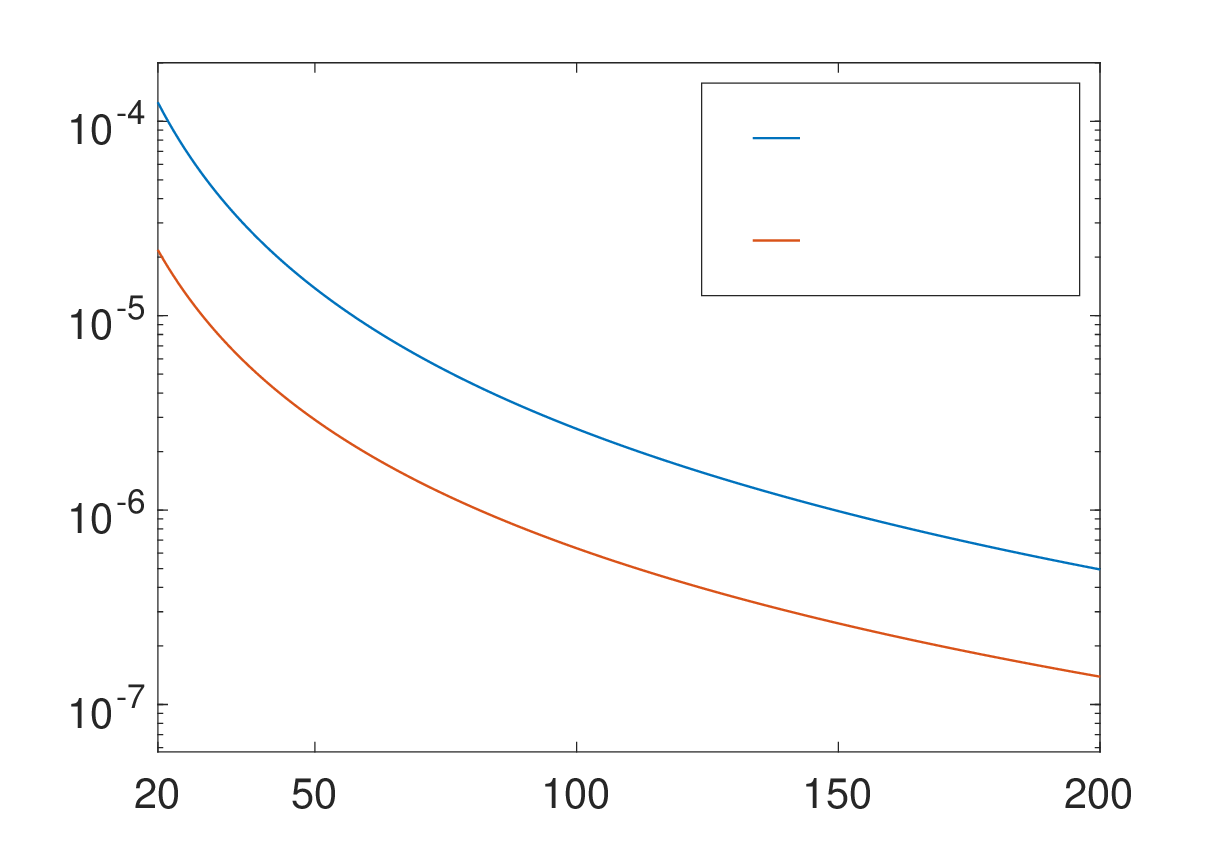}
     \put(-115,-5){$N$}
    \put(-248,80){$\Delta \phi_\text{ion}$}
    \put(-78,134){$\Delta \phi_\text{ion(max)}$}
    \put(-78,114){$\Delta \phi_\text{ion(mean)}$}
    \put(-230,160){$\textbf{b)}$}

    \caption{Maximal and mean relative errors for: a) the potential of electron phase, $\phi_\text{el} $, b) the potential of ion phase, $\phi_\text{ion} $, for different numbers of nodal points, $N$. }

\label{bledy_3}
\end{figure}

\begin{figure}[h!]

    \includegraphics [scale=0.40]{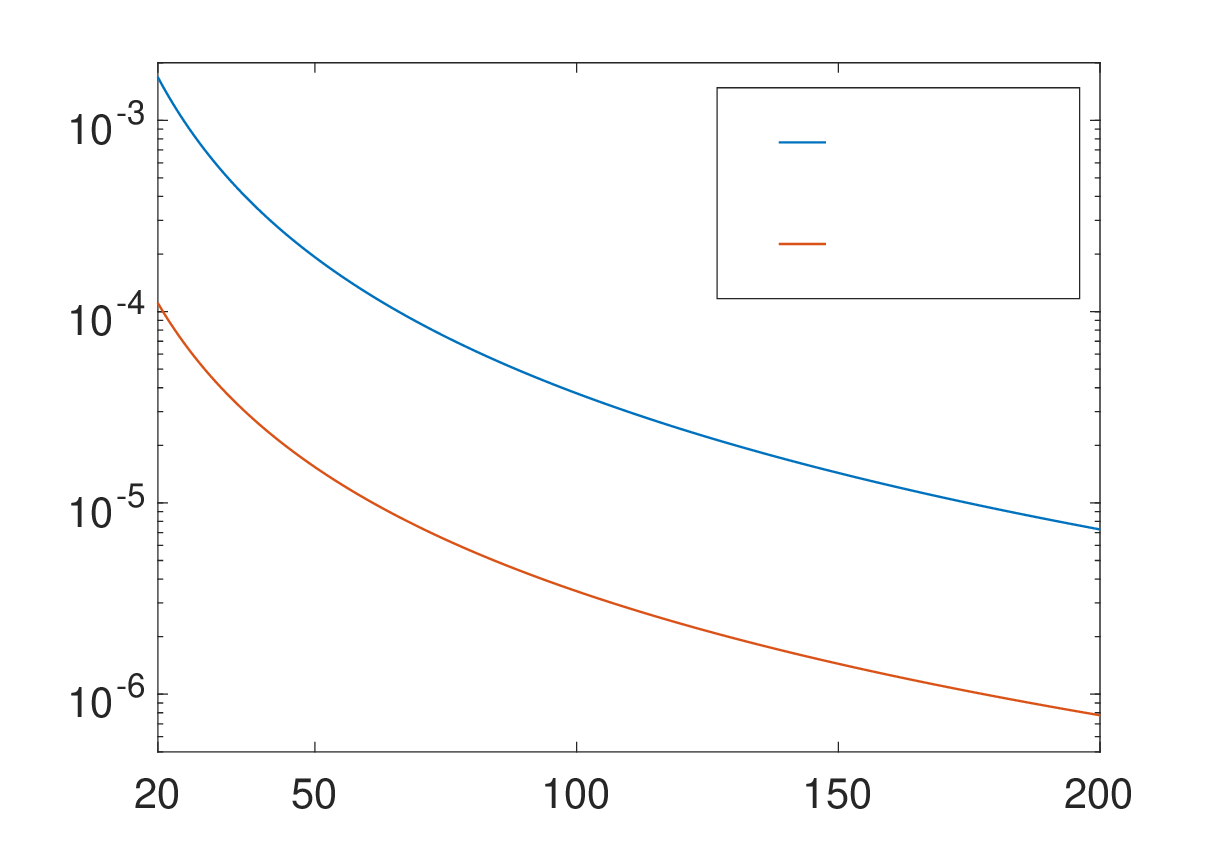}
    \put(-115,-5){$N$}
    \put(-245,80){$\Delta j_\text{el}$}
    \put(-72,132){$\Delta j_\text{el(max)}$}
    \put(-72,112){$\Delta j_\text{el(mean)}$}
    \put(-230,160){$\textbf{a)}$}
    \hspace{2mm}
    \includegraphics [scale=0.40]{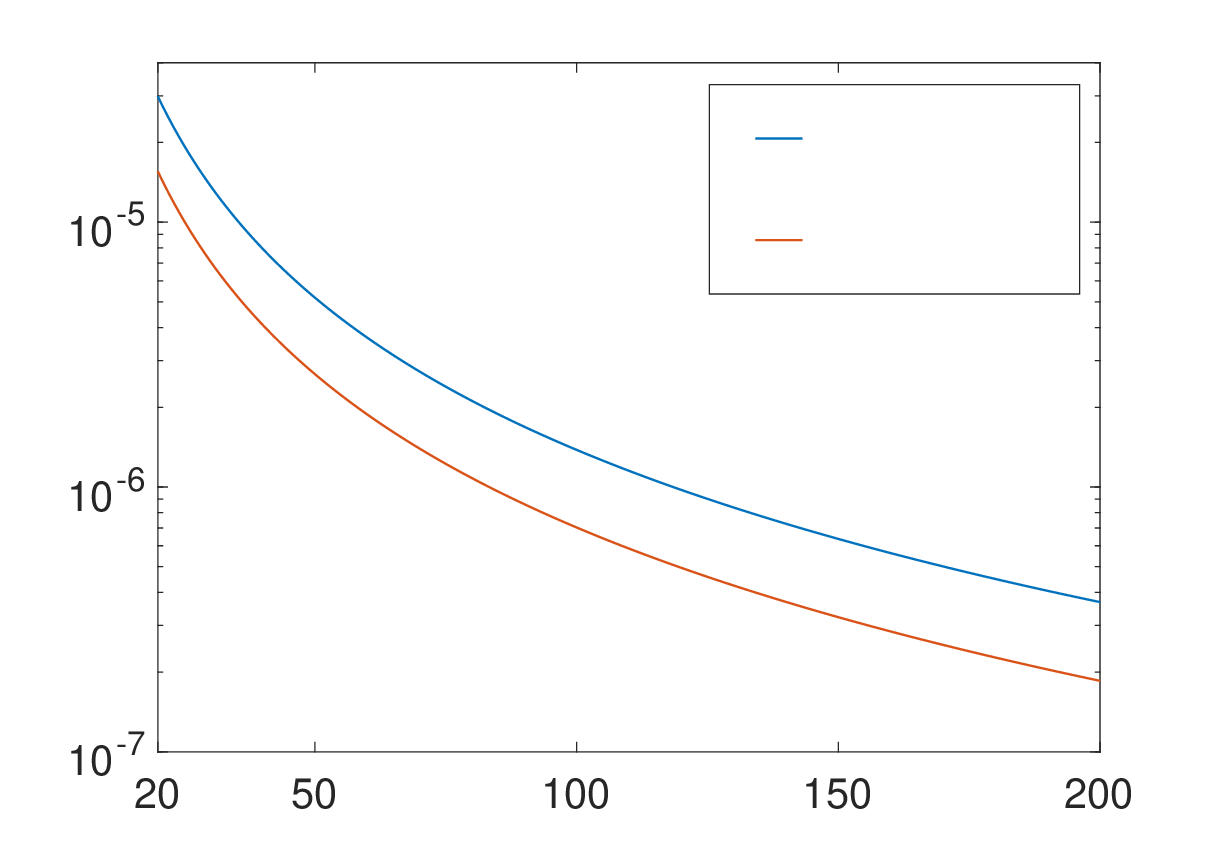}
     \put(-115,-5){$N$}
    \put(-248,80){$\Delta C_{\text{O}_2}$}
    \put(-78,134){$\Delta C_{\text{O}_2\text{(max)}}$}
    \put(-78,114){$\Delta C_{\text{O}_2\text{(mean)}}$}
    \put(-230,160){$\textbf{b)}$}

    \caption{Maximal and mean relative errors for: a) the electron current, $j_\text{el} $, b) the mass fraction of oxygen, $C_{\text{O}_2} $, for different numbers of nodal points, $N$.  }

\label{bledy_4}
\end{figure}

\begin{figure}
\begin{center}
\includegraphics[scale=0.40]{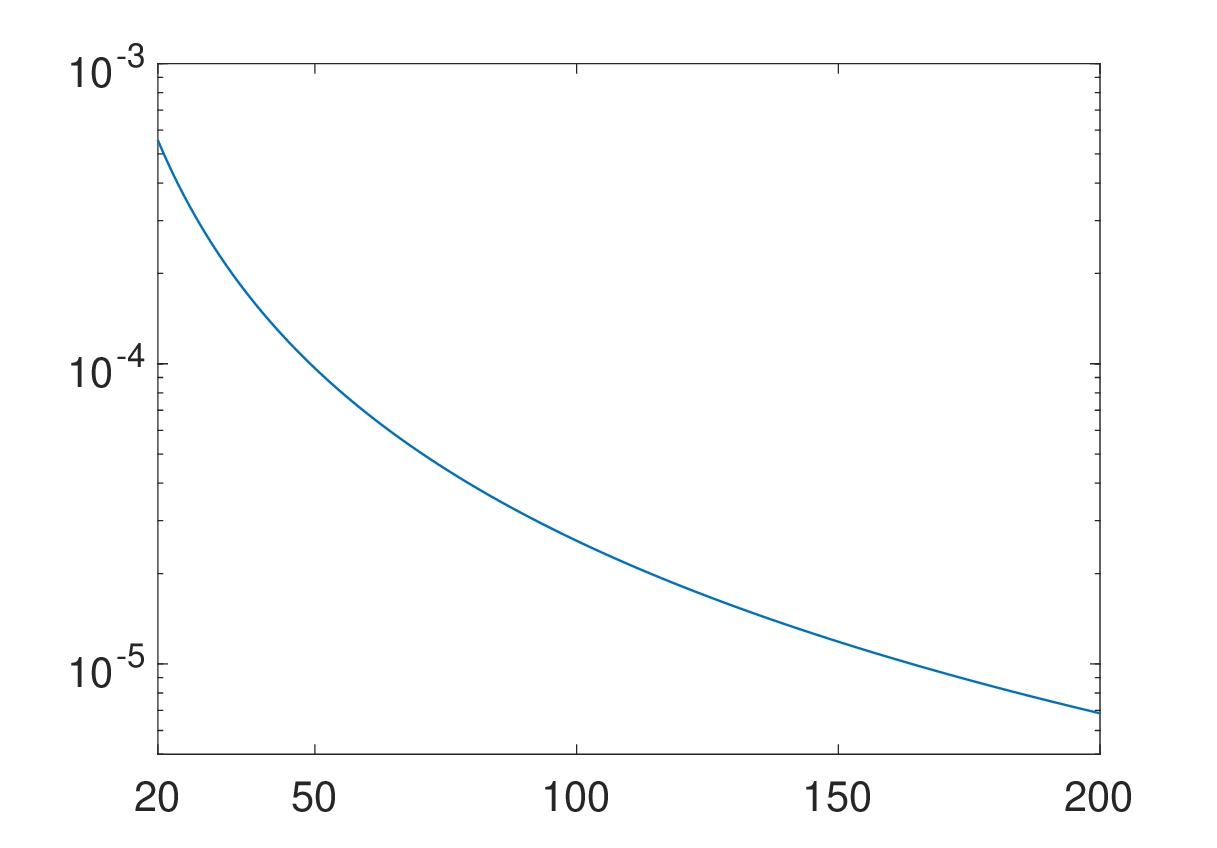}
\put(-240,80){$\Delta \delta_\text{c}$}
\put(-120,-5){$N$}

\caption{The relative error of the active layer thickness, $\delta_\text{c}$, for different numbers of nodal points, $N$.}
\label{bledy_5}
\end{center}
\end{figure}

In order to illustrate the latter tendency we depict the maximal and mean (over the active layer thickness) errors of solution as functions of the number of nodal points, $N$. The approximations of error curves are depicted in Figs. \ref{bledy_3}-\ref{bledy_5}.
It shows that, even for 200 nodal points, the error saturation is not achieved and accuracy can be further improved by increasing  mesh density. Depending on the considered solution component, the average error can be up to one order of magnitude lower than the maximal one. Note a very good accuracy to which the active layer thickness is determined. Even for 20 nodal points $\Delta \delta_\text{c}$ is of the order $10^{-4}$.

\subsection{LSCF cathode - numerical analysis}
\label{LSCF_num}

Having verified the credibility and accuracy of results produced by the developed solver, we will perform the numerical computations for the LSCF cathode whose geometrical and microstructural parameters are given in \cite{Brus_et_al_2017}. The problem of active layer thickness will be discussed based on this example.

In the analyzed problem the cathode thickness was $h_\text{cat}=5\cdot 10^{-5}$ m. The terminal voltage was set to $V_2=0.3$ V. Standard composition of air was assumed (by volume: 79$\%$ of nitrogen and 21$\%$ of oxygen). The computations were carried out for a mesh composed of $N=100$ uniformly spaced nodes. Below we present the results obtained for different values of the operating temperature, $T$, and cell current density, $j_\text{cell}$.

\begin{figure}[h!]

    \includegraphics [scale=0.40]{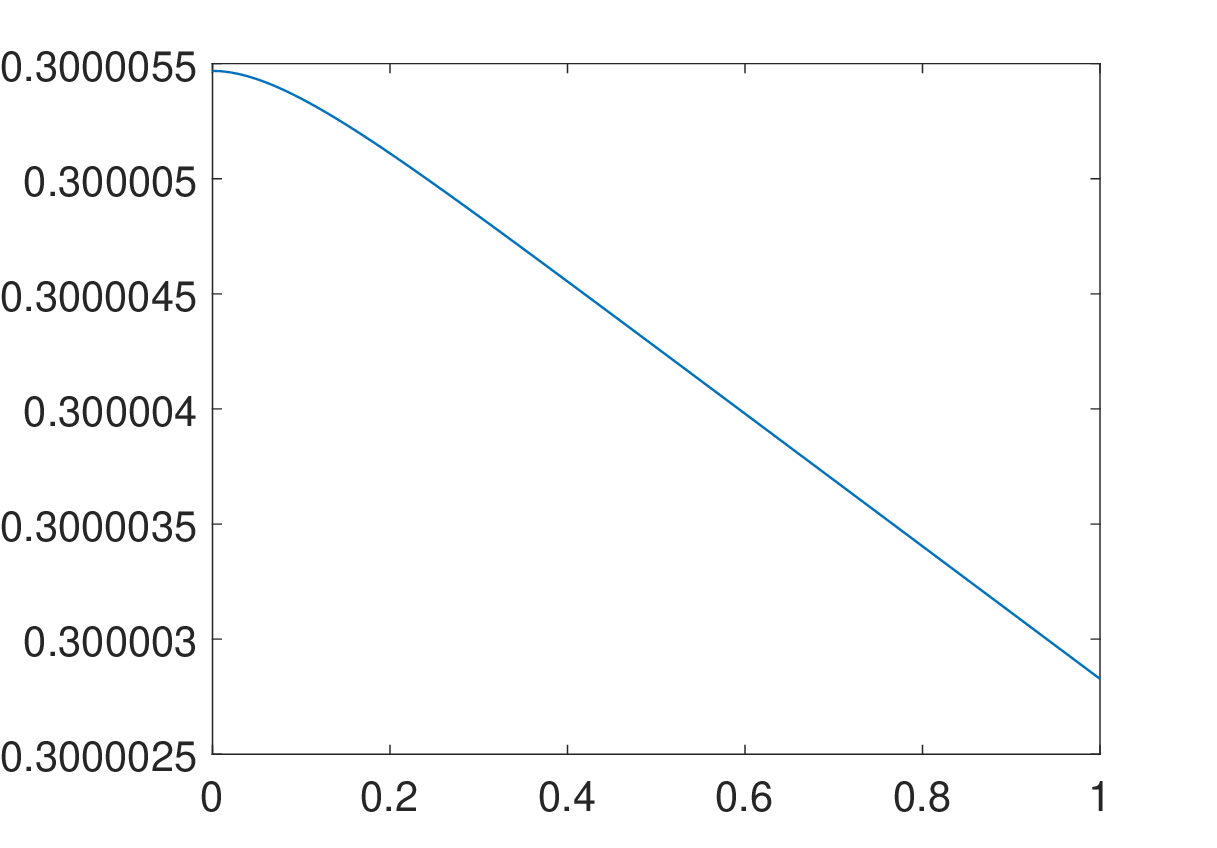}
    \put(-115,-5){$z$}
    \put(-245,80){$\phi_\text{el}$}
    \put(-230,160){$\textbf{a)}$}
    \hspace{2mm}
    \includegraphics [scale=0.40]{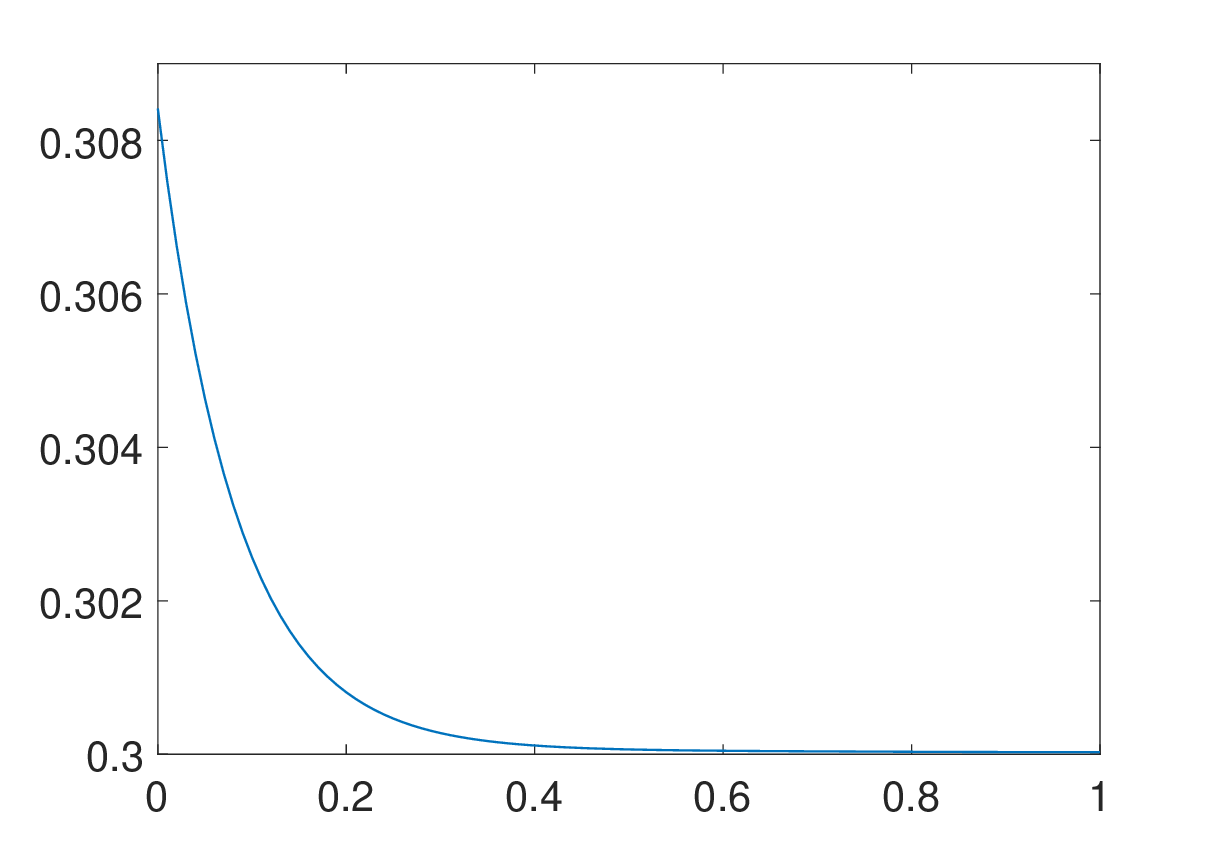}
     \put(-115,-5){$z$}
    \put(-248,80){$\phi_\text{ion}$}
    \put(-230,160){$\textbf{b)}$}

    \caption{Distributions of: a) the electron potential, $\phi_\text{el} $, b) the ion potential, $\phi_\text{ion}$, over the active layer thickness for $j_\text{cell}=2000$ $\text{A}\cdot \text{m}^{-2}$ and $T=800$ $^\circ \text{C}$.  }

\label{pot_LSCF}
\end{figure}

\begin{figure}[h!]

    \includegraphics [scale=0.40]{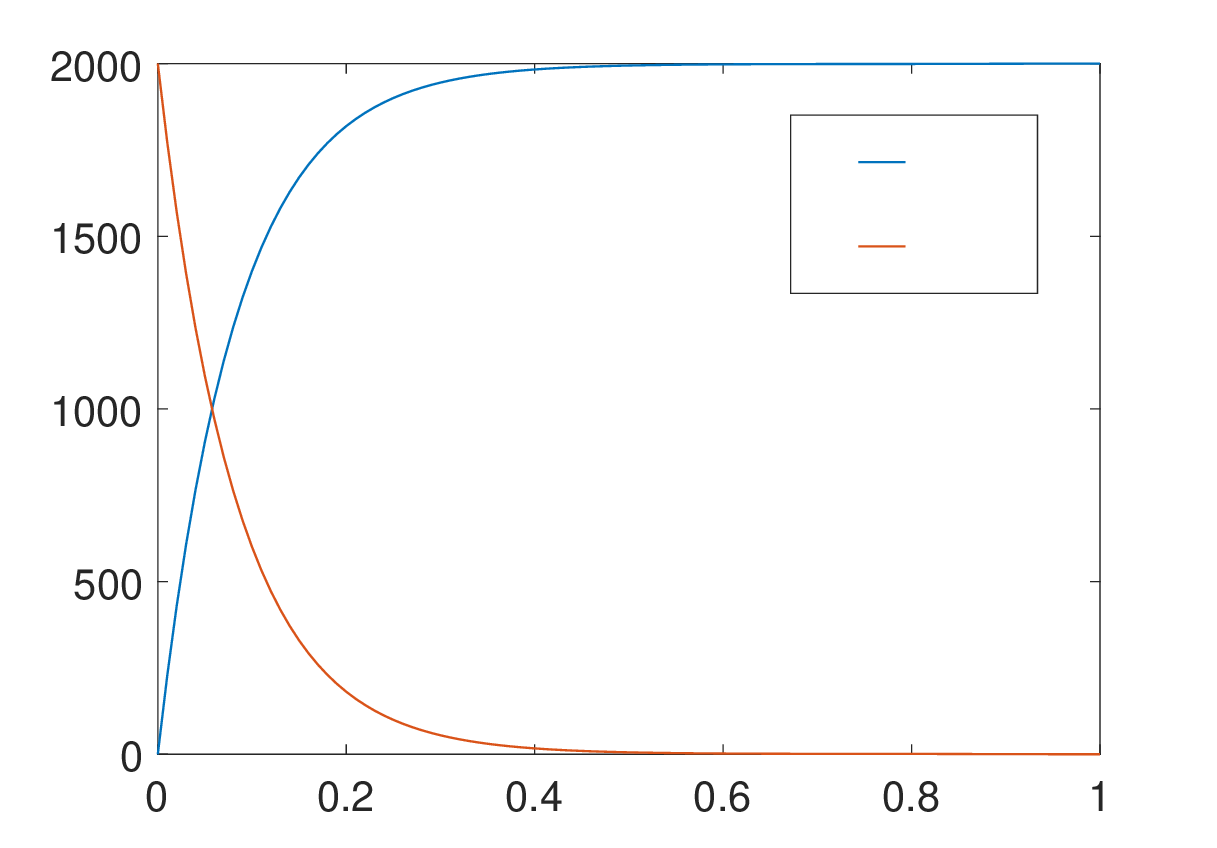}
    \put(-115,-5){$z$}
    \put(-55,130){$j_\text{el}$}
    \put(-55,115){$j_\text{ion}$}
    \put(-230,160){$\textbf{a)}$}
    \hspace{2mm}
    \includegraphics [scale=0.40]{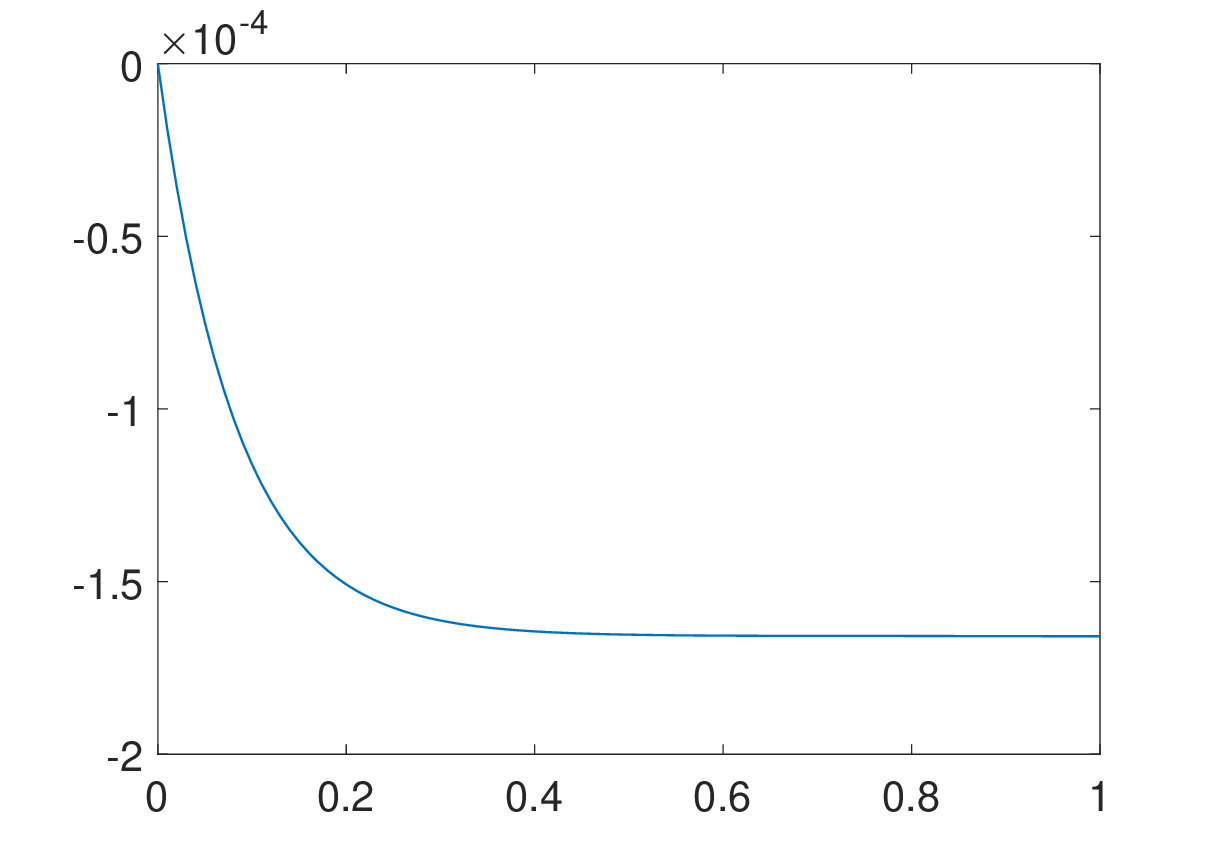}
     \put(-115,-5){$z$}
    \put(-232,80){$J_{\text{O}_2}$}
    \put(-230,160){$\textbf{b)}$}

    \caption{Distributions of: a) the electron and ion currents, $j_\text{el} $ and $j_\text{ion} $, b) the oxygen mass flux, $J_{\text{O}_2}$, over the active layer thickness for $j_\text{cell}=2000$ $\text{A}\cdot \text{m}^{-2}$ and $T=800$ $^\circ \text{C}$.  }

\label{j_J_LSCF}
\end{figure}

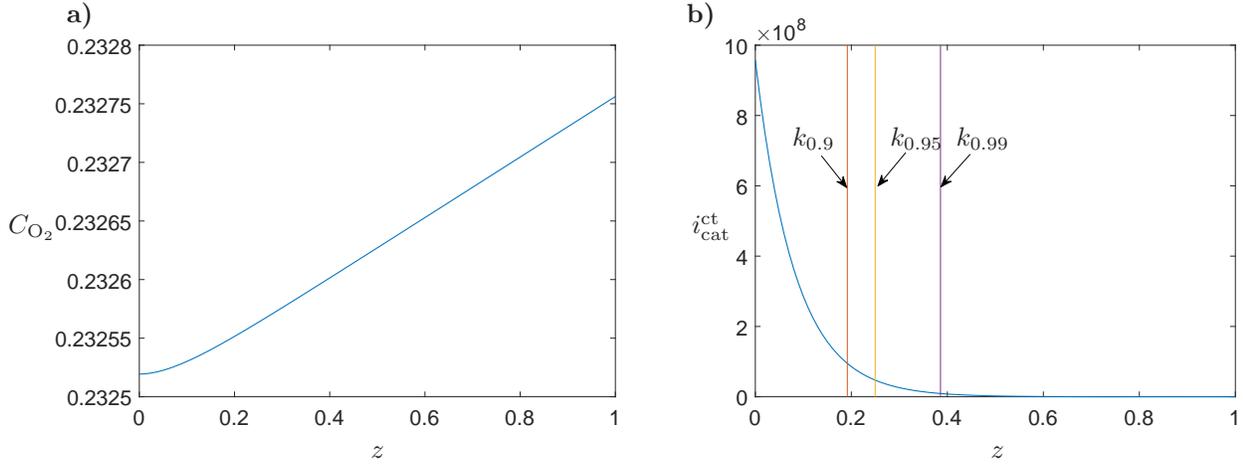
\begin{figure}[h!]

    \includegraphics [scale=0.40]{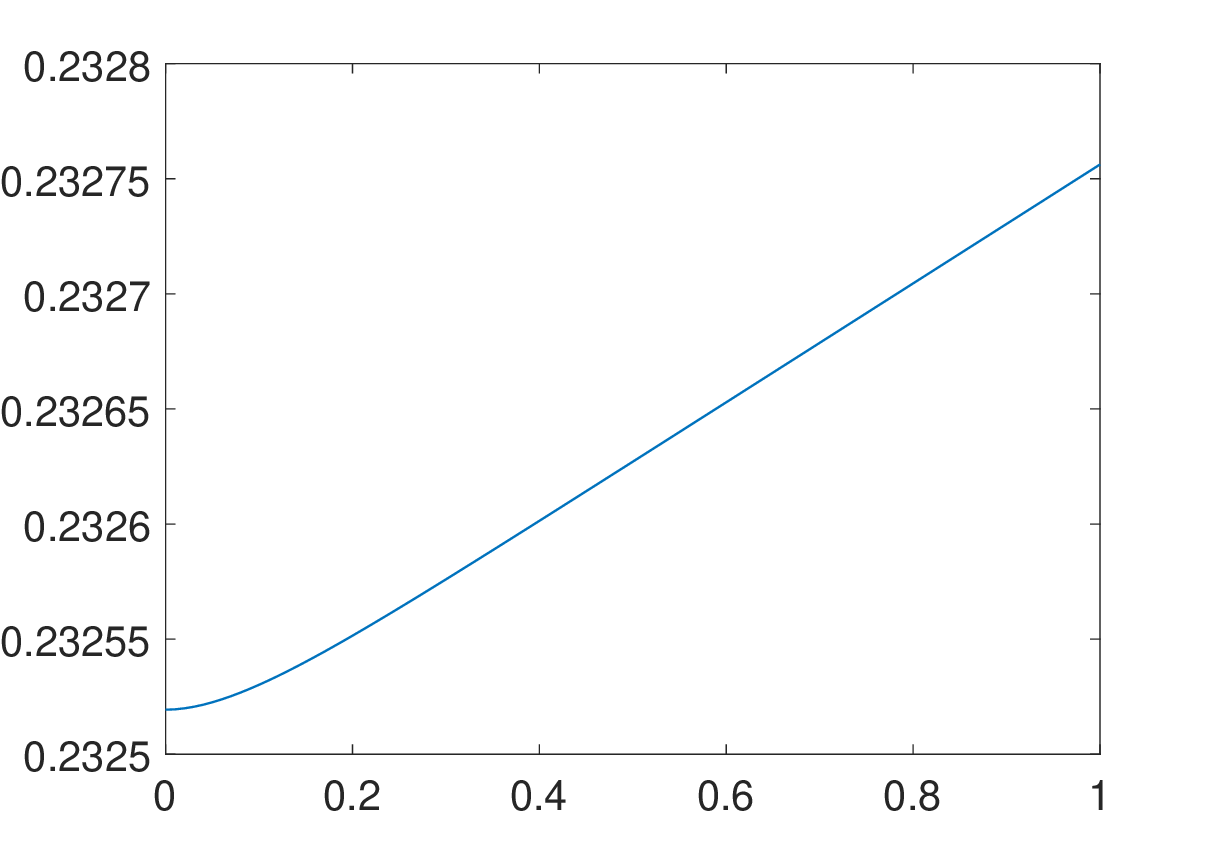}
    \put(-115,-5){$z$}
    \put(-252,80){$C_{\text{O}_2}$}
    \put(-230,160){$\textbf{a)}$}
    \includegraphics [scale=0.40]{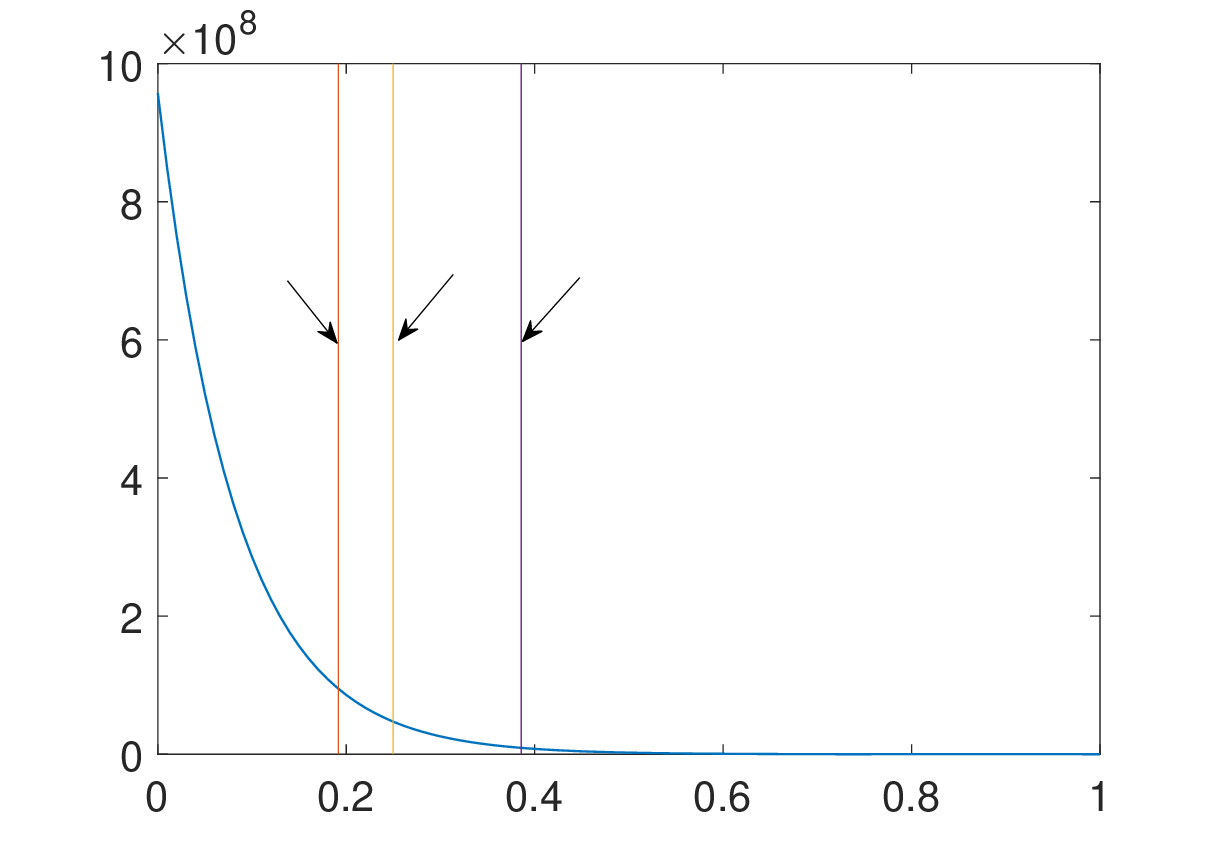}
     \put(-115,-5){$z$}
    \put(-228,80){$i_\text{cat}^\text{ct}$}
    \put(-190,113){$k_{0.9}$}
    \put(-153,113){$k_{0.95}$}
    \put(-128,113){$k_{0.99}$}
    \put(-230,160){$\textbf{b)}$}

    \caption{Distributions of: a) the mass fraction of oxygen, $C_{\text{O}_2} $, b) the charge transfer current, $i_\text{cat}^\text{ct}$, over the active layer thickness for $j_\text{cell}=2000$ $\text{A}\cdot \text{m}^{-2}$ and $T=800$ $^\circ \text{C}$. The vertical lines in figure b) refer to the section of the active layer over which 90$\%$, 95$\%$ or 99$\%$ of the current is produced, respectively.}

\label{co2_ict_LSCF}
\end{figure}

First, in Figs. \ref{pot_LSCF} -- \ref{co2_ict_LSCF}, we show spatial distributions of respective dependent variables for $T=800$$^{\circ}\text{C}$ and $j_\text{cell}=2000$ $ \text{A}\cdot \text{m}^{-2}$. For such operating conditions, the computed active layer thickness was $\delta_\text{c}=2.25\cdot 10^{-5}$ m. As can be seen in the figures, all the analyzed component physical fields are smooth, and no peculiarities of the type signalized in Figs. \ref{j_fi_2}--\ref{J_c_2} are reported (which is a consequence of employing the modification \eqref{eta_act_cat_abs}). Moreover, it shows that most of the electrochemical activity takes place only over some section of the catalyst zone. In Fig. \ref{co2_ict_LSCF}b) we have marked by vertical lines the sections of the catalyst layer which produce 90$\%$, 95$\%$ and 99$\%$ of the cell current, respectively. If we introduce the following definition for the thicknesses of these sections:
\begin{equation}
\label{delta_k_def}
\delta_j=k_j\delta_\text{c}, \quad j=\{0.9,0.95,0.99\},
\end{equation}
then the corresponding coefficients $k_j$ yield for the considered case (with the accuracy of two decimal digits):
\[
k_{0.9}=0.19, \quad k_{0.95}=0.25, \quad k_{0.99}=0.39.
\]
We believe that expression \eqref{delta_k_def} (possibly for $j=0.99$) can deliver better agreement with the real thickness of the active layer than the original definition accepted in this paper for the sake of computations.

\begin{figure}[h!]

    \includegraphics [scale=0.40]{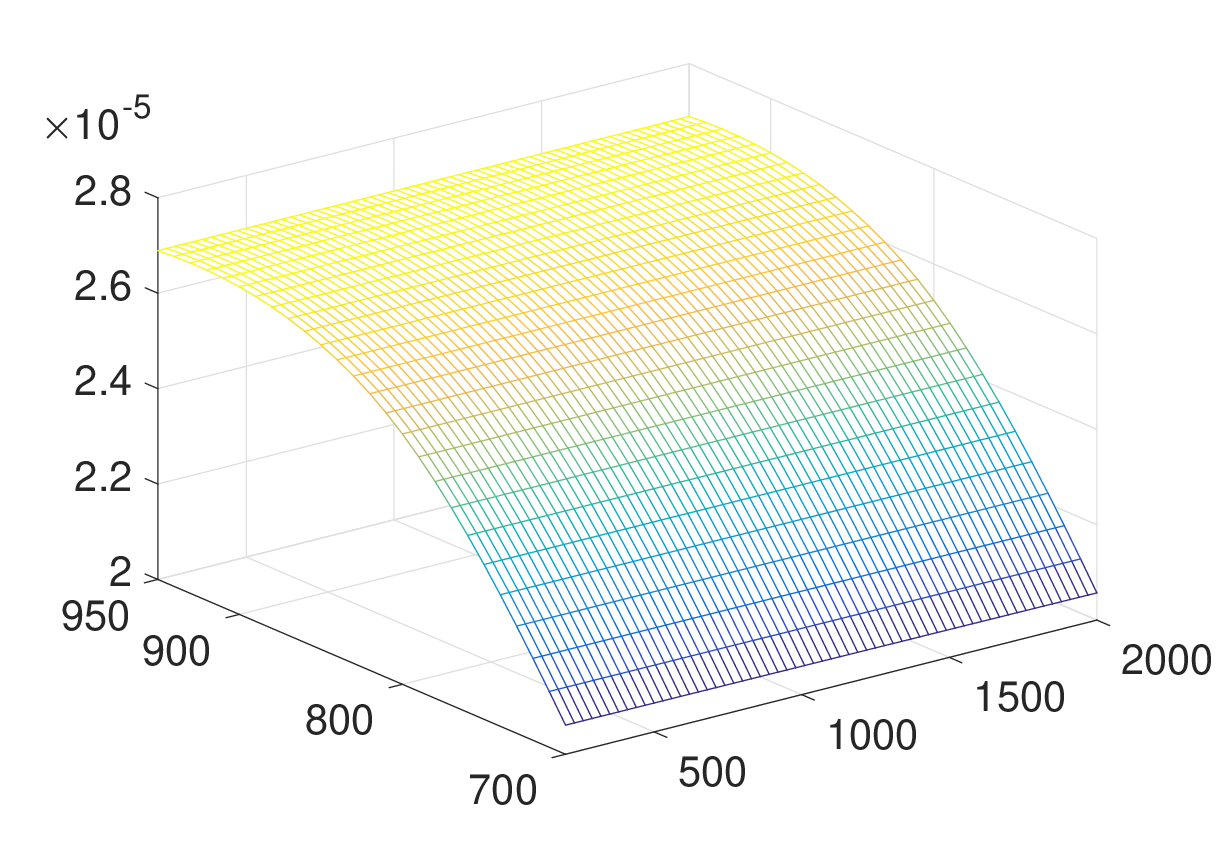}
    \put(-60,10){$j_\text{cell}$}
    \put(-190,15){$T$}
    \put(-235,90){$\delta_{\text{c}}$}
    \put(-230,160){$\textbf{a)}$}
    \hspace{2mm}
    \includegraphics [scale=0.40]{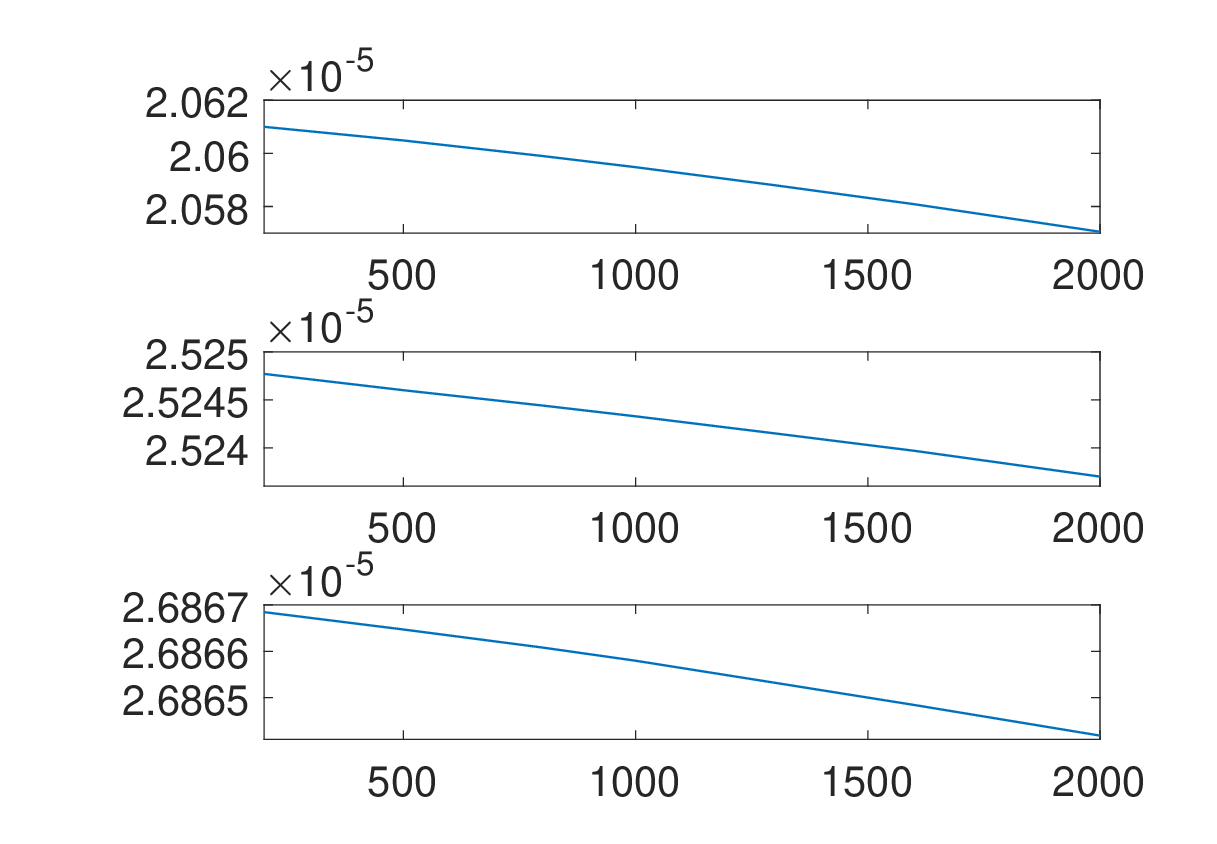}
     \put(-115,-5){$j_\text{cell}$}
     \put(-225,35){$\delta_\text{c}$}
     \put(-225,83){$\delta_\text{c}$}
     \put(-225,130){$\delta_\text{c}$}
     \put(-75,133){$T=700^\circ\text{C}$}
     \put(-75,85){$T=800^\circ\text{C}$}
     \put(-75,36){$T=900^\circ\text{C}$}
    \put(-230,160){$\textbf{b)}$}

    \caption{The thickness of active layer, $\delta_\text{c}$, as a function of the cell current density, $j_\text{cell}$, and the operating temperature, $T$. Figure b) depicts respective trends for constant temperature.}

\label{delta_i_T}
\end{figure}
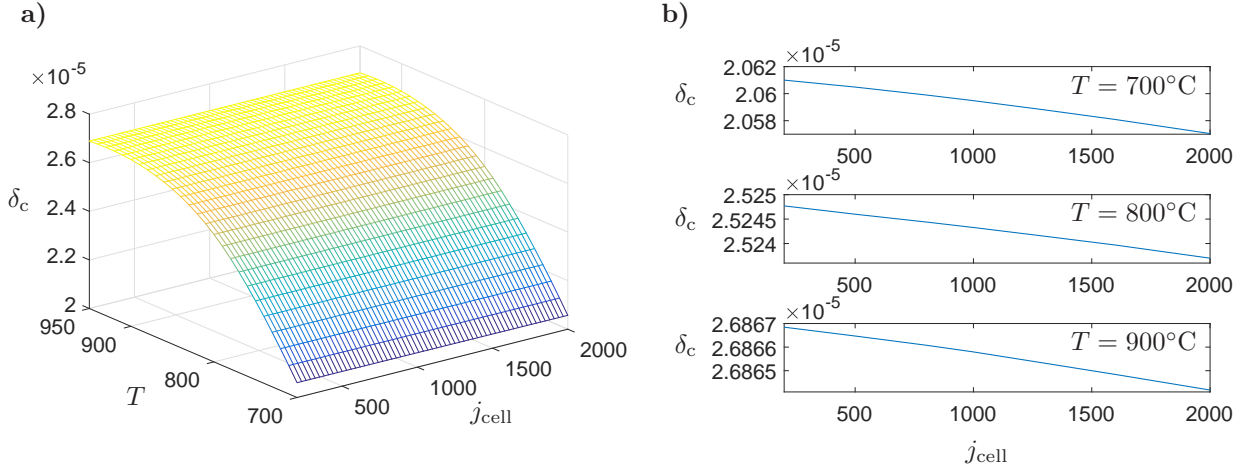

In order to complete this part of our analysis, we investigate the influence of the overall cell current density and the operating temperature on the thickness of the active layer. To this end, a series of simulations was run for different values of  $j_\text{cell}$, ranging from 200 $\text{A}\cdot \text{m}^{-2}$ to 2000 $\text{A}\cdot \text{m}^{-2}$, and the temperature, $T$, varying between 700 $^\circ \text{C}$ and 950 $^\circ \text{C}$. Respective results are shown in Fig.\ref{delta_i_T} a) (a 3D plot) and Fig.\ref{delta_i_T} b) ($\delta_\text{c}(j_\text{cell})$ for three different values of $T$). As can be seen, the thickness of the active layer is much more sensitive to the changes of the operating temperature than to the variation of the cell current. For constant $j_\text{cell}$ one can expect 30$\%$ increase in $\delta_\text{c}$ when changing the temperature from  $700^\circ \text{C}$ to $950^\circ \text{C}$. The respective impact of the cell current is much less pronounced, being reduced with the temperature magnification. However, a distinct decreasing trend is visible regardless of the value of $T$ (see Fig.\ref{delta_i_T} b).

\subsection{Active layer thickness - numerical verification of the uniqueness of solution. Elements of sensitivity analysis.}

In subsection \ref{cat_sub} it has been proved analytically that the thickness of the active layer is an element of solution, and as such cannot be predefined in an arbitrary way. As mentioned previously, the fixed values of $h_\text{a}$ or $h_\text{b}$ can be introduced only if some other conditions are relaxed. In the following, we will verify this conclusion numerically. 	

Let us consider the problem whose solution is depicted in Figs.\ref{pot_LSCF}-\ref{co2_ict_LSCF}. However, this time we shall modify the problem formulation by assuming that the thickness of the active layer, $\delta_\text{c}$, is known (and taken from the recalled numerical solution), while the mass fraction of oxygen at the external interface of the cathode, $C_{\text{O}_2}^\text{bulk}$, constitutes an element of solution to be found numerically. Obviously, according to the analysis provided in subsection \ref{cat_sub}, both solutions (i.e. the one from the subsection \ref{cat_sub} and that computed for the modified formulation) have to be equivalent.

In order to solve the new system of equations we need to modify the computational algorithm  described in subsection \ref{comp_alg}. To this end, we introduce a function of reduced mass fraction of oxygen defined as:
\begin{equation}
\label{co2_red}
\tilde C_{\text{O}_2}=\frac{C_{\text{O}_2}}{C_{\text{O}_2}^\text{bulk}}.
\end{equation}
Note that, when using the above function, the boundary condition \eqref{C_BCs}$_3$ converts to:
\begin{equation}
\label{co2_red_cond}
\tilde C_{\text{O}_2}(h_2)=1,
\end{equation}
while respective boundary and transmission conditions that employ the spatial derivative of $C_{\text{O}_2}$ should be rescaled by the factor $C_{\text{O}_2}^\text{bulk}$.

The expression for concentration overpotential \eqref{eta_conc_cat} yields now:
\begin{equation}
\label{eta_conc_new}
\eta_\text{cat}^\text{conc}=-\frac{RT}{4F}\ln\left(\tilde C_{\text{O}_2}\right).
\end{equation}
When combining \eqref{i_0_cat_def} with \eqref{i_cat_def} and \eqref{co2_red} one arrives at the following formula for the charge transfer current:
\begin{equation}
\label{ict_alt}
i^\text{ct}_\text{cat}=\left(C_{\text{O}_2}^\text{bulk}\right)^{0.2}\psi(y),
\end{equation}
where:
\begin{equation}
\label{psi_def}
\psi(y)=1.47\cdot 10^6 \cdot \left(\tilde C_{\text{O}_2}\frac{\rho_\text{a}RT}{M_{\text{O}_2}}\right)^{0.2}\exp \left(-\frac{85859}{RT}\right)A^\text{dpb}_\text{cat}\left[\exp\left(\frac{1.2F\eta^\text{act}_\text{cat}}{RT} \right)-\exp\left(-\frac{F\eta^\text{act}_\text{cat}}{RT} \right) \right].
\end{equation}
Next, by integrating equation \eqref{ict_alt} over the active layer thickness and taking into account relation \eqref{j_bal_int} we obtain after simple algebra an expression to compute $C_{\text{O}_2}^\text{bulk}$:
\begin{equation}
\label{co2_comp}
C_{\text{O}_2}^\text{bulk}=\left( \frac{j_\text{cell}}{\int_{h_1}^{h_\text{b}}\psi(y)dy}\right)^5.
\end{equation}

Having the above relations we are able to modify accordingly the computational algorithm from subsection \ref{comp_alg}. Respective changes include:
\begin{enumerate}
\item{Using a constant, predefined value of the active layer thickness, $\delta_\text{c}$.}
\item{Substituting the oxygen mass fraction function, $C_{\text{O}_2}$, by the reduced mass fraction of oxygen, $\tilde C_{\text{O}_2}$ \eqref{co2_red}.}
\item{Employing relation \eqref{co2_comp} to compute $C_{\text{O}_2}^\text{bulk}$. This point replaces point 4 of the original algorithm. Note that, as previously, at this stage of the scheme the global charge balance is identically satisfied. }
\end{enumerate}

\begin{figure}[h!]

    \includegraphics [scale=0.40]{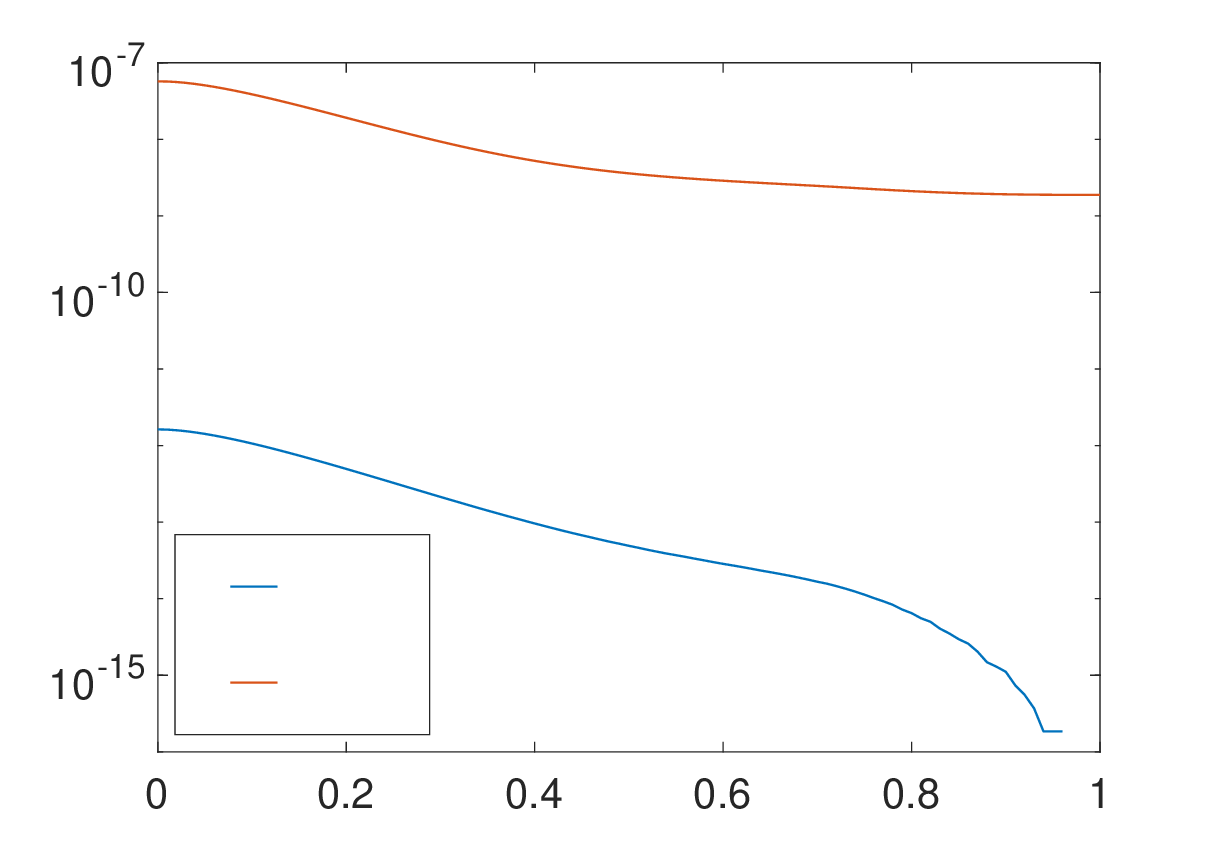}
    \put(-115,-5){$z$}
    \put(-175,47){$\delta \phi_\text{el}$}
    \put(-175,28){$\delta \phi_\text{ion}$}
    \put(-230,160){$\textbf{a)}$}
    \hspace{2mm}
    \includegraphics [scale=0.40]{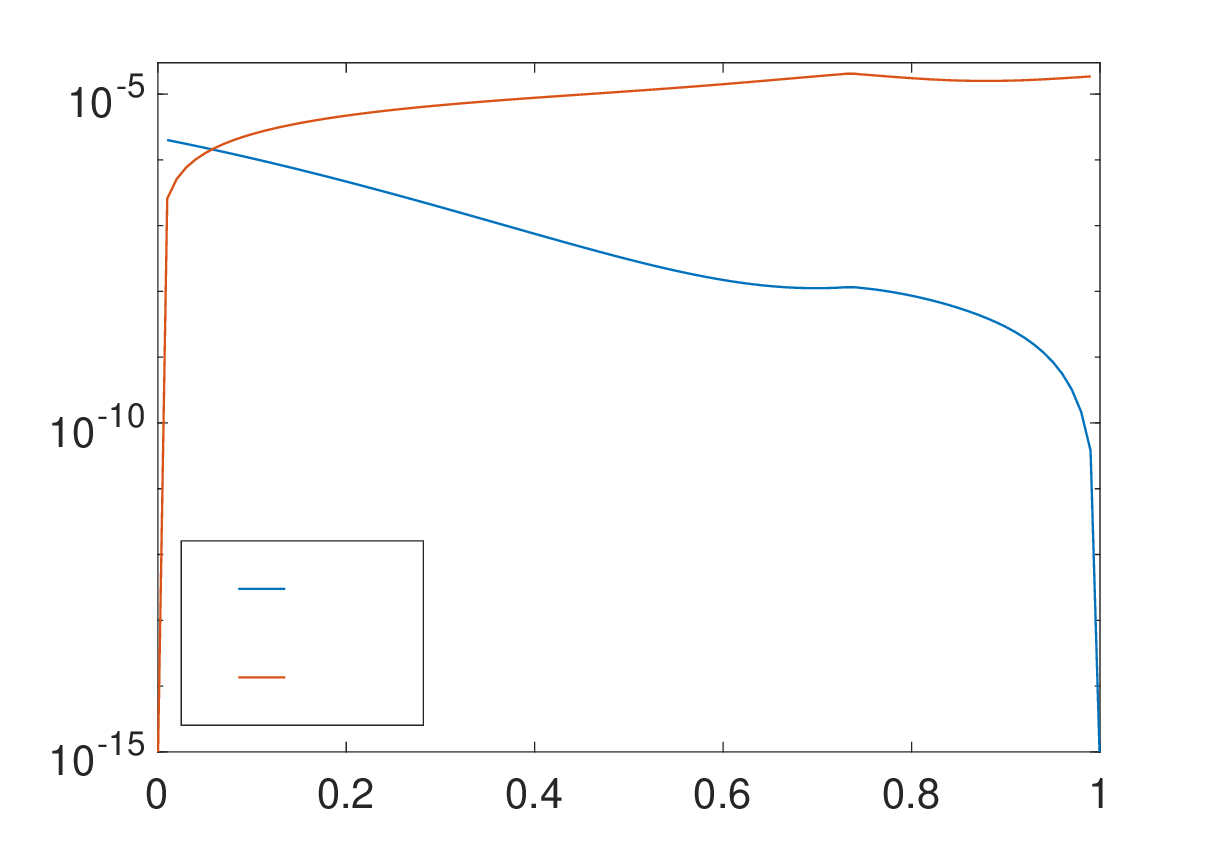}
     \put(-115,-5){$z$}
     \put(-175,47){$\delta j_\text{el}$}
    \put(-175,28){$\delta j_\text{ion}$}
    \put(-230,160){$\textbf{b)}$}

    \caption{The relative deviations between respective components of solution in terms of: a) the electron and ion potentials, b) the electron and ion currents.}

\label{delta_dfi_jedn}
\end{figure}

\begin{figure}
\begin{center}
\includegraphics[scale=0.40]{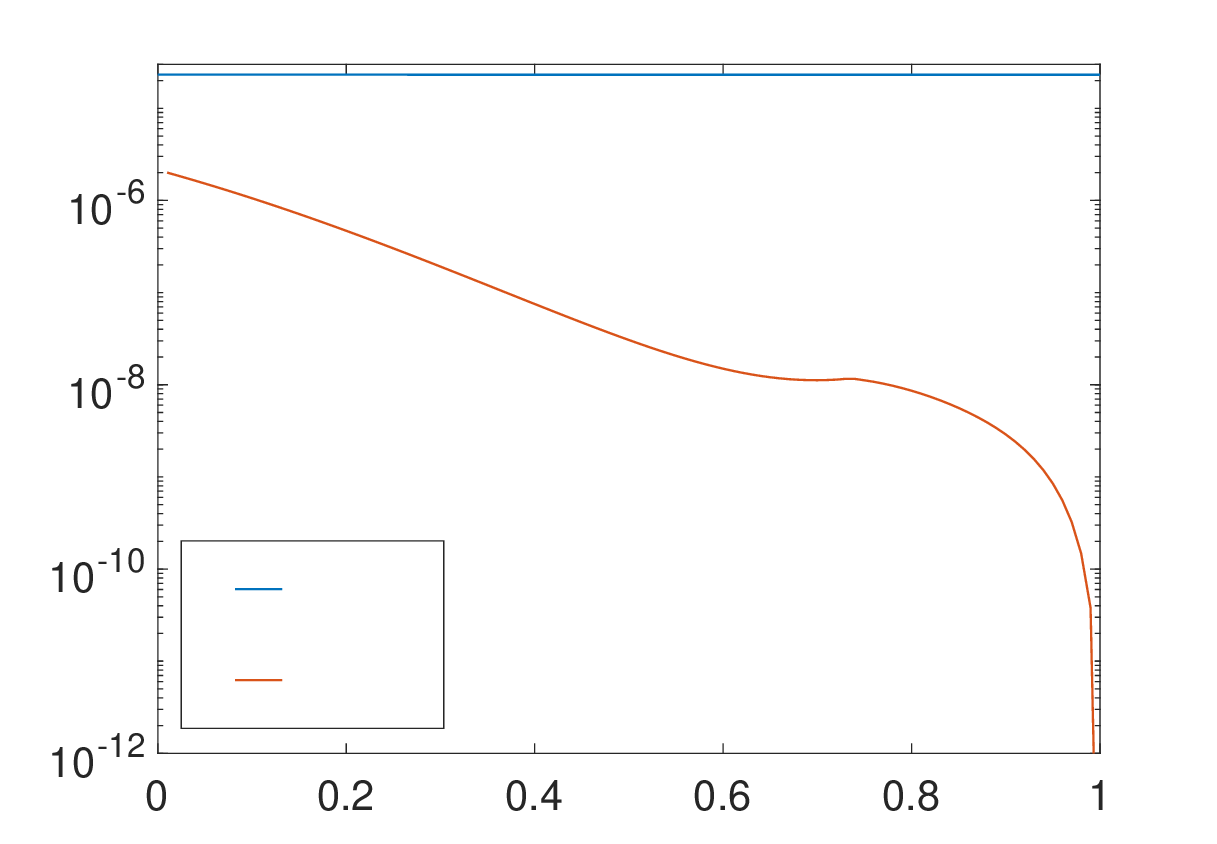}
\put(-120,-5){$z$}
\put(-175,47){$\delta C_{\text{O}_2}$}
\put(-175,28){$\delta J_{\text{O}_2}$}
\caption{The relative deviations between respective components of solution in terms of mass fraction of oxygen and mass flux of oxygen.}
\label{delta_co2_jedn}
\end{center}
\end{figure}

The comparison of solution obtained by means of the modified algorithm with that computed previously in subsection \ref{LSCF_num} is depicted in Figs. \ref{delta_dfi_jedn}--\ref{delta_co2_jedn}, where relative deviations between respective dependent variables are shown (the deviations are marked by symbol $\delta$ combined with corresponding dependent variable). As can be seen, a very good agreement between the analyzed solutions is obtained. The maximal deviations' order ranges from $10^{-12}$ (for the electron potential) to $10^{-5}$ (for the mass fraction of oxygen). The accuracy to which $C_{\text{O}_2}^\text{bulk}$ was recreated amounts to $2.24 \cdot 10^{-5}$, and indeed it is the mass fraction of oxygen which yields the biggest discrepancy between the analyzed results.  However, from the presented comparison it follows that both solutions are equivalent (identical up to the level of accuracy of computations). In this way, the conclusion of analysis delivered in subsection \ref{cat_sub}, that for the predefined parameters of the cell operation there exists a unique value of the active layer thickness, has been numerically confirmed.

\begin{figure}
\begin{center}
\includegraphics[scale=0.40]{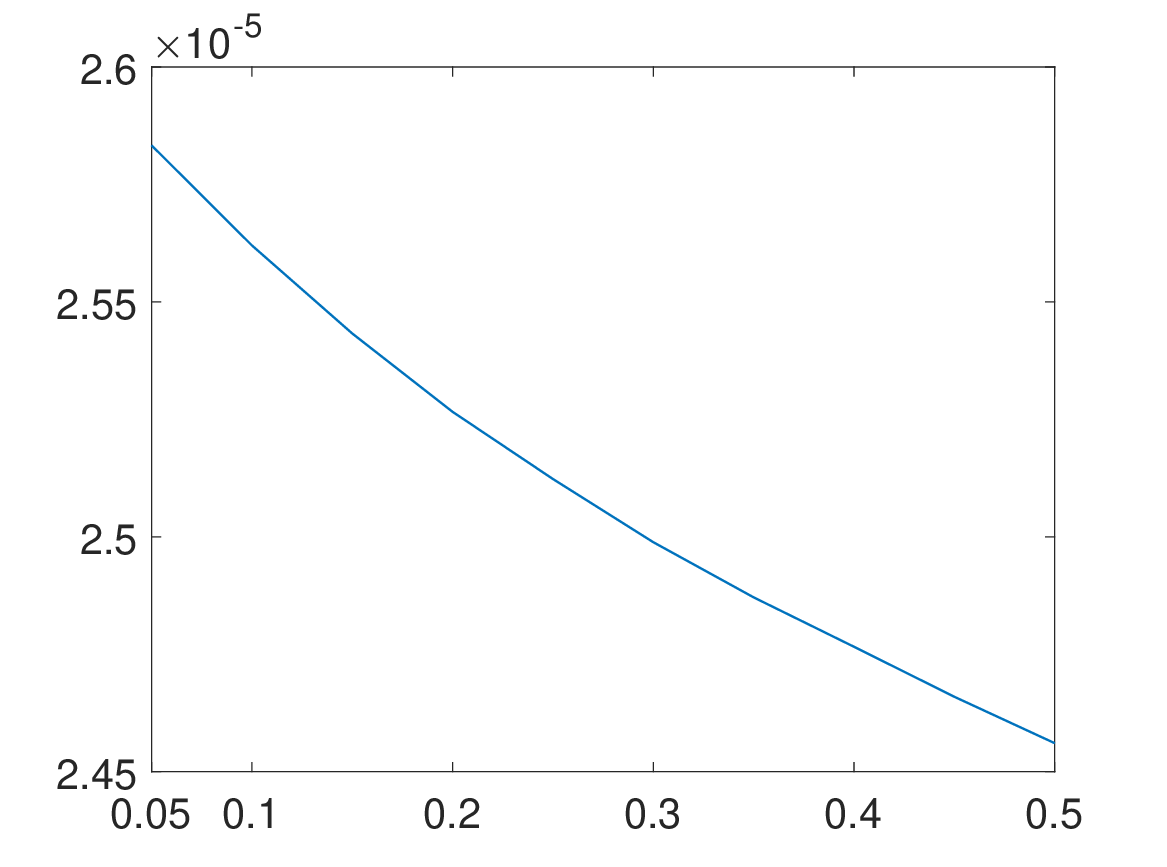}
\put(-120,-5){$x_{\text{O}_2}^\text{bulk}$}
\put(-230,80){$\delta_\text{c}$}
\caption{The relation between the oxygen volume fraction in the air, $x_{\text{O}_2}^\text{bulk}$, and the active layer thickness, $\delta_\text{c}$. }
\label{w_o2_delta}
\end{center}
\end{figure}

To complement this subsection we would like to discuss one aspect of the solution sensitivity. Namely, as it results from the provided analysis, one can implement the computational strategy where the active layer thickness, $\delta_\text{c}$, is predefined with some other cell operation parameter to be found as an element of solution (in our case it was $C_{\text{O}_2}^\text{bulk}$). However, even though different strategies are supposed to yield the same results, such an approach may not be computationally reasonable (or be in fact prohibitive). In Fig. \ref{w_o2_delta} we present a numerically obtained relation between the oxygen volume fraction in the air, $x_{\text{O}_2}^\text{bulk}$, and the resulting active layer thickness, $\delta_\text{c}$. The characteristic was delivered for $T=800^\circ \text{C}$ and $j_\text{cell}=2000$ $\text{A}\cdot \text{m}^{-2}$.

It shows that when increasing the oxygen volume fraction ten times (from $5\%$ to $50\%$), the active layer thickness changes by a mere $5\%$. Thus, any small deviation from the actual value of $x_{\text{O}_2}^\text{bulk}$ has a negligible impact on $\delta_\text{c}$. On the other hand, if one assumes that it is the active layer thickness to be predefined in computations, any small inaccuracy in defining its value leads to very large errors of solution and may cause computational instability. Thus, the latter strategy should be avoided.

\subsection{Modified formulation of the SOFC problem - comparison with experimental data}

We started this paper by recalling the classical formulation of the mathematical model of SOFC. Next, in a way of careful investigation of the desired behaviour of the component physical fields and the resulting mathematical structure of solution, we came to the conclusion that standard mathematical model produces locally unphysical results (negative activation overpotentials and charge transfer current). A simple remedy to these peculiarities was proposed in the form of modifications \eqref{eta_act_ano_abs}-\eqref{eta_act_cat_abs} or \eqref{i_ct_gen_abs}. Next, we proved that the thickness of active layer is a component of solution to be found in computations. All the aforementioned novelties modify substantially the classical SOFC model. In the following we will check whether such a modified formulation of the problem is conducive to credible modelling of the underlying physics, by comparing the numerical results with available experimental data.

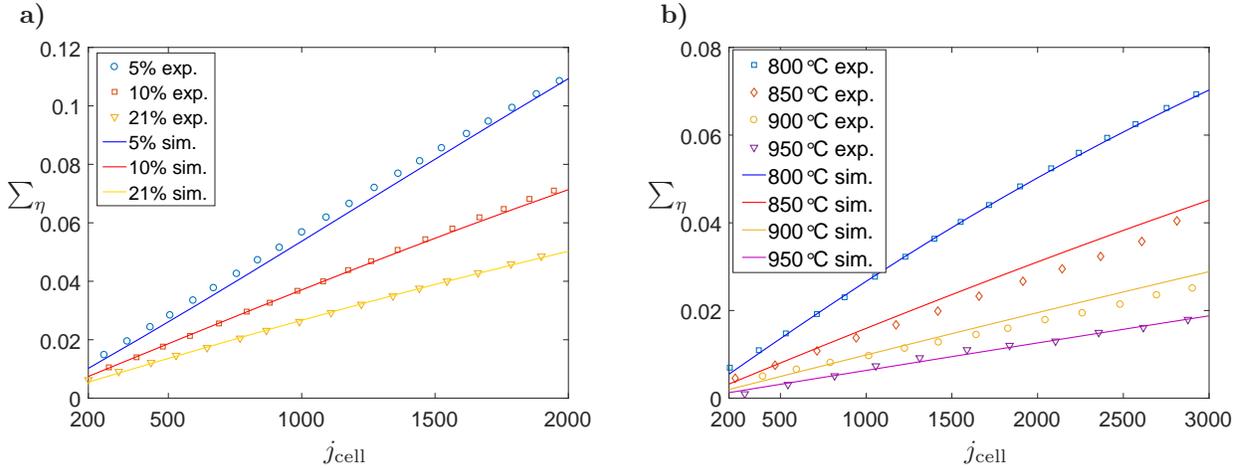
\begin{figure}[h!]

    \includegraphics [scale=0.40]{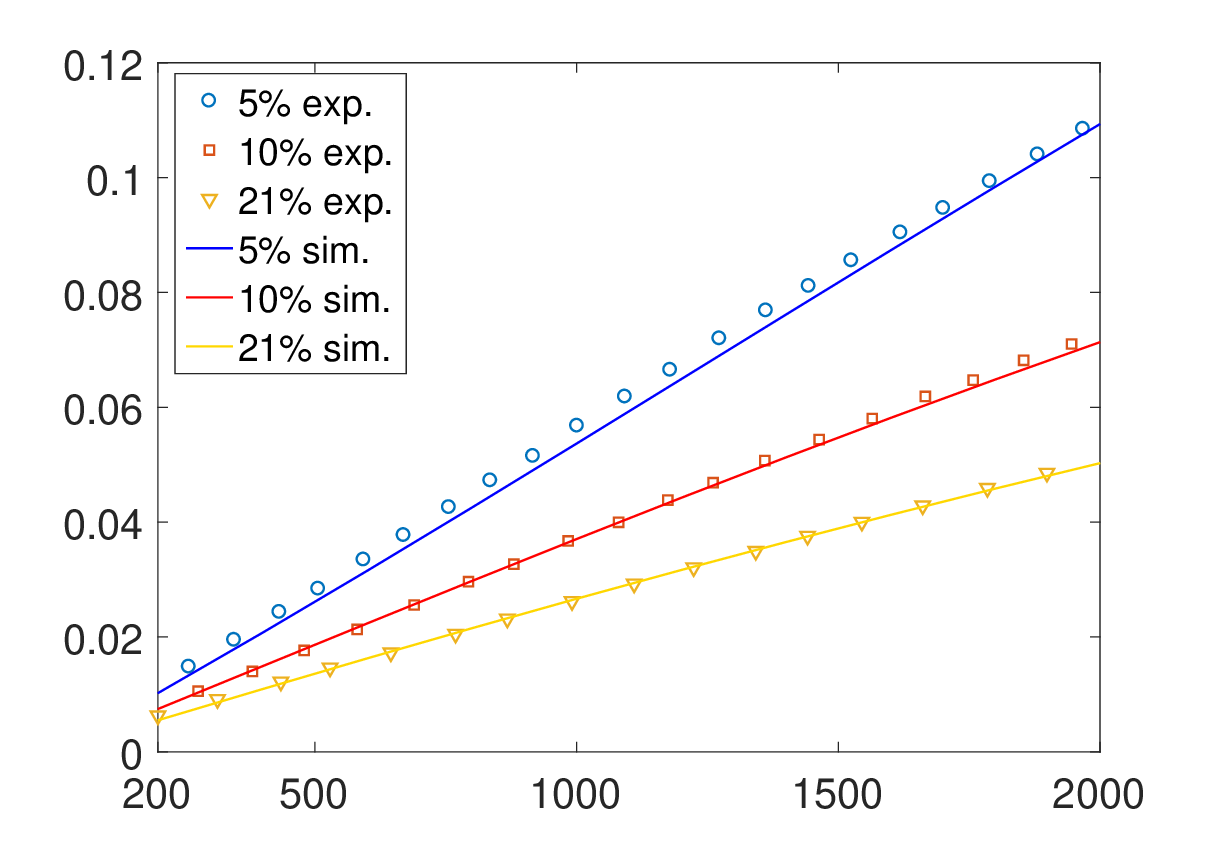}
    \put(-115,-5){$j_\text{cell}$}
		\put(-235,92){$\sum_\eta$}
    \put(-230,160){$\textbf{a)}$}
    \hspace{2mm}
    \includegraphics [scale=0.40]{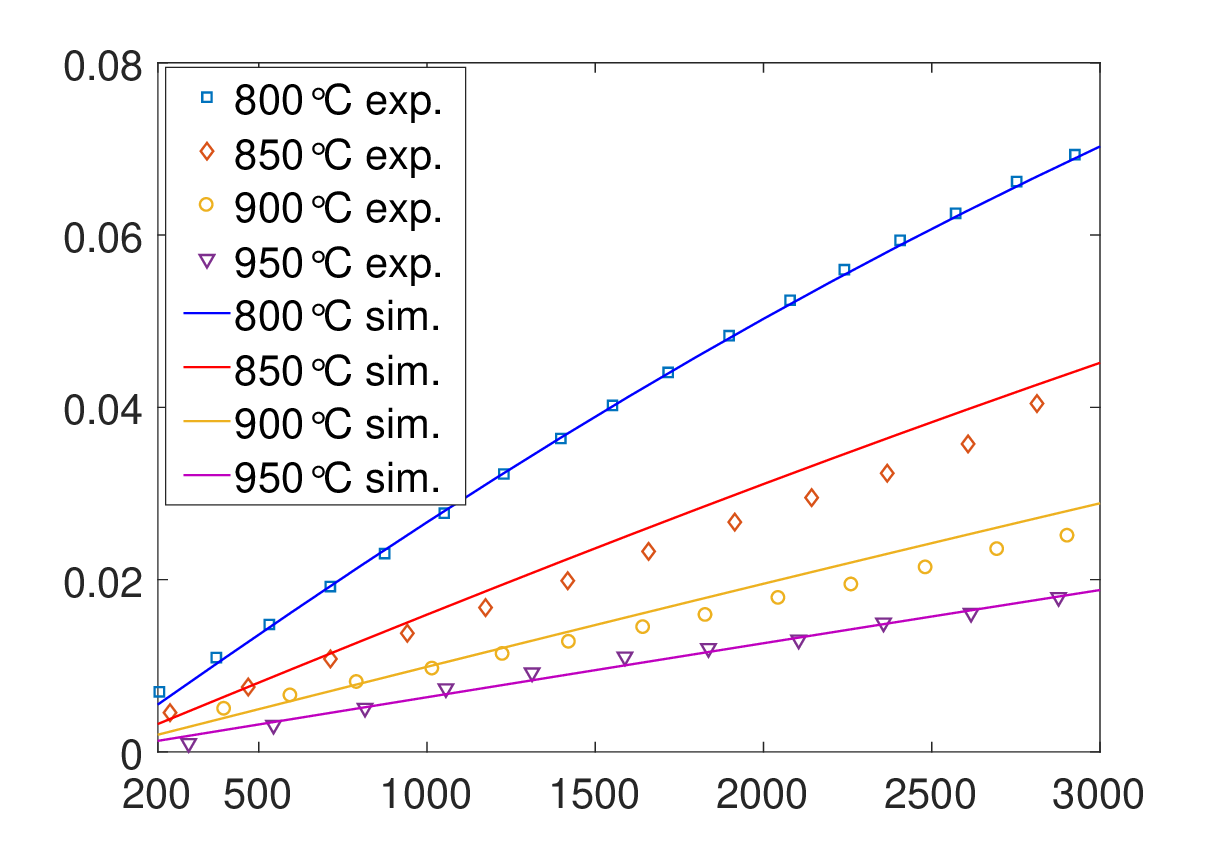}
    \put(-115,-5){$j_\text{cell}$}
		\put(-235,92){$\sum_\eta$}
    \put(-230,160){$\textbf{b)}$}

    \caption{Comparison of the experimental results from \cite{Miyoshi_2016} with numerical simulations in terms of the sum of activation and concentration overpotentials, $\sum _\eta$, for: a) different values of $x_{\text{O}_2}^\text{bulk}$ at $800^\circ \text{C}$, b) different values of temperature under $x_{\text{O}_2}^\text{bulk}=0.21$.}

\label{fig_Mio}
\end{figure}

In \cite{Miyoshi_2016} one can find experimental data for the LSM cathode of SOFC. The authors presented a measured sum of activation and concentration overpotentials, ${\sum}_\eta$, for different values of the cell current density, oxygen concentration and temperature.
Numerical simulations were performed with two different empirical relations  for the exchange current density, $i^\text{dpb}_\text{0,cat}$ (equations (6) and (7) therein for porous and patterned LSM, respectively).

Based on the information on the geometry, microstructural properties and operating conditions given in \cite{Miyoshi_2016} we carried out analogous numerical simulations and compared the obtained results  with the respective experimental data. The outcome of this comparison for the model of patterned LSM (which agreed better with the experiment) is shown in Fig.\ref{fig_Mio}. As can be seen, our numerical results mimic the experimental characteristics to a very large degree. It is only for $x_{\text{O}_2}^\text{bulk}=0.21$ and $T=850^\circ \text{C}$, $T=900^\circ \text{C}$ that the simulation deviates noticeably from the experiment. The level of coincidence between the computed and measured data is  better than that originally reported in \cite{Miyoshi_2016}. This allows us to conclude that the modifications introduced by us to the classical formulation of SOFC problem, facilitate the credible and efficient numerical simulation of the underlying physical phenomena.

\section{Discussion and final conclusions}
\label{disc}

In the paper an isothermal 1D mathematical model of the PEN structure of Solid Oxide Fuel Cell has been analyzed. The classical formulation has been carefully investigated in terms of the governing equations, the boundary and transmission conditions and the related asymptotic behaviour of solution. It has been proved that employing the standard definition of the activation overpotential in the Butler-Volmer equation leads to locally unphysical effects (negative charge transfer current and resulting peculiarities of the respective dependent variables). In order to remediate these issues a modification to the definition of the activation overpotetntial (or alternative definition of the charge transfer current) has been proposed. Next, it has been analytically substantiated that for the problem under consideration (even in its classical formulation) there is only one thickness of the active layer in each electrode for which the solution exists. As such, the active layers' thicknesses are the elements of solution to be found, unless some other solution parameters are relaxed. All the abovementioned discoveries have been included in the modified formulation of the SOFC problem. In the remaining part of the paper numerical computations for the cathode sub-problem have been carried out. To this end, a dedicated integral solver based on the modified formulation of the problem has been constructed. Its performance has been verified against a newly introduced analytical benchmark solution (for the combined physical fields). Next, computations for the LSCF cathode have been performed. On the basis of the obtained results, certain aspects of the modified formulation and the active layer thickness have been discussed. The uniqueness of solution with respect to the latter has been numerically confirmed. Finally, comparison with the experimental results have been done, which proved very good coincidence of the numerical solution for the modified formulation with the measured data.

Based on the conducted research, the following conclusions can be drawn:
\begin{itemize}
\item{When defining respective source terms for the electrochemical reactions in the governing equations, a clear distinction between the backing (inactive) and catalyst (active) layers of the electrodes should be imposed.   }
\item{The standard formulation of the problem leads by definition to locally unphysical results originating from negative values of the activation overpotentials (and thus negative values of the charge transfer current). }
\item{The aforementioned unphysical behaviour can be counteracted by employing a modified definition of the activation overpotential or alternatively a modified definition of the charge transfer current.   }
\item{The active layers' thicknesses in respective electrodes are supposed to be elements of solution (and thus cannot be taken in an arbitrary way), unless some other solution parameters are relaxed. Proper values of these thicknesses are necessary to guarantee the existence and uniqueness of solution.}
\item{Simultaneously, when employing a strict definition of the active layer accepted in this paper, its thickness not necessarily has to coincide with that of the real physical layer. The latter is predicted to be smaller. It can  possibly be found by applying the criterion of 99$\%$ of produced current.}
\item{The modified formulation of the problem facilitates credible simulation of the underlying physical phenomena, which has been confirmed by comparison with the experimental data. }
\item{The modified formulation opens a new avenue for modelling of the transient SOFC problems where the active layers thicknesses change with the operating conditions. The mathematical and computational scheme proposed in this paper gives measures for automatic identification of the respective active zones under dynamic operating regimes.}
\end{itemize}

\appendix

\section{Auxiliary transformations}
\label{apA}
\subsection{The cathode}

By combining equations \eqref{phi_el_ODE} and \eqref{co2_ODE} one obtains the following relation for the catalyst layer:
\begin{equation}
\label{cat_ax_1}
\frac{\text{d}}{\text{d}y}\left( \sigma_\text{el}\frac{\text{d} \phi_\text{el}}{\text{d}y}\right)=-\frac{4F}{M_{\text{O}_2}}\frac{\text{d}}{\text{d}y}\left(\rho_\text{a}D_\text{2}\frac{\text{d}C_{\text{O}_2}}{\text{d}y}\right).
\end{equation}
By integrating \eqref{cat_ax_1} from $h_1$ to $y$ under the  conditions \eqref{J_BCs}$_3$ and \eqref{j_el_ion_zero}$_1$ we have:
\begin{equation}
\label{cat_ax_2}
\sigma_\text{el}\frac{\text{d} \phi_\text{el}}{\text{d}y}=-\frac{4F\rho_\text{a}D_\text{2}}{M_{\text{O}_2}}\frac{\text{d}C_{\text{O}_2}}{\text{d}y}.
\end{equation}
Note that, according to \eqref{cat_ax_2}, the oxygen flux is proportional (with a minus sign) to the electron current.
After one more integration (this time from $y$ to $h_\text{b}$) and simple algebra, \eqref{cat_ax_2} allows to define the potential of the electron phase in the following form:
\begin{equation}
\label{fi_el_cat}
\phi_\text{el}(y)=V_2+\frac{j_\text{cell}}{\sigma_\text{el}}(h_2-h_\text{b})+\frac{4F\rho_\text{a}D_2}{M_{\text{O}_2}\sigma_\text{el}}(C_{\text{O}_2}^{\text{b}}-C_{\text{O}_2}(y)).
\end{equation}
In the above operations we accounted for the formula obtained by integrating \eqref{j_el_def}$_1$ with the boundary condition \eqref{j_el_BC}$_2$:
\begin{equation}
\label{cat_ax_3}
\phi_\text{el}(h_\text{b})=V_\text{b}=V_2+\frac{j_\text{cell}}{\sigma_\text{el}}(h_2-h_\text{b}).
\end{equation}

The value of the mass fraction of oxygen at $y=h_\text{b}$ can be derived by double integration of \eqref{co2_ODE} over the backing zone with the boundary conditions \eqref{C_BCs}$_3$ and \eqref{J_ext_BCs}$_3$:
\begin{equation}
\label{C^b_def}
C_{\text{O}_2}(h_\text{b})=C_{\text{O}_2}^{\text{b}}=C_{\text{O}_2}^{\text{bulk}}-\frac{M_{\text{O}_2}j_\text{cell}}{4F\rho_\text{a}D_2}(h_2-h_\text{b}).
\end{equation}
After substitution of \eqref{C^b_def} into \eqref{fi_el_cat} we arrive at the final formula for the potential of the electron phase:
\begin{equation}
\label{fi_el_cat_fin}
\phi_\text{el}(y)=V_2+\frac{4F\rho_\text{a}D_2}{M_{\text{O}_2}\sigma_\text{el}}(C_{\text{O}_2}^{\text{bulk}}-C_{\text{O}_2}(y)).
\end{equation}

In a similar way to that presented above, we can derive respective equation for the potential of ionic phase in the active zone. When comparing equations \eqref{phi_ion_ODE} and \eqref{co2_ODE} one has:
\begin{equation}
\label{cat_ax_4}
\frac{\text{d}}{\text{d}y} \left(\sigma_\text{ion}\frac{\text{d} \phi_\text{ion}}{\text{d}y}\right)=\frac{4F}{M_{\text{O}_2}}\frac{\text{d}}{\text{d}y}\left(\rho_\text{a}D_\text{2}\frac{\text{d}C_{\text{O}_2}}{\text{d}y}\right).
\end{equation}
Integration of \eqref{cat_ax_4} from $h_1$ to $y$ under the conditions \eqref{J_BCs}$_3$ and \eqref{j_el_ion_zero}$_2$ yields:
\begin{equation}
\label{cat_ax_5}
\sigma_\text{ion}\frac{\text{d} \phi_\text{ion}}{\text{d}y}+j_\text{cell}=\frac{4F}{M_{\text{O}_2}\rho_\text{a}D_\text{2}}\frac{\text{d}C_{\text{O}_2}}{\text{d}y}.
\end{equation}
By subsequent integration (from $y$ to $h_\text{b}$) we arrive at the expression for the ionic potential in the form:
\begin{equation}
\label{fi_ion_cat}
\phi_\text{ion}(y)=V_2+\frac{RT}{4F}\ln\left(\frac{C^\text{bulk}_{\text{O}_2}}{C^\text{b}_{\text{O}_2}}\right)-\frac{4F\rho_\text{a}D_2}{M_{\text{O}_2}\sigma_\text{ion}}(C_{\text{O}_2}^\text{bulk}-C_{\text{O}_2}(y))+\frac{j_\text{cell}}{\sigma_\text{ion}}(h_\text{b}-y)+j_\text{cell}(h_2-h_\text{b})\left(\sigma_\text{el}^{-1}+\sigma_\text{ion}^{-1}\right),
\end{equation}
where the value $\phi_\text{ion}(h_\text{b})$ was defined as:
\begin{equation}
\label{fi_ion_b}
\phi_\text{ion}(h_\text{b})=\phi_\text{b}=\phi_\text{el}(h_\text{b})+\eta^\text{conc}_\text{cat}(h_\text{b}),
\end{equation}
with the respective terms in the right hand side taken from \eqref{cat_ax_3} and \eqref{eta_conc_cat}.

\subsection{The anode}
By combining equations \eqref{phi_el_ODE} with \eqref{ch2_ODE} we obtain the following relation:
\begin{equation}
\label{an_ax_1}
\frac{\text{d}}{\text{d}y}\left(\sigma_\text{el}\frac{\text{d}\phi_\text{el}}{\text{d}y}\right)=\frac{2F}{M_{\text{H}_2}}\frac{\text{d}}{\text{d}y}\left(\rho_\text{f}D_1\frac{\text{d}C_{\text{H}_2}}{\text{d}y}\right).
\end{equation}
After integration from $y$ to $-h_1$ (under the conditions \eqref{J_ext_BCs}$_1$ and \eqref{j_el_ion_zero}$_1$) equation \eqref{an_ax_1} yields:
\begin{equation}
\label{an_ax_2}
\sigma_\text{el}\frac{\text{d}\phi_\text{el}}{\text{d}y}=\frac{2F\rho_\text{f}D_1}{M_{\text{H}_2}}\frac{\text{d}C_{\text{H}_2}}{\text{d}y}.
\end{equation}
Subsequent integration of \eqref{an_ax_2}, this time from $-h_\text{a}$ to $y$, leads to the formula:
\begin{equation}
\label{an_ax_3}
\phi_\text{el}=\phi_\text{el}(-h_\text{a})+\frac{2F\rho_\text{f}D_1}{M_{\text{H}_2}\sigma_\text{el}}\left(C_{\text{H}_2}-C_{\text{H}_2}^\text{a}\right),
\end{equation}
where the limiting value of the electron potential, $\phi_\text{el}(-h_\text{a})$, is obtained by integration \eqref{j_el_def}$_1$ under the condition \eqref{fi_el_BC}$_1$:
\begin{equation}
\label{V_a_def}
\phi_\text{el}(-h_\text{a})=V_\text{a}=V_3-j_\text{cell}\frac{h_3-h_\text{a}}{\sigma_\text{el}},
\end{equation}
while $C_{\text{H}_2}^\text{a}$ can be derived when integrating \eqref{ch2_ODE} with account for the boundary condition \eqref{C_BCs}$_1$ and \eqref{J_ext_BCs}$_1$:
\begin{equation}
\label{C_h2_a_def}
C_{\text{H}_2}^\text{a}=C_{\text{H}_2}(-h_\text{a})=C_{\text{H}_2}^\text{bulk}-\frac{M_{\text{H}_2}j_\text{cell}}{2F\rho_\text{f}D_1}(h_3-h_\text{a}).
\end{equation}
After substitution of \eqref{V_a_def} and \eqref{C_h2_a_def} into \eqref{an_ax_3} we obtain a final relation for the potential of the electron phase:
\begin{equation}
\label{fi_el_an_fin}
\phi_\text{el}(y)=V_3+\frac{2F\rho_\text{f}D_1}{M_{\text{H}_2}\sigma_\text{el}}\left(C_{\text{H}_2}(y)-C_{\text{H}_2}^\text{bulk}\right).
\end{equation}

In order to derive respective expression for the ionic potential we start by merging equations \eqref{phi_ion_ODE} and \eqref{ch2_ODE}:
\begin{equation}
\label{an_ax_4}
\frac{\text{d}}{\text{d}y}\left(\sigma_\text{ion}\frac{\text{d}\phi_\text{ion}}{\text{d}y}\right)=-\frac{2F}{M_{\text{H}_2}}\frac{\text{d}}{\text{d}y}\left(\rho_\text{f}D_1\frac{\text{d}C_{\text{H}_2}}{\text{d}y}\right).
\end{equation}
By integrating \eqref{an_ax_4} from $y$ to $-h_1$ with conditions \eqref{J_BCs}$_1$ and \eqref{j_el_ion_zero}$_2$ we have:
\begin{equation}
\label{an_ax_5}
\sigma_\text{ion}\frac{\text{d}\phi_\text{ion}}{\text{d}y}=-j_\text{cell}-\frac{2F\rho_\text{f}D_1}{M_{\text{H}_2}}\frac{\text{d}C_{\text{H}_2}}{\text{d}y}.
\end{equation}
Consecutive integration (from $-h_\text{a}$ to $y$) and some algebra yield:
\begin{equation}
\label{an_ax_6}
\phi_\text{ion}=\phi_\text{ion}(-h_\text{a})-\frac{j_\text{cell}}{\sigma_\text{ion}}(h_\text{a}+y)-\frac{2F\rho_\text{f}D_1}{M_{\text{H}_2}\sigma_\text{ion}}\left(C_{\text{H}_2}-C_{\text{H}_2}^\text{a}\right),
\end{equation}
where $\phi_\text{ion}(-h_\text{a})$ is computed as:
\begin{equation}
\label{fi_a_def}
\phi_\text{ion}(-h_\text{a})=\phi_\text{a}=\phi_\text{el}(-h_\text{a})-\eta_\text{conc}^\text{ano}(-h_\text{a}),
\end{equation}
with respective entries taken according to \eqref{V_a_def} and \eqref{eta_conc_ano}, while $C_{\text{H}_2}^\text{a}$ is being defined by \eqref{C_h2_a_def}. When accounting for the above formulae, \eqref{an_ax_6} assumes the following form:
\begin{equation}
\label{fi_ion_an_fin}
\phi_\text{ion}(y)=V_3-\frac{RT}{2F}\ln\left(\frac{C_{\text{H}_2}^\text{bulk}}{C_{\text{H}_2}^\text{a}}\frac{C_{\text{H}_2\text{O}}^\text{a}}{C_{\text{H}_2\text{O}}^\text{bulk}}\right)
-\frac{j_\text{cell}}{\sigma_\text{ion}}(h_\text{a}+y)-
\frac{2F\rho_\text{f}D_1}{M_{\text{H}_2}\sigma_\text{ion}}\left(C_{\text{H}_2}-C_{\text{H}_2}^\text{bulk}\right)-j_\text{cell}(h_3-h_a)\left(\sigma_\text{el}^{-1}+\sigma_\text{ion}^{-1}\right).
\end{equation}

\section{The benchmark example}
\label{benchmark_app}

In order to verify the accuracy of computations provided by the proposed numerical scheme we shall introduce a dedicated benchmark solution. It should preserve all the essential features of the original problem, including the type of boundary and transmission conditions, as well as the asymptotic behaviour of solution. Our strategy to built such a benchmark relies on: i) defining a desired form of the reference solution, ii) modifying the system of governing ODEs, by addition of a known predefined function, in a way that the resulting equations are identically satisfied for the selected solution. Note that the latter does not detract from the applicability and relevance of the results, as it simply amounts to adding additional (known) source terms in the right hand sides of respective equations. Obviously, for the 'real' problems, those terms turn to zero.

In the following we will use the spatial rescaling \eqref{z_def}. Only the formulae for the catalyst layer will be given (i.e. for $y\in(h_1,h_\text{b})$), as respective relations for the backing zone of the cathode can be taken from the subsection \ref{comp_rels}.  Let the benchmark solution for the ionic current be of the form:
\begin{equation}
\label{j_ion_bench}
j_\text{ion}(z)=j_\text{cell}(1-z)^\alpha,
\end{equation}
where $\alpha>1$ is a parameter to be taken as convenient. Consequently, the electron current is taken in accordance with \eqref{j_cell} to yield:
\begin{equation}
\label{j_el_bench}
j_\text{el}(z)=
j_\text{cell}[1-(1-z)^\alpha].
\end{equation}
The above representation complies with the desired qualitative behavior of the current functions and respective transmission conditions.

The expressions for the ion and electron potentials can be obtained by integration of \eqref{j_ion_bench}-\eqref{j_el_bench}:
\begin{equation}
\label{fi_ion_bench}
\phi_\text{ion}(z)=\phi_\text{b}+\frac{\delta_\text{c}j_\text{cell}}{\sigma_\text{ion}(1+\alpha)}(1-z)^{1+\alpha},
\end{equation}
\begin{equation}
\label{fi_el_bench}
\phi_\text{el}(z)=V_\text{b}+\frac{\delta_\text{c}j_\text{cell}}{\sigma_\text{el}}\left[1-\frac{1}{1+\alpha}(1-z)^{1+\alpha}\right].
\end{equation}
Again, the desired qualitative behavior of dependent variables and transmission conditions are satisfied.

The benchmark function for the oxygen flux is assumed as follows:
\begin{equation}
\label{J_o2_bench}
J_{\text{O}_2}=-\frac{M_{\text{O}_2}j_\text{cell}}{4F}[1-(1-z)^\alpha].
\end{equation}
When integrated with respect to the spatial variable, equation \eqref{J_o2_bench} yields:
\begin{equation}
\label{C_o2_bench}
C_{\text{O}_2}(z)=C_{\text{O}_2}^\text{b}-\frac{M_{\text{O}_2}j_\text{cell}\delta_\text{c}}{4F\rho_\text{a}D_2}(1-z)\left[1-\frac{1}{1+\alpha}(1-z)^\alpha\right].
\end{equation}
Note that also here the proper qualitative bahviour and transmission conditions for respective functions hold automatically.

Having the above components of solution we can compute the activation overpotential and charge transfer current according to \eqref{eta_act_cat_abs} and \eqref{i_cat_def}. The active layer thickness, $\delta_\text{c}$, can be found iteratively by employing equation \eqref{delta_c_obl}.

Unfortunately, in general the introduced functions for dependent variables do not satisfy the governing ODEs (see \eqref{phi_el_ODE},\eqref{phi_ion_ODE}, \eqref{co2_ODE}) in the active layer. That is why, for the sake of benchmark computations, we will modify slightly respective equations by adding an additional known source term. The governing equations (written for simplicity in terms of the ion and electron currents and the oxygen flux) have now the following forms in the active layer:
\begin{equation}
\label{ODE_ion_bench}
\frac{\text{d}j_\text{ion}}{\text{d}y}=i^\text{ct}_\text{cat}+g(y),
\end{equation}
\begin{equation}
\label{ODE_el_bench}
\frac{\text{d}j_\text{el}}{\text{d}y}=-i^\text{ct}_\text{cat}-g(y),
\end{equation}
\begin{equation}
\label{ODE_Jo2_bench}
\frac{\text{d}J_{\text{O}_2}}{\text{d}y}=-\frac{M_{\text{O}_2}}{4F}\left(i^\text{ct}_\text{cat}+g(y)\right),
\end{equation}
where $g(y)$ is a known predefined function:
\begin{equation}
\label{g_def}
g(y)=i^\text{ct(b)}_\text{cat}+j_\text{cell}\frac{\alpha}{\delta_\text{c}}\left(\frac{\delta_\text{c}+h_1-y}{\delta_\text{c}}\right)^{\alpha-1}.
\end{equation}
$i^\text{ct(b)}_\text{cat}$ is computed by substitution of \eqref{fi_ion_bench},\eqref{fi_el_bench} and \eqref{C_o2_bench} into \eqref{i_cat_def}.

After the above modifications, the proposed benchmark solution satisfies the governing system of ODEs together with respective boundary and transmission conditions. Thus, it constitutes a very useful reference solution to investigate the accuracy and efficiency of computations by the numerical scheme employed (or any other numerical scheme).




\vspace{3mm}
\noindent
{\bf Acknowledgements:} The presented research is the part of the Easy-to-Assemble Stack Type (EAST): Development of solid oxide fuel cell stack for the innovation in Polish energy sector project, carried out within the FIRST TEAM program (project number First TEAM/2016-1/3) of the Foundation for Polish Science, co-financed by the European Union under the European Regional Development Fund. The authors are grateful for the support.

\end{document}